\documentclass[12pt,a4paper,titlepage,twoside]{book}
\usepackage[italian]{babel}
\usepackage{amsfonts,graphics,graphicx,hyperref,bm,url,color}
\newtheorem{teo}{Teorema}[chapter]
\newtheorem{cor}{Corollario}[chapter]
\newtheorem{lem}{Lemma}[chapter]
\newtheorem{oss}{Osservazione}[chapter]
\newtheorem{defi}{Definizione}[chapter]
\newtheorem{eser}{Esercizio}[chapter]
\newtheorem{algo}{Algoritmo}[chapter]
\def\QED{\mbox{\ \vbox{\hrule\hbox{\vrule height1.3ex\hskip0.8ex\vrule}\hrule}}}
\def\no{\noindent}
\def\nulit{\item[~]}
\def\RR{{\rm I\!R}}
\def\NN{{\rm I\!N}}
\def\CC{{\mathchoice {\setbox0=\hbox{$\displaystyle\rm C$}\hbox{\hbox
to0pt{\kern0.4\wd0\vrule height0.9\ht0\hss}\box0}}
{\setbox0=\hbox{$\textstyle\rm C$}\hbox{\hbox
to0pt{\kern0.4\wd0\vrule height0.9\ht0\hss}\box0}}
{\setbox0=\hbox{$\scriptstyle\rm C$}\hbox{\hbox
to0pt{\kern0.4\wd0\vrule height0.9\ht0\hss}\box0}}
{\setbox0=\hbox{$\scriptscriptstyle\rm C$}\hbox{\hbox
to0pt{\kern0.4\wd0\vrule height0.9\ht0\hss}\box0}}}}
\def\pmatrix{ \left( \begin{array} }
\def\endpmatrix{ \end{array} \right) }
\def\proof{{\em \underline{Dim.}}\quad}
\def\rank{\mbox{\rm rank}}
\def\aa{\alpha}
\def\bb{\beta}
\def\cc{\gamma}
\def\dd{\delta}
\def\om{\omega}
\def\lam{\lambda}
\def\bfa{{\bm{a}}}
\def\bfb{{\bm{b}}}
\def\bfc{{\bm{c}}}

\def\bfe{{\bm{e}}}

\def\bfg{{\bm{g}}}
\def\bfl{{\bm{l}}}
\def\bfel{{\bm{\ell}}}
\def\bfcc{{\bm{\gamma}}}
\def\bfp{{\bm{p}}}
\def\bfq{{\bm{q}}}
\def\bfr{{\bm{r}}}
\def\bfs{{\bm{s}}}
\def\bft{{\bm{t}}}
\def\bfu{{\bm{u}}}
\def\bfv{{\bm{v}}}
\def\bfw{{\bm{w}}}
\def\bfx{{\bm{x}}}
\def\bfy{{\bm{y}}}
\def\bfz{{\bm{z}}}
\def\bfo{{\bm{0}}}
\def\eps{\varepsilon}
\def\phi{\varphi}
\def\hT{\hat{T}}
\def\bz{\bar{z}}
\def\hv{\hat{v}}
\def\tL{\tilde{L}}
\def\tD{\tilde{D}}
\def\td{\tilde{d}}
\def\tG{\tilde{G}}

\begin{document}

\pagenumbering{roman}

\begin{titlepage}

\begin{center}
\rule{6cm}{.125mm}\\
\rule{8cm}{.25mm}\\
\rule{10cm}{.5mm}\\
\rule{12cm}{1mm} \\[1cm] {\LARGE\bf
\begin{tabular}{l}
METODI\\METODI I\\METODI IT\\METODI ITE\\
METODI ITER\\METODI ITERA\\METODI ITERAT\\METODI ITERATI\\METODI
ITERATIV\\ METODI ITERATIVI\end{tabular}\\[3mm]per\\[3mm] SISTEMI
LINEARI}\\[13mm] {\large\bf Luigi Brugnano \qquad Cecilia Magherini}\\[13mm]
\rule{12cm}{1mm}\\
\rule{10cm}{.5mm}\\
\rule{8cm}{.25mm}\\
\rule{6cm}{.125mm}
\end{center}

\end{titlepage}

\thispagestyle{empty}
\mbox{~}
\begin{flushright}\em In memoria di Cecilia\end{flushright}

\vskip 2cm
\parbox{7cm}{\em Tutti gli avvenimenti sono concatenati nel
migliore dei mondi possibili: \dots\\[1mm] -- \`E giusto, rispose Candido,
ma bisogna coltivare il no\-stro giardino. \\[3mm] \rm Voltaire, \em
Candido.}

\setcounter{page}{2}
\tableofcontents
\chapter*{Prefazione}
\addcontentsline{toc}{chapter}{Prefazione}
\markboth{PREFAZIONE}{~}

Queste note intendono dare le nozioni di base riguardo ai metodi
iterativi per la risoluzione di sistemi di equazioni lineari.

Il materiale presentato \`e sufficiente per un modulo di 40-45 ore
di lezione. Inoltre, gli esercizi proposti possono essere
argomento di ulteriori ore di esercitazione.

La trattazione, sebbene non esaustiva, ha comunque l'obiettivo di
fornire, in un contesto il pi\`u possibile omogeneo ed
autoconsistente, le idee principali che hanno ispirato le tecniche
di risoluzione pi\`u note. Di queste ultime si daranno molti
dettagli, sebbene l'aspetto implementativo non sar\`a, per
brevit\`a, approfondito.

Anche la bibliografia \`e ispirata da analoghi criteri: sono
infatti riportate solo le fonti bibliografiche essenziali, dove i
princ\`\i pi che hanno portato alla derivazione delle principali
metodologie sono descritti. Nel fare questo, tuttavia, non vi sono
pretese di esaustivit\`a. A questo riguardo, si sottolinea che
alcuni dei testi riportati in bibliografia contengono riferimenti
bibliografici molto pi\`u approfonditi.

Le note sono suddivise in 7 capitoli ed una appendice:

\begin{enumerate}
\setlength{\itemsep}{0cm}

\item un breve capitolo introduttivo;

\item un capitolo contenente le nozioni di base;

\item un capitolo che illustra i metodi iterativi di base;

\item un capitolo dedicato al caso di sistemi con matrici
simmetriche e definite positive;

\item un capitolo riguardante il precondizionamento
di matrici simmetriche e definite positive;

\item un capitolo dedicato al caso di sistemi con matrici solo
simmetriche;

\item un capitolo finale dedicato al caso di matrici nonsingolari
generiche;

\item[A.] una appendice riguardante la memorizzazione di matrici
sparse.

\end{enumerate}

Gli esempi riportati sono stati codificati in Matlab.

\subsection*{Ringraziamenti} Si ringraziano gli studenti del corso di
Analisi Numerica della scuola estiva SMI 2002, Perugia, e quelli
del corso di Analisi Numerica per Informatica, a.a. 2001/02,
presso l'Universit\`a degli Studi di Firenze, per le loro
osservazioni e gli utili commenti.

\subsection*{Edizione corretta}
Questa \`e una versione delle note pubblicate da Pitagora Editrice Bologna nel 2003 (ISBN 88-371-1367-6) 
in cui si fissano alcuni errori nel testo, oltre ad apportare delle piccole modifiche.

\chapter*{Notazioni}
\addcontentsline{toc}{chapter}{Notazioni}
\markboth{NOTAZIONI}{~}

\begin{itemize}

\item $\NN$ denota l'insieme dei numeri naturali;

\item $\RR$ ($\CC$) denota il campo reale (complesso);

\item $\RR^+, \RR^-$ indicano la semiretta reale positiva e negativa, rispettivamente;

\item $\RR^m$ \`e lo spazio vettoriale dei vettori reali di dimensione $m$;

\item $\RR^{m\times n}$ \`e lo spazio vettoriale delle matrici
reali $m\times n$;

\item $\nabla f$ denota il gradiente di una funzione $f$;

\item $\bfv^\top$ denota il trasposto del vettore $\bfv$, $\bfv^*$ ne denota
il trasposto e coniu\-gato; simile notazione \`e utilizzata per le
matrici;

\item $\sigma(A)$ denota lo spettro della matrice quadrata $A$;

\item $\rho(A)$ denota il raggio spettrale di $A$, $\max_{\lam\in\sigma(A)}|\lam|$;

\item $\rank(A)$ denota il rango della matrice $A$;

\item $A\sim B$ denota che le due matrici (quadrate) $A$ e $B$
sono simili;

\item $\Pi_k$ denota lo spazio vettoriale dei polinomi a coefficienti reali di grado
al pi\`u $k$;

\item $\Pi_k'$ denota l'\,insieme dei polinomi monici a coefficienti reali di
grado $k$;

\item $A \otimes B$ denota il prodotto di Kronecker, o tensoriale,
delle matrici $A$ e $B$. Se $A=(a_{ij})$, si ottiene la matrice a
blocchi $A\otimes B=(a_{ij}B)$;

\item $I$ denota la matrice identit\`a, le cui dimensioni sono
deducibili dal contesto;

\item $I_N$ denota la matrice identit\`a di dimensione $N\times
N$;

\item $O$ denota la matrice nulla di dimensione opportuna
deducibile dal contesto;

\item $\bfo$ denota il vettore (colonna) nullo, di dimensione
opportuna deducibile dal contesto;

\item {\em sdp} \`e l'acronimo di (matrice) {\em simmetrica e definita
positiva};

\item $\|\cdot\|$ indica la norma su vettore o, a seconda del contesto,
quella indotta su matrice. Quando non altrimenti specificato, la
norma \`e la norma 2;

\item $\kappa(A)$ indica il numero di condizione della
matrice $A$, $\|A\|\cdot\|A^{-1}\|$. Se non altrimenti
specificato, la norma \`e la norma 2;

\item $E_i$ indica l'$i$-esimo versore in $\RR^n$, dove la
dimensione $n$ dello spazio \`e deducibile dal contesto;

\item $E_i^{(k)}$ indica l'$i$-esimo versore in $\RR^k$;

\item $[\bfu_1,\dots,\bfu_k]$ denota il sottospazio generato dai
vettori racchiusi tra pa\-rentesi. In altre parole, il {\em
range}, ${\rm ran}(U_k)$, della matrice $U_k$, avente per colonne
i vettori $\{\bfu_i\}$;

\item $\leftarrow$ indica, nella descrizione degli algoritmi, la
riscrittura del membro sinistro con il risultato del membro
destro.

\item {\em matvec, scal, axpy}: denotano, rispettivamente, il
{\em prodotto matrice-vettore}, il {\em prodotto scalare} tra vettori, e
l'operazione\\ \centerline{``{\em vettore $\leftarrow$
scalare\,$\times$\,vettore\,$+$\,vettore}\,''.}

\end{itemize}

~

\newpage

\setcounter{page}{0}
\pagenumbering{arabic}
%
%
\chapter{Introduzione e motivazioni}\label{cap1}

{\em Ci occuperemo della risoluzione di sistemi
lineari ``sparsi''. In questo breve capitolo introduttivo introduciamo 
il concetto di matrice sparsa, dando una motivazione pratica per la
trattazione successiva.}

\section{Matrici sparse}
La risoluzione di molti problemi applicativi comporta la
risoluzione di sistemi di equazioni lineari,
\begin{equation}\label{Axb}
A\bfx = \bfb,
\end{equation}

\no in cui $\bfx,\bfb\in\RR^n$, e $A\in\RR^{n\times n}$ \`e una matrice
nonsingolare, avente la caratte\-ristica aggiuntiva di essere {\em
sparsa}. Queste caratteristiche di $A$ saranno sempre assunte nel
seguito.

\begin{defi} Diremo che una matrice $n\times n$ \`e {\em sparsa}
quando il numero dei suoi elementi non nulli \`e proporzionale ad
$n$, invece che a $n^2$.\end{defi}

Per matrici sparse, esistono efficienti tecniche di memorizzazione
che consentono di memorizzare, essenzialmente, solo gli elementi
significativi della matrice (o poco pi\`u di questi). In tal modo,
la occupazione di memoria, all'interno di un calcolatore, pu\`o
essere notevolmente ridotta, rispetto alle cano\-niche $n^2$ celle
di memoria richieste da una matrice piena di pari dimensioni. In
tal caso, i metodi diretti di risoluzione, basati sul calcolo di
opportune fattorizzazioni della matrice $A$, risultano essere in
generale non convenienti, sia per quanto riguarda l'occupazione di
memoria, che in termini di operazioni da effettuare.

Infatti, i fattori di una matrice sparsa sono difficilmente
matrici sparse a loro volta. Pertanto, salvo casi particolari,
essi necessitano di una occupazione di memoria molto pi\`u
elevata, rispetto alla matrice originaria.

Inoltre, per quanto riguarda il costo computazionale, i metodi di
fattorizzazione diretta richiedono generalmente un numero di
operazioni proporzionale al cubo della dimensione di $A$, con un
costo quindi molto elevato, anche per valori di $n$ non
eccessivamente grandi. D'altro canto, osserviamo che se $A$
possiede $O(n)$ elementi significativi, allora il prodotto di $A$
per un vettore pu\`o essere convenientemente formulato in modo da
richiedere $O(n)$ operazioni.

\begin{figure}[t]
\begin{center}
\includegraphics[width=13cm,height=10cm]{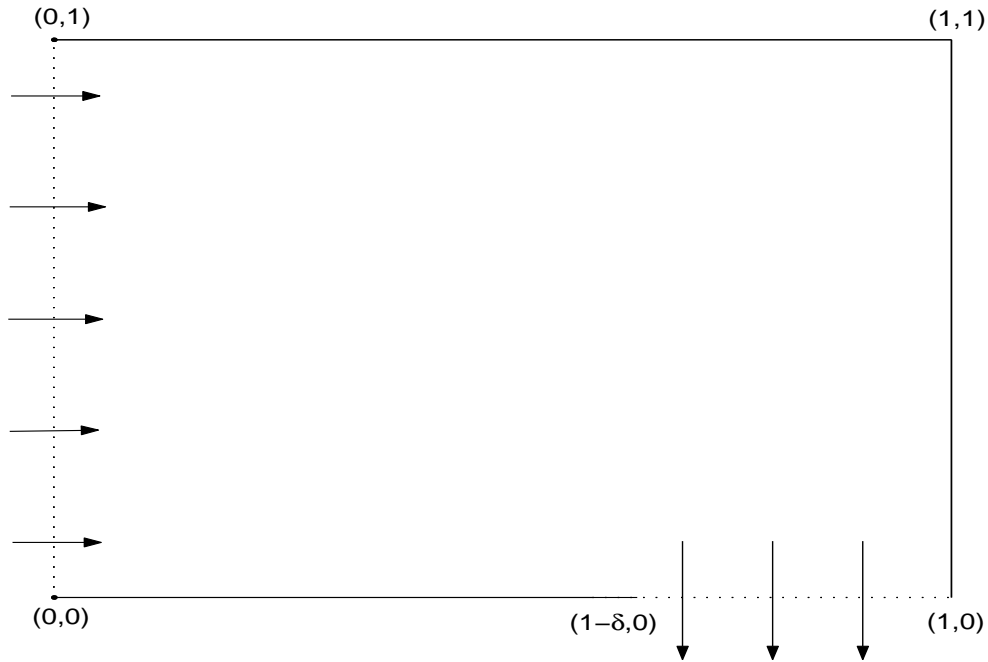}
\caption{\protect\label{fig0} Flusso di un fluido in un mezzo
poroso.}
\end{center}
\end{figure}

\section{Un esempio}

Si consideri il flusso di un fluido incomprimibile attraverso una
cavit\`a, contenente un mezzo poroso omogeneo, disposta
parallelamente al livello di riferimento della gravit\`a. Si
considerer\`a una sezione orizzontale della cavit\`a. Supporremo
inoltre che il fluido abbia raggiunto un regime di flusso
stazionario e che non vi siano, all'interno della cavit\`a,
sorgenti o pozzi. In tal caso, supponendo che il fluido entri
nella cavit\`a da sinistra e fuoriesca in basso da destra (vedi
Figura~\ref{fig0}), le equazioni che descrivono il flusso sono
date, dopo opportune normalizzazioni, da:

\begin{eqnarray} \nonumber
u(x,y) &=& -\frac{\partial p}{\partial x}(x,y), \\[1mm] \nonumber
v(x,y) &=&-\frac{\partial p}{\partial y}(x,y), \\[1mm] \label{cont}
-\triangle p(x,y) &=& 0,
\qquad (x,y)\in(0,1)\times(0,1),\\[1mm] \nonumber p(0,y) &\equiv& 1, \qquad
y\in(0,1),\\[1mm] \nonumber
p(x,0) &\equiv& 0,\qquad x\in(1-\dd,1),\\[1mm] \nonumber
\frac{\partial p}{\partial n}(x,0) &\equiv& 0, \qquad
x\in(0,1-\dd],\\[1mm] \nonumber
\frac{\partial p}{\partial n}(x,1) &\equiv&
0,\qquad x\in[0,1],\\[1mm] \nonumber
\frac{\partial p}{\partial n}(1,y) &\equiv&
0,\qquad y\in[0,1],
\end{eqnarray}

\no in cui $u$ e $v$ sono le componenti della velocit\`a del
fluido, $p$ \`e il campo delle pressioni, $\dd$ \`e l'ampiezza
della apertura da cui il fluido fuoriesce, ed $n$ \`e il versore
normale alla frontiera del dominio.

Un modo per approssimare il campo delle pressioni consiste nel
discretizzare il dominio $[0,1]\times[0,1]$, ad esempio nel
seguente modo:
$$(ih,jh), \quad i,j=0,\dots,N+1, \qquad h=\frac{1}{N+1}.$$

\no Su tale dominio discreto, indicando con $p_{ij}$ la
approssimazione di $p(ih,jh)$, l'equazione di Laplace,
$$-\triangle p \equiv -\left(\frac{\partial^2}{\partial
x^2}p+\frac{\partial^2}{\partial y^2}p\right) = 0,$$

\no viene approssimata mediante delle differenze finite,
\begin{equation}\label{discr}
\frac{-p_{i,j-1}-p_{i-1,j}+4p_{ij}-p_{i+1,j}-p_{i,j+1}}{h^2}=0,\qquad
i,j=1,\dots,N.\end{equation}

\no Infine, le condizioni sulla frontiera si traducono in:
\begin{eqnarray*}
p_{0j}&=&1,\qquad j =1,\dots,N,\\[1mm] p_{i0}&=&0,\qquad i
=\nu+1\equiv\frac{1-\dd}h+1,\dots,N,\\[1mm] 
p_{i0}&=&p_{i1}, \qquad i =
1,\dots,\nu,\\[1mm] p_{i,N+1}&=& p_{iN},\qquad i=1,\dots,N,\\[1mm]
p_{N+1,j}&=&p_{Nj},\qquad j=1,\dots,N.
\end{eqnarray*}

\no Tenendo conto di queste condizioni, le $N^2$ equazioni in
(\ref{discr}) coinvolgono solo le incognite $p_{ij}$,
$i,j=1,\dots,N$. Definendo il vettore
$$\bfp=(p_{11},\dots,p_{N1},\dots,p_{1N},\dots,p_{NN})^\top  \in\RR^{N^2},$$

\no le (\ref{discr}) possono essere riscritte come un sistema
lineare,
\begin{equation}\label{Apb}
A\bfp = \bfb,
\end{equation}

\no la cui matrice dei coefficienti $A$, di dimensione $N^2\times
N^2$, risulta essere tridiagonale a blocchi,
\begin{equation}\label{A}
A=\pmatrix{rrrrr}
A_1 &-I_N\\
-I_N &A_2 &-I_N\\
     &\ddots &\ddots &\ddots\\
     &       & \ddots &\ddots &-I_N\\
     &       &        &-I_N   &A_N\endpmatrix,\end{equation}

\no dove i blocchi diagonali $A_i$ sono tridiagonali e simmetrici.
Pertanto, la matrice $A$ risulta essere sparsa e simmetrica. Si
dimostra, inoltre, che essa \`e anche definita positiva. Un'altra
propriet\`a importante della matrice $A$ si vede essere quella di
avere gli elementi diagonali positivi e quelli extra-diagonali non
positivi. Matrici con questa caratteristica saranno studiate in
maggior dettaglio in Sezione~\ref{Mmat}.

\begin{figure}[t]
\begin{center}
\includegraphics[width=13cm,height=10cm]{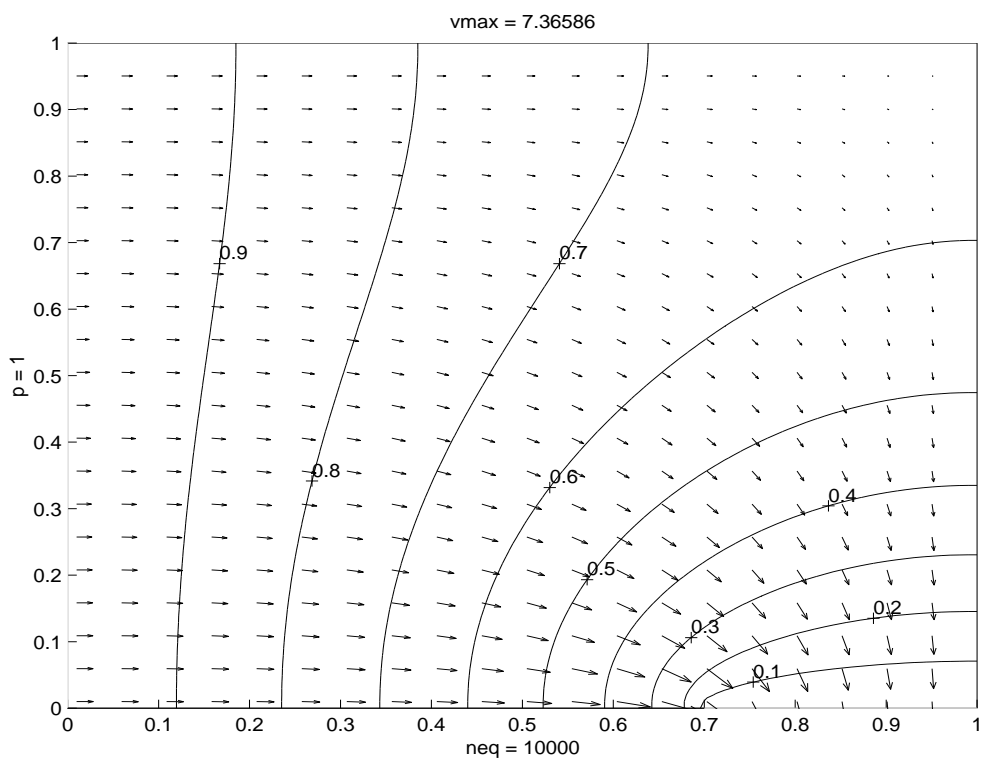}
\caption{\protect\label{fig1} Soluzione dell'equazione (\ref{A}),
approssimante (\ref{cont}), $N=100$.}
\end{center}
\end{figure}

Si pu\`o dimostrare che la soluzione discreta approssima quella
continua, proiettata sui punti del dominio discreto, con un errore
che tende a 0 quando $h$ tende a 0, ovvero, quando $N$ tende ad
infinito. Pertanto la dimensione di $A$ cresce notevolmente, al
crescere di $N$, anche se il numero di elementi non nulli \`e
minore di $5N^2$. In Figura~\ref{fig1} riportiamo il flusso
ottenuto nel caso $\dd=0.3$, in corrispondenza di $N=10^2$. In
figura sono anche indicate le isobare corrispondenti alle
pressioni $p=0.1,0.2,\dots,0.9$.

Questo \`e, quindi, un esempio al quale ben si applicano gli
argomenti che andiamo ad illustrare. Si osserva che, se il dominio
fosse stato il cubo tridimensionale unitario, una discretizzazione
uniforme con passo $h=\frac{1}{N+1}$ sulle 3 dimensioni avrebbe
portato ad un sistema lineare di $N^3$ equazioni, invece che
$N^2$, anche esso, a sua volta, sparso.

\bigskip
\begin{eser}\label{ex1}
Dimostrare che la discretizzazione dell'equazione di Poisson,

\medskip
$$-\triangle u(x,y) = x+y, \qquad
(x,y)\in\Omega=[0,1]\times[0,1],$$

\medskip
\no con la condizione $$u(x,y)= 0, \qquad(x,y)\in\partial\Omega,$$

\medskip
\no utilizzando un passo costante $$h=\frac{1}{N+1}$$ su entrambe le
direzioni spaziali, porta alla risoluzione del sistema lineare
\begin{equation}\label{test}
A\bfu = \bfb,
\end{equation}

\no dove
\begin{eqnarray*}
\bfu &=& (u_{11},\dots,u_{N1},\dots,u_{1N},\dots,u_{NN})^\top  ,\\[2mm] \bfb
&=& h^3((1+1),\dots,(N+1),\dots,(1+N),\dots,(N+N))^\top  ,\\[2mm] A &=& T_N
\otimes I_N + I_N \otimes T_N,
\end{eqnarray*}

\medskip
\no con $$u_{ij}\approx u(ih,jh),\qquad i,j=1\dots,N,$$ 

\medskip
\no e la matrice $T_N$ definita come

$$T_N = \pmatrix{rrrr} 2 &-1\\ -1 &\ddots &\ddots\\
   &\ddots &\ddots &-1\\
   &       &-1     &2\endpmatrix_{N\times N}.$$
\end{eser}

%
%
\chapter{Nozioni preliminari}\label{cap2}

{\em
In questo capitolo verranno date alcune nozioni preliminari,
riguardo a que\-stioni di algebra lineare ed equazioni alle
differenze, utili alla trattazione fatta nei successivi cinque
capitoli. In questo modo, gli unici prerequisiti che assumeremo
sono nozioni di base di algebra lineare e di teoria delle matrici,
quest'ultima confinata al solo caso delle matrici hermitiane (per
questi argomenti si vedano, ad esempio, \cite{Ba,LaTi}).}

\section{Norme indotte su matrice}

Sia $\bfx=(x_i)\in\RR^n$. La norma $p$ di $\bfx$ \`e definita da
$$\|\bfx\|_p = \left( \sum_{i=1}^n |x_i|^p \right)^{\frac{1}p}.$$

\no Norme molto utilizzate in Analisi Numerica sono:

\begin{itemize}

\item $\|\bfx\|_1 = \sum_{i=1}^n |x_i|$;

\item $\|\bfx\|_2 = \sqrt{\bfx^\top \bfx}$ ~(indotta dall'usuale prodotto
scalare);

\item $\|\bfx\|_\infty = \lim_{p\rightarrow\infty}  \left( \sum_{i=1}^n |x_i|^p
\right)^{\frac{1}p}= \max_i |x_i|$.

\end{itemize}

\no Una generica norma $\|\cdot\|$ soddisfa alle seguenti propriet\`a
caratteristiche:

\begin{description}
\item[N1:~] $\|\bfx\|\ge0,\qquad \|\bfx\|=0~\Leftrightarrow~\bfx=0$;

\item[N2:~] per ogni ~$\aa\in\RR:\|\aa\bfx\|=|\aa|\cdot\|\bfx\|$;

\item[N3:~] $\|\bfx+\bfy\|\le\|\bfx\|+\|\bfy\|$.
\end{description}

\begin{eser} Verificare che, se $Q\in\RR^{n\times n}$ \`e una
matrice ortogonale, allora $\|\bfx\|_2=\|Q\bfx\|_2$.\end{eser}

Sia ora $A\in \RR^{m\times n}$ e sia $\|\cdot\|_p$ la norma $p$ su
vettore. Definiamo {\em norma $p$ della matrice} $A$, indotta
dalla norma su vettore assegnata, la quantit\`a:
$$\|A\|_p = \sup_{\|\bfx\|_p=1} \|A\bfx\|_p.$$

\no Osserviamo che, nella precedente definizione, $\bfx\in\RR^n$ e
$A\bfx\in\RR^m$.

La definizione appena vista definisce effettivamente una norma, in
quanto si verifica facilmente che, se $A,B\in\RR^{m\times n}$ (per
brevit\`a, nel seguito ometteremo il pedice $p$ della norma, ove
non strettamente necessario):

\begin{enumerate}

\item[1.] $\|A\|\ge 0, \qquad \|A\|=0 ~\Leftrightarrow~ A=O$,~ la
matrice nulla;

\item[2.] per ogni $\aa\in\RR:~\|\aa A\| = |\aa|\cdot\|A\|$;

\item[3.] $\|A+B\| \le \|A\| + \|B\|$.

\end{enumerate}

\no Ulteriori propriet\`a, la cui semplice dimostrazione \`e
lasciata come esercizio, sono le seguenti:

\begin{enumerate}

\item[4.] la norma indotta su matrice e la corrispondente norma su vettore
sono {\em compatibili}: $\|A\bfx\|\le\|A\| \cdot\|\bfx\|$;

\item[5.] pi\`u in generale, se $A$ e $B$ sono matrici per le quali \`e definito il
prodotto $AB$: $\|AB\|\le\|A\|\cdot\|B\|$;

\item[6.] se $B$ \`e una sottomatrice di $A$, allora
$\|B\|\le\|A\|$.

\end{enumerate}

Vediamo alcune delle norme indotte di utilizzo pi\`u comune in
Analisi Numerica; siano $A=(a_{ij})\in\RR^{m\times n}$,
$\bfx=(x_i)\in\RR^n$:

\begin{itemize}

\item se $\|\bfx\|_1 = \sum_i |x_i|$, allora $\|A\|_1
= \max_j \sum_i |a_{ij}|$;

\item se $\|\bfx\|_\infty = \max_i |x_i|$, allora
$\|A\|_\infty = \max_i \sum_j |a_{ij}|$;

\item se $\|\bfx\|_2 = \sqrt{\bfx^\top \bfx}$, allora $\|A\|_2 =
\sqrt{\rho(A^\top A)}$.

\end{itemize}

\begin{eser} Verificare che, se $Q$ e $U$ sono matrici ortogonali
di dimensione $m$ ed $n$, rispettivamente, allora
$\|QAU\|_2=\|A\|_2$.\end{eser}

Sulla scorta delle precedenti propriet\`a e definizioni,
dimostriamo i seguenti risultati che riguardano il caso in cui la
matrice $A$ sia quadrata.

\begin{teo} Se $A\in\RR^{n\times n}$, allora per ogni
$\lam\in\sigma(A): |\lam|\le\|A\|$, per ogni norma indotta su
matrice.\end{teo}

\proof Sia $\lam$ un autovalore di $A$, e sia $\bfx$ un
corrispondente autovettore tale che $\|\bfx\|=1$. Si ottiene:
$$|\lam| = \|\lam\bfx\| = \|A\bfx\| \le \|A\|\cdot\|\bfx\| =
\|A\|.$$\QED

\begin{cor}\label{ro1} $\rho(A)\le\|A\|$.\end{cor}

\begin{cor} Se $A$ \`e nonsingolare, allora, per ogni
$\lam\in\sigma(A):$ $$\|A^{-1}\|^{-1}\le|\lam|.$$\end{cor}

Il seguente teorema \`e solo enunciato, sebbene la sua
dimostrazione di\-scenda da semplici argomenti di teoria delle
matrici.

\begin{teo}\label{conv} $A^k\rightarrow O$, per $k\rightarrow\infty$, se e
solo se $\rho(A)<1$.\end{teo}

\begin{defi} Una matrice $A$ che soddisfi la condizione
$\rho(A)<1$ si dice {\em convergente}.\end{defi}

\no Osserviamo che la prima implicazione, nel Teorema~\ref{conv},
pu\`o essere dimostrata ricorrendo alle norme indotte su matrice.
Infatti, se per assurdo, esistesse $\lam\in\sigma(A)$,
$|\lam|\ge1$, allora, detto $\bfx$ un corrispondente autovettore
di norma unitaria, si avrebbe, per $k\rightarrow\infty$:
$$0\leftarrow \|A^k\| \ge \|A^k\bfx\| = \|\lam^k\bfx\| =
|\lam^k|\ge1.$$

\no Inoltre, dal Corollario~\ref{ro1} discende facilmente il
seguente risultato.

\begin{cor} Se, per qualche norma indotta, $\|A\|<1$, allora la
matrice $A$ \`e convergente.\end{cor}

\section{Direzioni coniugate e norme indotte}\label{Anorm}

Gli argomenti di questa sezione sono da riguardare come una
generalizzazione del concetto di ortogonalit\`a.

\begin{defi} Data una matrice simmetrica $A$ e vettori reali
$\{\bfv_1,\dots,\bfv_k\}$ tali che $$\bfv_i^\top A\bfv_j=0, \qquad
i\ne j,$$ diremo che i vettori sono tra loro {\em
$A$-ortogonali}.\end{defi}

Evidentemente, quando $A=I$, la matrice identit\`a, allora la
nozione di $A$-ortogonalit\`a coincide con quella classica di
ortogonalit\`a.

\begin{teo}\label{Aort} Se $A$ \`e simmetrica e definita positiva (sdp) e
$\{\bfv_1,\dots,\bfv_k\}$ sono $k$ vettori non banali e tra di
loro $A$-ortogonali, allora tali vettori sono linearmente
indipendenti.\end{teo}

\proof Bisogna dimostrare che una combinazione lineare dei
predetti vettori, che dia il vettore banale, deve essere a
coefficienti tutti nulli. Infatti, se
$$\aa_1\bfv_1+\dots+\aa_k\bfv_k = \bfo,$$

\no allora, moltiplicando a sinistra per $\bfv_i^\top A$, si ottiene,
tenendo conto della $A$-ortogonalit\`a dei vettori,
$$0 = \bfv_i^\top A(\aa_1\bfv_1+\dots+\aa_k\bfv_k) = \bfv_i^\top A\bfv_i
\aa_i,$$

\no e di conseguenza $\aa_i=0$, dato che $A$ \`e sdp e che
$\bfv_i\ne\bfo$.\QED
\bigskip

Assumendo, nel seguito, che $A$ sia sdp di dimensione $n$,
osserviamo che la forma bilineare
$$\bfx^\top A\bfy$$

\no definisce un prodotto scalare. Infatti, si verifica facilmente
che, per ogni $\bfx,\bfy,\bfz\in\RR^n$ e $\aa\in\RR$:

\begin{itemize}

\item $\bfx^\top A\bfy = \bfy^\top A\bfx$;

\item $\bfx^\top A(\bfy+\bfz) = \bfx^\top A\bfy + \bfx^\top A\bfz$;

\item $(\aa\bfx)^\top A\bfy = \bfx^\top A(\aa\bfy) = \aa(\bfx^\top A\bfy)$;

\item $\bfx^\top A\bfx\ge0; \quad \bfx^\top A\bfx=0 \Leftrightarrow \bfx=\bfo$.

\end{itemize}

Come l'usuale prodotto scalare in $\RR^n$ induce la norma
euclidea, cos\`\i\, il prodotto scalare definito da $A$ induce una
norma corrispondente, che chiameremo {\em $A$-norma}:
$$\|\bfx\|_A = \sqrt{ \bfx^\top A\bfx }.$$

\no Verifichiamo che, in effetti, questa \`e una norma:

\begin{itemize}

\item $\|\bfx\|_A\ge 0$ e $\|\bfx\|_A=0 \Leftrightarrow
\bfx=\bfo$. Questo discende dal fatto che $A$ \`e sdp;

\item se $\aa\in\RR$, $\|\aa\bfx\|_A = |\aa|\cdot\|\bfx\|_A$.
Anche questa propriet\`a \`e evidente;

\item $\|\bfx+\bfy\|_A\le \|\bfx\|_A+\|\bfy\|_A$. Infatti, dalla
diseguaglianza di Cauchy-Schwarz e tenendo conto che se $A$ \`e
sdp tale \`e anche $A^{1/2}$ (che \`e una matrice reale), segue
che:

\begin{eqnarray*}
\|\bfx+\bfy\|_A^2 &=& (\bfx+\bfy)^\top A(\bfx+\bfy) = \bfx^\top A\bfx
+\bfy^\top A\bfy +2\bfx^\top A\bfy\\[1mm]
                  &=& \|\bfx\|_A^2 +\|\bfy\|_A^2
                  +2(\bfx^\top A^{1/2})(A^{1/2}\bfy)\\[1mm]
                  &\le& \|\bfx\|_A^2 +\|\bfy\|_A^2
                  +2\sqrt{(\bfx^\top A\bfx)(\bfy^\top A\bfy)}\\[1mm]
                  &=& \|\bfx\|_A^2
                  +\|\bfy\|_A^2+2\|\bfx\|_A\|\bfy\|_A =
                  (\|\bfx\|_A+\|\bfy\|_A)^2.
                  \end{eqnarray*}

\end{itemize}

\begin{oss} Chiaramente, la precedente norma su vettore ne pu\`o
indurre una corrispondente su matrice, per matrici che abbiano
le stesse dimensioni di $A$.\end{oss}

\section{Polinomi di matrici}\label{polmat}

Sia assegnata una matrice $A\in\RR^{n\times n}$ ed un polinomio
$p(x)\in\Pi_k$,
$$p(x) = \sum_{i=0}^k c_i x^i.$$

\no Il {\em polinomio di matrice} $p(A)$ \`e quindi definito da
$$p(A) = \sum_{i=0}^k c_i A^i.$$

\no Per tali polinomi, e pi\`u in generale per funzioni di
matrice, valgono molte e interessanti propriet\`a che dipendono
dalla struttura di $A$ (si vedano, ad esempio, \cite{LaTr,LaTi}).
Noi siamo interessati solo ad alcune semplici propriet\`a che
riportiamo nel seguito. Innanzitutto, se la forma canonica di
Jordan di $A$ \`e $$A=V\Lambda V^{-1},$$

\no allora si dimostra facilmente per induzione che $$A^i =
V\Lambda^iV^{-1}, \qquad i=0,1,2,\dots.$$

\no Ne consegue che $$p(A) = \sum_{i=0}^k c_i A^i = \sum_{i=0}^k
c_i V\Lambda^i V^{-1} = V\left(\sum_{i=0}^k c_i
\Lambda^i\right)V^{-1} = Vp(\Lambda)V^{-1}.$$

\no Nel caso in cui $A$ sia diagonalizzabile, ovvero quando
$$\Lambda=\pmatrix{ccc}\lam_1\\ &\ddots\\&&\lam_n\endpmatrix,$$

\no si perviene alla, pi\`u semplice, espressione $$p(A) =
Vp(\Lambda)V^{-1} = V\pmatrix{ccc}
p(\lam_1)\\&\ddots\\&&p(\lam_n)\endpmatrix V^{-1}.$$

\no Questo sicuramente avviene quando $A$ \`e simmetrica, nel
quale caso la matrice $V$ \`e una matrice ortogonale e, pertanto,
$$p(A) = Vp(\Lambda)V^\top  = V\pmatrix{ccc}
p(\lam_1)\\&\ddots\\&&p(\lam_n)\endpmatrix V^\top .$$

\bigskip
\begin{oss} $A$ e $p(A)$ commutano, ovvero $Ap(A)=p(A)A$. \end{oss}

\bigskip
\begin{oss} Se $A$ \`e simmetrica, allora la matrice $p(A)$
\`e ancora simmetrica.\end{oss}

\bigskip
\begin{oss} Se $A$ \`e simmetrica ed $f(x)$ \`e una funzione
analitica in un dominio contenente lo spettro, $\sigma(A)$, di
$A$, allora $f(A)\equiv p(A)$, dove $p(x)$ \`e il polinomio che
interpola $f(x)$ sui punti di $\sigma(A): \lam\in\sigma(A)
\Rightarrow p(\lam) = f(\lam)$. \\Per approfondimenti sulle
funzioni di matrice, si vedano \cite{LaTr,LaTi}.\end{oss}

\bigskip
\begin{eser} Se $\lam\in\sigma(A)$, allora
$p(\lam)\in\sigma(p(A))$.\end{eser}

\bigskip
\begin{eser} Sia $L\in\RR^{n\times n}$ una matrice strettamente
triangolare. Dimostrare che $$(I-L)^{-1} = I
+L+L^2+\dots+L^{n-1}.$$ Concludere, quindi, che
$$(I-L)^{-1}L=L(I-L)^{-1}.$$\end{eser}

\section{Fattorizzazione LU di una matrice}\label{fattLUmat}

In questa sezione si richiamano le definizioni e le propriet\`a
essenziali (ai fini della trattazione seguente) riguardo alla
fattorizzazione $LU$ di una matrice nonsingolare $A\in\RR^{n\times
n}$.

\begin{defi} Se esistono due matrici triangolari, rispettivamente
inferiore e superiore,
\begin{eqnarray*}
L &=& \pmatrix{ccc}
\ell_{11}\\
\vdots &\ddots\\
\ell_{n1}&\dots&\ell_{nn}\endpmatrix, \qquad \ell_{ii}=1, \quad
i=1,\dots,n,\\[2mm]
U &=& \pmatrix{ccc}
u_{11} &\dots &u_{1n}\\
       &\ddots &\vdots\\
       &       &u_{nn}\endpmatrix,
       \end{eqnarray*}

\no tali che \begin{equation}\label{lu} A = LU,\end{equation}
allora $A$ dicesi {\em fattorizzabile $LU$}.\end{defi}

Valgono le seguenti propriet\`a:

\begin{description}

\item[LU1:~] se $A$ \`e fattorizzabile $LU$, la sua
fattorizzazione (\ref{lu}) \`e unica;

\item[LU2:~] $A$ \`e fattorizzabile $LU$ se e solo se tutti i suoi
minori principali sono non nulli.

\end{description}

\no L'ultima propriet\`a \`e senz'altro vera nel caso di matrici
sdp, i cui minori principali sono tutti positivi. In tal caso, si
vede facilmente dalla LU1 che $A$ \`e fattorizzabile nella forma
\begin{equation}\label{Aldlt}A = L D L^\top ,\end{equation}
\no con $D$ matrice diagonale. Ovvero la matrice $U$ della
fattorizzazione $LU$ di $A$ \`e scrivibile nella forma $U=DL^\top $.

\begin{eser} Verificare che, nel caso in cui $A$ sia tridiagonale
e sdp,
$$A = \pmatrix{cccc}
a_1 &b_1\\
b_1 &\ddots &\ddots\\
    &\ddots &\ddots &b_{n-1}\\
    &       &b_{n-1} &a_n\endpmatrix,$$

\no allora essa \`e fattorizzabile nella forma (\ref{Aldlt}) con
$$
L = \pmatrix{cccc}
1\\
\mu_1 &1\\
&\ddots &\ddots\\
&       &\mu_{n-1} &1\endpmatrix, \qquad
D = \pmatrix{ccc} d_1 \\ &\ddots\\ &&d_n\endpmatrix,$$

\no dove $d_1=a_1$ e
$$\mu_i = \frac{b_i}{d_i}, \qquad d_{i+1}=a_{i+1} -b_i\mu_i,
\qquad i=1,\dots,n-1.$$
\end{eser}

\section{Fattorizzazione QR di una matrice}

Tratteremo ora la fattorizzazione $QR$ di una matrice
$A\in\RR^{m\times n}, m\ge n$. Gli argomenti esposti in questa
sezione riguardano pi\`u gli aspetti teorici che non quelli
computazionali. Per una completa trattazione di questi ultimi, si
veda, ad esempio, \cite{GoVL}.

Preliminarmente, poniamoci il problema, assegnato un vettore
$\bfv\in\RR^m$, di annullarne tutte le componenti a partire dalla
$(k+1)$-esima, per un assegnato $k\in\{1,\dots,m-1\}$, preservando
la norma euclidea del vettore.

A questo fine, osserviamo che il problema pu\`o essere decomposto
in $m-k$ sottoproblemi consecutivi. Infatti, se
$\bfv=(v_1,\dots,v_m)^\top $, allora possiamo pensare di annullare,
innanzitutto, l'ultima componente, considerando le ultime due. Lo
strumento per ottenere questo sar\`a una {\em matrice di
rotazione} tale che: $$\pmatrix{rr}c_m &s_m\\-s_m &c_m\endpmatrix
\pmatrix{c} v_{m-1}\\v_m\endpmatrix = \pmatrix{c} \hv_{m-1} \\
0\endpmatrix,$$

\no dove gli scalari $c_m, s_m$ sono individuati (a meno del
segno) dalle equazioni $$c_m^2+s_m^2=1, \qquad
c_mv_m=s_m v_{m-1}.$$ Si ricava:
$$\rho_m=\|(v_{m-1},v_m)^\top \|_2, \qquad c_m = v_{m-1}/\rho_m,
\qquad s_m = v_m/\rho_m.$$ Chiaramente, se $v_{m-1}=v_m=0$, \`e
sufficiente porre $c_m=1,s_m=0$. In questo modo, abbiamo
trasformato il vettore $(v_1,\dots,v_m)^\top $ nel vettore
$(v_1,\dots,v_{m-2},\hv_{m-1},0)^\top $.

\begin{eser} Verificare che i due precedenti vettori hanno la
stessa norma euclidea.\end{eser}

La procedura pu\`o essere, quindi, applicata iterativamente al
sottovettore di componenti $(v_{m-2},\hv_{m-1})^\top $ che, mediante
moltiplicazione a sinistra per la matrice di rotazione
$$\pmatrix{cc} c_{m-1} &s_{m-1}\\ -s_{m-1}&c_{m-1}\endpmatrix,$$
\`e trasformato nel vettore $(\hv_{m-2},0)^\top $, e cos\`\i\, via.

In notazione vettoriale, avremo operato la trasformazione
\begin{equation}\label{giv}
G_{k+1}\cdot G_{k+2}\cdots G_m \bfv =
(v_1,\dots,v_{k-1},\hv_k,0,\dots,0)^\top \equiv \hat\bfv,
\end{equation}

\no con le {\em matrici di Givens} $G_j$ definite da $$G_j =
\pmatrix{ccc} I_{j-2} \\ &\pmatrix{rr} c_j &
s_j\\-s_j&c_j\endpmatrix \\ &&I_{m-j}\endpmatrix, \qquad
j=k+1,\dots,m.$$

\bigskip
\begin{eser} Verificare che le matrici $G_j$ in (\ref{giv}) sono
ortogonali e che i vettori $\bfv$ e $\hat\bfv$ hanno la stessa
norma euclidea.\end{eser}

Sia data, ora, una matrice rettangolare $A\in\RR^{m\times
n}$, in cui supporremo che $m\ge n\equiv \rank(A)$. Vale il seguente
risultato.

\begin{teo}\label{QRfatt} $A$ \`e fattorizzabile nella forma
$A=QR$, dove $Q\in\RR^{m\times m}$ \`e ortogonale e

\begin{equation}\label{R}
R=\pmatrix{c}R_1 \\O\endpmatrix
\in\RR^{m\times n},\end{equation}

\no con $R_1$ matrice $n\times n$ triangolare superiore e
nonsingolare.\end{teo}

\proof Che $R_1$ sia nonsingolare discende dal fatto che $A$ ha
rango massimo e $Q$ \`e nonsingolare. Pertanto rimane solo da
dimostrare che la fattorizzazione $A=QR$ esiste. Usando argomenti
simili a quelli visti per derivare la (\ref{giv}), \`e possibile
definire $m-1$ matrici di Givens $G_m^{(1)},\dots,G_2^{(1)}$ tali
che $$G_2^{(1)}\cdots G_m^{(1)} A = \pmatrix{cccc} a_{11}^{(1)}
&a_{12}^{(1)} &\dots &a_{1n}^{(1)}\\ 0 &a_{22}^{(1)} &\dots
&a_{2n}^{(1)}\\ \vdots &\vdots & &\vdots\\ 0 &a_{m2}^{(1)} &\dots
&a_{mn}^{(1)}\endpmatrix \equiv A^{(1)}.$$ Similmente, sar\`a
possibile definire $m-2$ matrici di Givens
$G_m^{(2)},\dots,G_3^{(2)}$ tali che $$G_3^{(2)}\cdots G_m^{(2)}
A^{(1)} = \pmatrix{ccccc} a_{11}^{(1)} &a_{12}^{(1)} &a_{13}^{(1)}
&\dots &a_{1n}^{(1)}\\ 0 &a_{22}^{(2)} &a_{23}^{(2)}&\dots
&a_{2n}^{(2)}\\ \vdots &0& a_{33}^{(2)} &\dots &a_{3n}^{(2)}\\
\vdots &\vdots &\vdots& &\vdots\\ 0 &0&a_{m3}^{(2)} &\dots
&a_{mn}^{(2)}\endpmatrix \equiv A^{(2)}.$$ Ragionando per
induzione, si definiranno ulteriori matrici di Givens $$G_i^{(j)},
\qquad i=j+1,\dots,m,\quad j=3,\dots,n,$$ tali che, infine, $$Q^\top A
\equiv G_{n+1}^{(n)}\cdots G_m^{(n)} \cdots G_2^{(1)} \cdots
G_m^{(1)} A =\pmatrix{ccc} a_{11}^{(1)} &\dots &a_{1n}^{(1)} \\ 0
&\ddots &\vdots \\ \vdots &\ddots &a_{nn}^{(n)}
\\ \vdots & &0\\ \vdots & &\vdots\\ 0&\dots &0\endpmatrix\equiv
\pmatrix{c} R_1
\\O\endpmatrix.$$ \QED

\begin{oss}\label{QRT} \`E da sottolineare il fatto che la generica matrice di Givens
$G_i^{(j)}$ assolve alla funzione di azzerare, in modo selettivo,
l'elemento in posizione $(i,j)$ della matrice corrente. Nel fare
questo, modifica il contenuto delle righe $i-1$ ed $i$, senza
reintrodurre {\em fill-in} negli elementi gi\`a azzerati, a patto
che l'ordine delle matrici sia quello specificato. Questo permette
di ridurre il numero delle matrici necessarie per ottenere la
fattorizzazione, quando $A$ \`e in forma di {\em Hessemberg}
(superiore), $$A = \pmatrix{ccc} a_{11} &\dots  & a_{1n} \\ a_{21}
&\ddots &\vdots  \\ 0      &\ddots &a_{nn} \\ \vdots &\ddots
&a_{n+1,n}\\ \vdots &       &0      \\ \vdots &       &\vdots \\ 0
&\dots  &0\endpmatrix.$$ Infatti, si verifica agevolmente che, in
questo caso, sono richieste, nell'ordine, le sole matrici
$$G_2^{(1)},G_3^{(2)},\dots,G_{n+1}^{(n)}$$ per ottenere la
fattorizzazione $QR$ di $A$.

Per ultimo, osserviamo che, se $A$ \`e anche {\em tridiagonale},
$$A = \pmatrix{ccccc}
a_{11} &a_{12} & 0     &\dots   &0 \\
a_{21} &\ddots &\ddots &\ddots  &\vdots  \\
0      &\ddots &\ddots &\ddots  &0 \\
\vdots &\ddots &\ddots &\ddots  &a_{n-1,n}\\
\vdots &       &\ddots &a_{n,n-1}&a_{nn}\\
\vdots &       &       &\ddots  &a_{n+1,n}\\
\vdots &       &       &        &0\\
\vdots &       &       &        &\vdots\\
0      &\dots  &\dots  &\dots   &0\endpmatrix,$$
allora tale \`e anche il suo fattore $R$. Infatti, si verifica che
(vedi (\ref{R})) $$R_1 = \pmatrix{cccccc}
r_{11} &r_{12} &r_{13} &0      &\dots  &0\\
       &\ddots &\ddots &\ddots &\ddots &\vdots\\
       &       &\ddots &\ddots &\ddots &0\\
       &       &       &\ddots &\ddots &r_{n-2,n}\\
       &       &       &       &\ddots &r_{n-1,n}\\
       &       &       &       &       &r_{nn}\endpmatrix.$$
\end{oss}

\section{Il problema dei minimi quadrati}\label{minquad}

Sia assegnata una funzione
$$f:\bfx\in\RR^m\rightarrow\RR,$$

\no di classe $C^2$, di cui si vuole ricercare il minimo. Se non
vi sono vincoli sulla $\bfx$, allora ricercheremo la soluzione tra
i punti stazionari della $f$, ovvero risolvendo l'equazione
$$\nabla f(\bfx)=\bfo^\top .$$

\no Se $\bfx^*$ \`e soluzione di tale equazione, allora esso
sar\`a un punto di minimo stretto (locale, in generale) per $f$ se
la matrice Hessiana, $\nabla^2f(\bfx^*),$ \`e sdp.

Supponiamo, adesso, di voler ricercare un punto di minimo della
$f$ non in tutto $\RR^m$ ma, piuttosto, nella variet\`a affine
$$\Psi = \left\{\bfx\in\RR^m: \bfx = \bfx_0 +B\bfz,
~\bfz\in\RR^n\right\},$$

\no dove $\bfx_0$ \`e un punto assegnato e $B\in\RR^{m\times n}$,
$n\le m$, \`e una matrice avente rango massimo $n$. In tal caso,
$$\min_{\bfx\in\Psi} f(\bfx) = \min_{\bfz\in\RR^n} f(\bfx_0+B\bfz)
\equiv \min_{\bfz\in\RR^n} f(g(\bfz)).$$

\no Applicando la regola di derivazione per funzioni composte, si
ottiene che, se $\bfx^* = g(\bfz^*)$ \`e il punto di minimo,
allora dovr\`a aversi
$$\nabla f(\bfx^*) B = \bfo^\top .$$

\no Ovvero, il gradiente della $f$, nel punto di minimo, deve
essere ortogonale alle colonne della matrice $B$ che individua lo
spazio affine. Questa condizione, che \`e in generale solo
necessaria, diviene anche sufficiente nel caso del problema
lineare dei minimi quadrati. In tal caso, infatti, per un
assegnato vettore $\bfb\in\RR^m,$

\begin{equation}\label{lin}f(\bfx) = \frac{1}2\|\bfb-\bfx\|_2^2
\equiv \frac{1}2(\bfb-\bfx)^\top (\bfb-\bfx).\end{equation}

\no Al solito, se $\bfx$ potesse variare in tutto $\RR^m$, allora
la scelta banale, $\bfx=\bfb$, sarebbe la soluzione del problema
di minimo. Se, tuttavia, restringiamo $\bfx\in\Psi$, allora il
punto stazionario sar\`a soluzione dell'equazione
\begin{equation}\label{errort}
(\bfb-\bfx)^\top B = \bfo^\top .
\end{equation}

\no Tenendo conto che $\bfx=\bfx_0+B\bfz$, questo porta al sistema
di equazioni lineari
\begin{equation}\label{first}
B^\top B \bfz = B^\top (\bfb-\bfx_0).
\end{equation}

\no Poich\`e la matrice dei coefficienti, $B^\top B$, \`e sdp (vedi
Esercizio~\ref{sdpbtb}), ne consegue che l'\underline{unica} soluzione \`e
\begin{equation}\label{zstar}\bfz^* =
(B^\top B)^{-1}B^\top (\bfb-\bfx_0).\end{equation}

Osserviamo che il corrispondente punto in $\Psi$,
\begin{equation}\label{sol}\bfx^* = \bfx_0+B\bfz^*,\end{equation}

\no sar\`a un punto di minimo, in quanto $\nabla^2f(\bfx^*)=B^\top B,$ che \`e
sdp. Quindi, come detto innanzi, nel caso del problema lineare dei minimi
quadrati, la condizione necessaria (\ref{errort}) \`e anche sufficiente.

\medskip
\begin{eser}\label{sdpbtb} Dimostrare che, se $B\in\RR^{m\times
n}, m\ge n,$ ed il rango di $B$ \`e uguale a $n$, allora $B^\top B$
\`e sdp.\end{eser}

\medskip
\begin{eser} Determinare l'approssimazione ai minimi quadrati del
vettore $(\aa~\aa~\bb)^\top $ nel sottospazio di $\RR^3$ generato dai
vettori $(1,~1,~1)^\top $ e $(1,~-1,~0)^\top $.
\end{eser}

\medskip
\begin{oss} Il problema (\ref{lin}) pu\`o essere generalizzato al caso di una
qualunque norma indotta da una matrice sdp $A$ (vedi
Sezione~\ref{Anorm}). Infatti, in tal caso, la (\ref{lin}) diviene
\begin{equation}\label{Anorm1}
f(\bfx) =
\frac{1}2\|\bfb-\bfx\|_A^2\equiv\frac{1}2(\bfb-\bfx)^\top A(\bfb-\bfx),
\end{equation}

\no ed il sistema lineare (\ref{first}) diviene
$$(B^\top AB)\bfz = B^\top A(\bfb-\bfx_0),$$

\no la cui matrice dei coefficienti \`e ancora sdp.
\end{oss}

\medskip
\begin{oss}\label{errort1} Osserviamo che la condizione
(\ref{errort}), che \`e necessaria e sufficiente per il punto di
minimo, nel caso del problema lineare dei minimi quadrati, pu\`o
essere interpretata come la condizione di ortogonalit\`a dell'{\em
errore}, $$\bfb-\bfx,$$ ai vettori che generano il sottospazio.
Nel caso della $A$-norma (\ref{Anorm1}), la corrispondente
condizione, $$(\bfb-\bfx)^\top AB = \bfo^\top ,$$ pu\`o essere
analogamente interpretata come la $A$-ortogonalit\`a dell'errore
alle colonne di $B$.\end{oss}

\medskip
\begin{oss} Il vettore $\bfz^*$ non \`e convenientemente
ottenuto risolvendo (\ref{lin}). \`E invece preferibile risolvere,
in alternativa, il problema di minimo:

$$\min_{\bfz\in\RR^n} \|B\bfz-(\bfb-\bfx_0)\|_2^2.$$

\no Infatti, fattorizzando $B=QR,$ con $R$ come in (\ref{R}), si
ottiene:

$$ \|B\bfz-(\bfb-\bfx_0)\|_2^2 = \|R\bfz -Q^\top (\bfb-\bfx_0)\|_2^2
\equiv \|R\bfz-\bfg\|_2^2.$$

\no Partizionando il vettore $\bfg=(\bfg_1,\bfg_2)^\top $, con
$\bfg_1\in\RR^n$, si ottiene, quindi,

$$\|R\bfz-\bfg\|_2^2 = \|R_1\bfz-\bfg_1\|_2^2 +\|\bfg_2\|_2^2,$$

\no che \`e minimizzato dalla scelta $\bfz^*=R_1^{-1}\bfg_1$,
essendo $R_1$ nonsingolare (vedi Teorema~\ref{QRfatt}).
\end{oss}

\medskip
\begin{oss} La soluzione (\ref{sol}) \`e invariante per scalamento
delle colonne della matrice $B$. Infatti, se
$$\Sigma = \pmatrix{ccc} \sigma_1 \\ &\ddots
\\&&\sigma_n\endpmatrix$$

\no \`e nonsingolare, allora
$$\bfx^* = \bfx_0 +B\bfz^* = \bfx_0 +(B\Sigma)(\Sigma^{-1}\bfz^*)
\equiv \bfx_0 +\hat{B}\hat\bfz,$$

\no con (vedi (\ref{zstar}))
\begin{eqnarray*}
\hat\bfz &\equiv& \Sigma^{-1}\bfz^* = \Sigma^{-1} (B^\top B)^{-1}
B^\top (\bfb-\bfx_0)\\
         &=& \left( (B\Sigma)^\top (B\Sigma) \right)^{-1}
         (B\Sigma)^\top (\bfb-\bfx_0)\\
         &\equiv& (\hat{B}^\top \hat{B})^{-1}\hat{B}^\top (\bfb-\bfx_0).
         \end{eqnarray*}

\no Cambia, quindi, solo la sua rappresentazione.
\end{oss}

\section{$M$-matrici}\label{Mmat}
In questa sezione riportiamo alcune propriet\`a relative alle
cosiddette $M$-matrici. Saranno esaminate solo le propriet\`a
strettamente attinenti alla trattazione seguente. Per maggiori
approfondimenti, si rimanda a \cite{Ba,LaTi,Va}.

\begin{defi} Diremo che una matrice $B=(b_{ij})$ \`e {\em non
negativa}, e scrive\-remo $B\ge0$, se per ogni $i,j: b_{ij}\ge0.$
\end{defi}

Valgono le seguenti, semplici, propriet\`a, per matrici quadrate
non negative.

\medskip
\begin{teo}\label{posmat} Se $A,B\ge0$, allora $AB\ge0$;
se $A\ge B\ge0$, allora $A^j\ge B^j\ge0$,\quad $j=0,1,2,\dots.$\end{teo}

\medskip
Si pu\`o inoltre dimostrare il seguente risultato, che discende dal
{\em Teorema di Perron-Frobenius} (vedi, ad esempio, \cite{LaTi}).

\medskip
\begin{teo}\label{PF} Se $A\ge0$, allora esiste $\lam\in\sigma(A)$ tale che
$\lam\equiv\rho(A)$. Ad esso corrisponde un autovettore $\bfv\ge0$.\end{teo}

\medskip
Per esemplificare il risultato del precedente teorema, si consideri $$A =
\pmatrix{cc} 0 &1\\1&0\endpmatrix\ge0.$$ Si ha:
$$\sigma(A)=\{-1,1\}, \qquad A\pmatrix{c}1\\1\endpmatrix =
\pmatrix{c}1\\1\endpmatrix\ge0.$$

\medskip
\begin{eser} Determinare il raggio spettrale della matrice $$A=
\pmatrix{ccccc}
1&2&3&4&5\\2&3&4&5&1\\3&4&5&1&2\\4&5&1&2&3\\5&1&2&3&4\endpmatrix.
$$
\end{eser}

\medskip
Vale, infine, la seguente propriet\`a di monotonia.

\medskip
\begin{teo}\label{monot} Se $0\le B \le C$, allora
$\rho(B)\le\rho(C)$.\end{teo}

\proof Infatti, si verifica facilmente che, per ogni $\eps>0$,
$$0\le\frac{C}{\rho(C)+\eps}$$

\no \`e una matrice convergente. Dalle ipotesi, e dal
Teorema~\ref{posmat}, segue quindi che, per ogni $\eps>0$, $$0\le
\left(\frac{B}{\rho(C)+\eps}\right)^j \le
\left(\frac{C}{\rho(C)+\eps}\right)^j \rightarrow 0,
\mbox{~~per~}j\rightarrow\infty.$$

\no Se ne conclude che, per ogni $\eps>0$, $\rho(B)<\rho(C)+\eps$.
La tesi segue passando al limite per $\eps\rightarrow0$.\QED
\bigskip

La definizione di matrice non negativa ci consente di dare anche
quella seguente.

\begin{defi} Una matrice $A\in\RR^{m\times m}$ \`e una $M$-matrice
se
\begin{equation}\label{Mm}
A = \aa I -B,\qquad B\ge0,\qquad\aa>\rho(B).\end{equation}

\no Se $A=A^\top $, ovvero $B=B^\top $, la matrice dicesi {\em matrice di
Stieltjes}.
\end{defi}

\medskip
\begin{oss} Gli elementi extradiagonali di una $M$-matrice sono non
positivi. La matrice $A$ in (\ref{A}) \`e una matrice di
Stieltjes.\end{oss}

\medskip
Valgono le seguenti propriet\`a, la cui dimostrazione \`e lasciata
come esercizio.

\medskip
\begin{teo} Una $M$-matrice \`e nonsingolare; una matrice di
Stieltjes \`e definita positiva.\end{teo}

\medskip
La propriet\`a fondamentale di una $M$-matrice \`e descritta nel
seguente teorema.

\medskip
\begin{teo} Se $A$ \`e una $M$-matrice, allora
$A^{-1}\ge0$.\end{teo}

\proof Se $A$ \`e una $M$-matrice, allora essa \`e scrivibile
nella forma (\ref{Mm}). Inoltre, la matrice $\aa^{-1}B$ risulta
essere non negativa e convergente. Ne consegue che
$$A^{-1} = (\aa I-B)^{-1} = \aa^{-1}\sum_{j=0}^\infty
(\aa^{-1}B)^j \ge0.$$\QED
\vspace{2mm}

Da questa propriet\`a discendono immediatamente le seguenti altre.

\begin{description}

\item[M1:~] se $A$ e $B$ sono $M$-matrici, $A\le B$, allora:
$$A^{-1}\ge B^{-1},\qquad B^{-1}A\le I\le A^{-1}B,\qquad
AB^{-1}\le I\le BA^{-1};$$

\item[M2:~] dato il sistema di diseguaglianze $A\bfx\le\bfb$,
allora, se $A$ \`e una $M$-matrice, segue che $\bfx\le
A^{-1}\bfb$;

\item[M3:~] gli elementi diagonali di una $M$-matrice sono
positivi.

\end{description}

\medskip
\begin{eser} Dimostrare le propriet\`a M1-M3.\end{eser}

Sulla scorta del Teorema~\ref{monot}, si pu\`o dimostrare il
seguente risultato.

\medskip
\begin{teo} Se $A_1 = \aa I-B_1 \le A_2 = \aa I-B_2$, $B_1,B_2\ge0$, e $A_1$ \`e
una $M$-matrice, allora lo \`e anche $A_2$.\end{teo}

\proof Infatti, si ha che $0\le B_2\le B_1$ e, quindi,
$\rho(B_2)\le \rho(B_1)<\aa$.\QED

\medskip
\begin{cor} Se la matrice $A_1$ \`e ottenuta da una $M$-matrice
$A$ azzerando qualche elemento extradiagonale di quest'ultima,
allora $A_1$ \`e, a sua volta, una $M$-matrice.\end{cor}

\section{Equazioni alle differenze lineari}

Una {\em equazione alle differenze lineare a coefficienti costanti
di ordine $k$} ha la forma

\begin{equation}\label{completa}
\sum_{i=0}^k \aa_i x_{n+i}=\bb_n, \qquad n=0,1,2,\dots,
\end{equation}

\no in cui i coefficienti $\{\aa_i\}$ sono fissati e la
successione $\{\bb_n\}$ \`e nota. Quando quest'ultima coincide con
la successione nulla, l'equazione si dice {\em omogenea}. Nel
seguito, considereremo il caso di equazioni omogenee: per una
trattazione esaustiva sulle equazioni alle differenze, si rimanda
a \cite{LaTr}. Se $\aa_0\aa_k\ne0$, la precedente equazione avr\`a
esattamente ordine $k$: questo sar\`a assunto si\-stematicamente.
Inoltre, poich\`e la equazione (\ref{completa}) \`e definita a
meno di una costante moltiplicativa non nulla, assumeremo la
normalizzazione $$\aa_k=1.$$ Una successione (eventualmente ad
elementi complessi) $\{x_n\}$ che soddisfi la precedente
equazione, ne costituisce una {\em soluzione}. Nel seguito
studieremo, quindi, le soluzioni della equazione
\begin{equation}\label{omo}
\sum_{i=0}^k \aa_i x_{n+i}=0, \quad n=0,1,2,\dots, \qquad
\aa_k=1,~\aa_0\ne0.
\end{equation}

\`E evidente che, essendo $\aa_k\ne0$, una qualunque soluzione di
(\ref{omo}) rimane univocamente individuata assegnandone le prime
$k$ componenti,
\begin{equation}\label{ini}
x_0,x_1,\dots,x_{k-1},
\end{equation}

\no ovvero, assegnando $k$ {\em condizioni iniziali}. Siamo
interessati, a riguardo, a stabilire la forma generale della
soluzione di (\ref{omo}), in funzione delle condizioni iniziali
(\ref{ini}). Preliminarmente, enunciamo il seguente risultato.

\medskip
\begin{teo} L'insieme delle soluzioni dell'equazione (\ref{omo})
\`e uno spazio vettoriale di dimensione $k$.\end{teo}

\proof Il fatto che l'insieme delle soluzioni sia uno spazio
vettoriale di\-scende facilmente dal fatto che, considerate due
qualunque soluzioni $\{x_n\}$ e $\{y_n\}$ di (\ref{omo}), e
comunque si fissino due scalari $\cc$ e $\dd$, la successione
$\{\cc x_n+\dd y_n\}$ \`e ancora una soluzione di (\ref{omo}). Il
fatto che la dimensione dello spazio sia $k$ deriva, infine, dal
fatto che per individuare univocamente una soluzione
dell'equazione si devono imporre le $k$ condizioni iniziali
(\ref{ini}).\QED \bigskip

Da questo risultato, segue che possiamo scrivere la generica
soluzione dell' equazione (\ref{omo}) come combinazione lineare di
$k$ soluzioni linearmente indipendenti. Se
$\{x_n^{(1)}\},\dots,\{x_n^{(k)}\}$ sono tali soluzioni, allora,
la generica soluzione di (\ref{omo}) sar\`a della forma $$x_n =
\sum_{i=1}^k c_i x_n^{(i)},$$ dove le costanti $c_1,\dots,c_k$
saranno univocamente determinate dalle $k$ condizioni iniziali
(\ref{ini}). Ci\`o premesso, andiamo a ricercare soluzioni di
(\ref{omo}) nella forma $$x_n = z^n,$$ in cui $z$ \`e uno scalare,
eventualmente complesso. Si vede facilmente che, se la soluzione
\`e non banale, allora $z$ deve essere una radice del {\em
polinomio caratteristico},
\begin{equation}\label{car}
p(z) = \sum_{i=0}^k \aa_i z^i,
\end{equation}

\no associato all'equazione (\ref{omo}). Nel seguito assumeremo
che le $k$ radici di tale polinomio siano semplici, in quanto
questo \`e il caso che ci riguarder\`a nel seguito. Per l'analisi
del caso pi\`u generale, in cui sia contemplato anche il caso di
radici multiple, si veda \cite{LaTr} oppure \cite[cap.\,2]{BrTr}.

Siano, dunque, $z_1,z_2,\dots,z_k$ le radici (distinte) del
polinomio (\ref{car}). Vale il seguente risultato.

\begin{teo} La soluzione generale di (\ref{omo}) \`e
\begin{equation}\label{gen}
x_n = \sum_{i=1}^k c_i z_i^n,\qquad n=0,1,\dots,\end{equation}
dove le costanti $c_1,\dots,c_k$ sono univocamente individuate
dalle condizioni iniziali (\ref{ini}).\end{teo}

\proof Infatti, la (\ref{gen}) \`e evidentemente soluzione di
(\ref{omo}). Viceversa, basta dimostrare che si possono scegliere
le costanti $c_1,\dots,c_k$ in modo da ottenere una qualunque
$k$-pla di valori iniziali (\ref{ini}). Dalla (\ref{gen}), per
$n=0,1,\dots,k-1$, si ottiene il sistema di equazioni lineare:
$$V\bfc = \bfx,$$

\no in cui $$\bfc=\pmatrix{c}c_1\\ \vdots\\c_k\endpmatrix, \qquad \bfx =
\pmatrix{c}x_0\\ \vdots\\x_{k-1}\endpmatrix, \qquad V = \pmatrix{ccc} z_1^0 & \dots
&z_k^0\\ \vdots & &\vdots\\ z_1^{k-1} &\dots
&z_k^{k-1}\endpmatrix.$$

\no La matrice $V$ risulta essere una matrice
di {\em Vandermonde}, che \`e noto essere nonsingolare se e solo
se gli scalari $z_1,\dots,z_k$ che la definiscono sono distinti,
come abbiamo assunto. Pertanto \begin{equation}\label{c}\bfc =
V^{-1}\bfx.\end{equation}\QED
\vspace{2mm}

Una immediata, e assai importante, conseguenza segue nel caso in
cui esi\-sta una {\em radice dominante}, ovvero di modulo massimo e
che separa in modulo tutte le altre.

\begin{cor}\label{initval} Se in (\ref{gen}) $c_k\ne0$ e
$$|z_1|\le|z_2|\le \dots \le |z_{k-1}|<|z_k|,$$ allora
$$x_n \approx c_k z_k^n, \qquad n\gg0.$$\end{cor}

\proof Infatti, per il problema discreto ai {\em valori iniziali}
(\ref{omo})-(\ref{ini}), vale la (\ref{c}), in cui le componenti
del vettore $\bfc$ sono indipendenti da $n$. La tesi, quindi, segue
dalla (\ref{gen}).\QED

\begin{oss} Le conclusioni del Corollario~\ref{initval} non sono
pi\`u applicabili nel caso di {\em problemi
discreti ai valori ai limiti}, cio\`e ottenuti associando
all'equazione (\ref{omo}) un insieme di condizioni ai limiti,
invece che ai valori iniziali, come fatto in (\ref{ini}). Si veda,
ad esempio, \cite{BrTr}.\end{oss}

\begin{eser} La {\em successione di Fibonacci,}
$$1,~1,~2,~3,~5,~8,~13,~\dots,$$
soddisfa l'equazione alle differenze
$$F_0=F_1=1, \qquad F_{n+1} = F_n+F_{n-1}, \qquad n = 1,2,\dots.$$
Determinare l'espressione del generico elemento $F_n$ della
successione, ed il limite del rapporto
$$\frac{F_{n+1}}{F_n},\qquad n\rightarrow\infty.$$
\end{eser}

\section{Polinomi di Chebyshev}\label{cebypol}

Gli argomenti che abbiamo appena visto riguardo alle equazioni
alle differenze, trovano una immediata applicazione nella
discussione riguardante famiglie di polinomi definite mediante una
equazione di ricorrenza. Nello specifico, siamo interessati alla
famiglia di polinomi $\{p_k(z)\}$ che soddisfano alla seguente
equazione:
\begin{equation}\label{ceby}
p_{k+1}(z)=2z\,p_k(z)-p_{k-1}(z),\qquad z\in\CC.
\end{equation}

\no In particolare, considereremo le successioni definite dalle
seguenti condizioni iniziali:

\begin{enumerate}

\item $p_0(z)\equiv 1$, $p_1(z)=z$, che origina i {\em polinomi di
Chebyshev di prima specie}, e che denoteremo con $\{T_k(z)\}$;

\item $p_{-1}(z)\equiv0$, $p_0(z)\equiv1$, che origina i {\em polinomi
di Chebyshev di seconda specie}, e che denoteremo con
$\{U_k(z)\}$.

\end{enumerate}

Vediamo di derivare alcune propriet\`a dei polinomi di Chebyshev,
in particolare di quelli di prima specie. Osserviamo che, per un
fissato $z\in\CC$, la successione $\{T_k(z)\}$ soddisfa
l'equazione alle differenze
\begin{equation}\label{cebye}
y_{k+1} = 2zy_k-y_{k-1}, \qquad k=1,2,\dots,
\end{equation}

\no con le condizioni iniziali
\begin{equation}\label{ceby1}
y_0 = 1, \qquad y_1=z.
\end{equation}

\no Le radici del polinomio caratteristico associato alla
(\ref{cebye}), $$\phi(w) = w^2-2zw+1,$$ sono $$w_{1/2} =
z\pm\sqrt{z^2-1}.$$ Osserviamo che $w_1w_2=1$ e, pertanto,
$w_2=w_1^{-1}$. Ricercando, quindi, la soluzione nella forma
$y_k=c_1w_1^k+c_2w_1^{-k}$, e imponendo le condizioni iniziali
(\ref{ceby1}), si ottiene, infine,

\begin{eqnarray}\nonumber
T_k(z) &=& \frac{1}2\left(w_1^k+w_1^{-k}\right)
    =\frac{1}2\left( e^{k\log w_1} +e^{-k\log w_1} \right)\\
    \label{genceb}
    &=& \cosh k\log(z+\sqrt{z^2-1}).
\end{eqnarray}

Dalla (\ref{genceb}) si ottengono agevolmente un certo numero di
propriet\`a, molte delle quali saranno utili nei prossimi
capitoli.

Se $z\in[-1,1]$, possiamo porre $z=\cos\theta$, con
$\theta\in[0,\pi]$. Ne consegue che $w_1 = \cos\theta
+i\sin\theta=e^{i\theta}$. Pertanto:

\begin{description}

\item[C1:~] $T_k(z)=T_k(\cos\theta) = \cosh k\log e^{i\theta}=\cos
k\theta;$

\item[C2:~] gli zeri di $T_k(z)$ sono dati da $z_j^{(k)} = \cos
\frac{(2j+1)\pi}{2k}, ~j=0,\dots,k-1;$

\item[C3:~] $\max_{z\in[-1,1]} |T_k(z)| =
1$, ~ $T_k\left(\cos\frac{j}k\pi\right)=(-1)^j$,
$j=0,\dots,k$.

\end{description}

Inoltre, dalla equazione di ricorrenza (\ref{ceby}), si ottengono
facilmente anche le seguenti altre propriet\`a:
\begin{description}

\item[C4:~] $T_k(z)$ ha grado esatto $k$ ed il suo coefficiente
principale \`e $2^{k-1};$

\item[C5:~] la famiglia di polinomi $\{\hT_k\}$, con
$\hT_k(z)\equiv2^{1-k}T_k(z)$, \`e una famiglia di polinomi
monici, con $\hT_k$ avente grado esatto $k$.

\end{description}

La famiglia di polinomi monici $\{\hT_k\}$ \`e molto importante in
Teoria dell'Approssimazione, in quanto fornisce la soluzione di
alcuni importanti problemi di {\em minimassimo}. Vale, infatti, la seguente
propriet\`a:
\begin{description}

\item[C6:~] $\max_{z\in[-1,1]}|\hT_k| = 2^{1-k} \equiv
\min_{p\in\Pi_k'}\max_{z\in[-1,1]}|p(z)|$, dove $\Pi_k'$ \`e l'insieme dei
polinomi monici di grado $k$.

\end{description}

\no Inoltre, fissato $\bz>1$, valgono le seguenti propriet\`a:
\begin{description}

\item[C7:~] $T_k(\bz) > \frac{1}2\left(\bz+\sqrt{\bz^2-1}\right)^k;$

\item[C8:~] il polinomio di grado $k$, definito da $$p_k(z;\bz) \equiv
\frac{T_k(z)}{T_k(\bz)},$$ \`e soluzione del problema di
minimassimo $$\min_{p\in\Pi_k,\, p(\bz)=1}\, \max_{z\in[-1,1]}\,
|p(z)|.$$

\end{description}

\medskip
\begin{eser} Dimostrare le propriet\`a C2--C8.\end{eser}

\medskip
\begin{eser} Per un assegnato $\bz>1$, calcolare il limite del
rapporto $$\frac{T_{k+1}(\bz)}{T_k(\bz)},\qquad k\rightarrow\infty.$$\end{eser}

\medskip
\begin{eser} Per un assegnato $\hat{z}<-1$, calcolare il limite del
rapporto $$\frac{T_{k+1}(\hat{z})}{T_k(\hat{z})},\qquad k\rightarrow\infty.$$\end{eser}

%
%
\chapter{Metodi iterativi di base}\label{cap3}

{\em
In questo capitolo studieremo i metodi iterativi di base che,
anche storicamente, sono stati i primi ad essere ottenuti. Sono
tutti metodi caratterizzati da una lodevole semplicit\`a
concettuale ma, in genere, hanno prestazioni inferiori, rispetto ai metodi
che esamineremo nei capitoli successivi.}

\section{I metodi di Jacobi, Gauss-Seidel e SOR}\label{base}

Avendo in mente le caratteristiche delle matrici viste
nell'esempio del primo capitolo, ovvero delle $M$-matrici, sar\`a
conveniente partizionare la matrice dei coefficienti del problema,
\begin{equation}\label{sist}
A\bfx=\bfb,
\end{equation}

\no nella forma
\begin{equation}\label{LDU}
A = D-L-U,
\end{equation}

\no dove la matrice $D$ \`e diagonale, $L$ \`e strettamente
triangolare inferiore ed $U$ \`e strettamente triangolare
superiore.

\begin{oss} Nel caso in cui $A$ sia una $M$-matrice, allora $D,L,U\ge0$.
\end{oss}

Tenendo conto delle (\ref{sist})-(\ref{LDU}), si ottiene che la
soluzione $\bfx^*$ del sistema lineare soddisfer\`a anche
all'equazione
\begin{equation}\label{Jacob1}
D\bfx^* = (L+U)\bfx^* +\bfb.
\end{equation}

\no Da questa, fissato un vettore iniziale $\bfx_0$, si pu\`o
pensare di definire l'iterazione
\begin{equation}\label{Jacob2}
D\bfx_{k+1} = (L+U)\bfx_k +\bfb, \qquad k=0,1,2,\dots,
\end{equation}

\no che descrive il {\em metodo di Jacobi}. In tal caso, il nuovo
punto $\bfx_{k+1}$ \`e definito formalmente da $$\bfx_{k+1} =
D^{-1}(L+U)\bfx_k +D^{-1}\bfb, \qquad k=0,1,2,\dots,$$

\no ovvero si ottiene risolvendo un sistema diagonale. Ne consegue
che la implementazione del metodo di Jacobi risulta essere
particolarmente semplice, e richiede che gli elementi diagonali di
$A$ siano non nulli.

Chiaramente noi siamo interessati al caso in cui la procedura
(\ref{Jacob2}) sia {\em convergente}, ovvero quando $$\bfx_k
\rightarrow\bfx^*, \mbox{\qquad per \quad} k\rightarrow\infty.$$

\no Per stabilire quando questo avvenga, si studia l'{\em equazione
dell'errore}, ottenuta sottraendo la (\ref{Jacob2}) dalla
(\ref{Jacob1}). In tal caso, denotando l'errore al passo $k$ col
vettore
$$\bfe_k = \bfx^* -\bfx_k,$$

\no si ottiene l'equazione dell'errore
$$ \bfe_{k+1} = D^{-1}(L+U)\bfe_k, \qquad k=0,1,2,\dots,$$

\no la cui soluzione si vede facilmente essere data da
$$ \bfe_k = \left( D^{-1}(L+U) \right)^k \bfe_0, \qquad k =
0,1,2,\dots.$$

\no Chiaramente, il metodo sar\`a convergente se e solo se
$$\bfe_k \rightarrow\bfo,\mbox{\qquad per \quad}
k\rightarrow\infty.$$

\no Se ne conclude che il metodo di Jacobi sar\`a convergente se e
solo se la corrispondente {\em matrice di iterazione},
\begin{equation}\label{GJ}
G_J = D^{-1}(L+U),
\end{equation}

\no \`e una matrice convergente, ovvero se (vedi Teorema~\ref{conv})
\begin{equation}\label{Jacob4}
\rho(D^{-1}(L+U))<1.
\end{equation}

Vedremo nel seguito delle condizioni di tipo generale per la
matrice $A$, sotto le quali la propriet\`a (\ref{Jacob4}) \`e
garantita. Per il momento, ci limitiamo ad enunciare particolari
propriet\`a sufficienti per garantire la convergenza del metodo di
Jacobi.

\begin{defi} Se $A=(a_{ij})\in\RR^{n\times n}$, allora $A$ si
dir\`a:

\begin{itemize}

\item {\em diagonale dominante per righe}, se $|a_{ii}|>\sum_{j\ne
i} |a_{ij}|,~i=1,\dots,n$;

\item {\em diagonale dominante per colonne}, se $|a_{ii}|>\sum_{j\ne
i} |a_{ji}|,~i=1,\dots,n$.

\end{itemize}
\end{defi}

Si lascia come semplice esercizio la dimostrazione della seguente
propriet\`a (vedi (\ref{GJ})).

\begin{teo} Se $A$ \`e diagonale dominante, per righe o per
colonne, allora $\rho(G_J)<1$.\end{teo}

\begin{oss}\label{stopjac} Nella implementazione pratica dei metodi iterativi
(come vedremo in modo pi\`u approfondito in Sezione~\ref{stopcg}), un criterio
di arresto molto utilizzato \`e basato sul controllo della norma del {\em
residuo} $$\bfr_k = \bfb - A\bfx_k.$$ In tal caso, la nuova approssimazione
$\bfx_{k+1}$ fornita al metodo di Jacobi (\ref{Jacob2}), si vede essere
ottenibile come (vedi (\ref{LDU})) $$\bfx_{k+1} = \bfx_k +
D^{-1}\bfr_k.$$\end{oss}

\begin{eser} Scrivere un codice che implementi il metodo di
Jacobi applicato al problema (\ref{test}) in Esercizio~\ref{ex1}.
\end{eser}

\bigskip
L'equazione (\ref{Jacob1}) pu\`o anche essere riscritta come
$$(D-L)\bfx^* = U\bfx^* +\bfb,$$

\no inducendo lo schema iterativo
\begin{equation}\label{GS2}
(D-L)\bfx_{k+1} = U\bfx_k +\bfb, \qquad k=0,1,2,\dots,
\end{equation}

\no che definisce il {\em metodo di Gauss-Seidel}. In tal caso, ad
ogni passo si deve risolvere un sistema triangolare. Anche in
questo caso, si richiede che la matrice $D$ sia nonsingolare.

Procedendo come nel caso del metodo di Jacobi, l'equazione
dell'errore si vede essere
$$\bfe_{k+1} = (D-L)^{-1}U\bfe_k,$$

\no per cui il metodo sar\`a convergente se e solo se la
corrispondente matrice di iterazione,
\begin{equation}\label{Ggs}
G_{GS} = (D-L)^{-1}U,
\end{equation}

\no ha raggio spettrale minore di 1. Riguardo alle caratteristiche
della matrice $A$ che garantiscono questa propriet\`a, si rimanda
all'analisi, pi\`u generale, che vedremo in Sezione~\ref{splreg}.
Per il momento, ci limitiamo ad enunciare un risultato parziale.

\begin{teo} Se la matrice $A$ \`e sdp, il metodo di Gauss-Seidel
\`e convergente.\end{teo}

\proof Se $A$ \`e sdp, allora (vedi (\ref{LDU})), $U=L^\top $.
Cominciamo con lo stabilire che $\lam=1$ non \`e autovalore di
$G_{GS}$. Infatti, se cos\`\i\, fosse, il corrispondente
autovettore, sia esso $\bfu$ sarebbe tale che $A\bfu=\bfo$ e,
pertanto, $A$ sarebbe singolare. D'altronde, se
$\lam\in\sigma(G_{GS})$, allora, detto $\bfv$ il corrispondente
autovettore, si ha:
\begin{equation}\label{gst1}
L^\top \bfv = \lam (D-L)\bfv.
\end{equation}

\no Da questa si ottiene che
$$
A\bfv \equiv (D-L-L^\top )\bfv = (1-\lam)(D-L)\bfv.
$$

\no Pertanto,
\begin{equation}\label{gst2}
\bfv^* A \bfv = (1-\lam) \bfv^* (D-L)\bfv =
(1-\bar\lam)\bfv^*(D-L)^\top \bfv = \bfv^* A^\top \bfv.
\end{equation}

\no Dalla (\ref{gst1}) si ottiene, quindi,
\begin{eqnarray*}
(1-\bar\lam)\bfv^*(D-L)^\top \bfv &=& (1-\lam) \bfv^*(D-L)\bfv =
(1-\lam)\left( \bfv^*D\bfv -\bfv^*L\bfv\right) \\
&=& (1-\lam)\left( \bfv^*D\bfv -\bar\lam\bfv^*(D-L)^\top \bfv\right).
\end{eqnarray*}

\no Segue che
$$(1-\bar\lam)\bfv^*(D-L)^\top \bfv = (|\lam|^2
-\bar\lam)\bfv^*(D-L)^\top \bfv +(1-\lam)\bfv^*D\bfv,$$

\no da cui
$$(1-\lam)\bfv^*D\bfv = (1-|\lam|^2)\bfv^*(D-L)^\top \bfv.$$

\no Moltiplicando per $(1-\bar\lam)$, e tenendo conto della
(\ref{gst2}),
$$ |1-\lam|^2\bfv^*D\bfv =
(1-|\lam|^2)(1-\bar\lam)\bfv^*(D-L)^\top \bfv =
(1-|\lam|^2)\bfv^*A\bfv.$$

\no Considerando che $A$ e $D$ sono sdp, e che $\lam=1$ non \`e
autovalore, segue infine che
$$1-|\lam|^2 = |1-\lam|^2\frac{\bfv^*D\bfv}{\bfv^*A\bfv}>0,$$

\no ovvero $|\lam|<1$.\QED

\begin{eser} Dimostrare che se $A$ \`e diagonale dominante, per
righe o per colonne, allora $\rho(G_{GS})<1$\quad (suggerimento: $(D-L)^{-1}U\sim
U(D-L)^{-1}$).\end{eser}

\begin{oss}\label{stopgs} In analogia con quanto detto nella
Osservazione~\ref{stopjac} per il metodo di Jacobi, anche per il metodo di
Gauss-Seidel, definito il residuo al passo $k$ come $$\bfr_k = \bfb -
A\bfx_k,$$ la nuova approssimazione (\ref{GS2}) \`e ottenibile risolvendo
$$(D-L)\bfu_k = \bfr_k, \qquad \bfx_{k+1}= \bfx_k +\bfu_k.$$
\end{oss}

\begin{eser} Scrivere un codice che implementi il metodo di
Gauss-Seidel applicato al problema (\ref{test}) in
Esercizio~\ref{ex1}.
\end{eser}

Una variante del metodo di Gauss-Seidel si ottiene definendo la
generica componente della nuova approssimazione come combinazione,
a pesi $\om$ e $1-\om$ rispettivamente, di quella che fornirebbe
Gauss-Seidel e di quella precedente. Poich\`e la matrice $L$ nella
(\ref{GS2}) \`e strettamente triangolare inferiore, questo
equivale, formalmente, a scrivere $$\bfx_{k+1} = \om D^{-1}\left(
L\bfx_{k+1} +U\bfx_k +\bfb \right) + (1-\om)\bfx_k, \qquad
k=0,1,2,\dots,$$

\no ovvero, $$ (D-\om L)\bfx_{k+1} = ( (1-\om)D+\om U)\bfx_k
+\om\bfb,\qquad k=0,1,2,\dots.$$

\no Si ottiene, in questo modo, il metodo di {\em
sovrarilassamento}, denominato {\em SOR} (Successive Over
Relaxation), che a volte scriveremo anche come $$
\frac{1}\om(D-\om L)\bfx_{k+1} = \frac{1}\om( (1-\om)D+\om
U)\bfx_k +\bfb,\qquad k=0,1,2\dots. $$

Osserviamo che l'ultima equazione \`e ben definita, in quanto
$\om=0$ non \`e un valore ammesso (porterebbe a
$\bfx_{k+1}=\bfx_k$). Inoltre, per $\om=1$, si riottiene il metodo
di Gauss-Seidel (\ref{GS2}). Pi\`u in generale, quando $0<\om<1$,
i due pesi sono entrambi positivi, ed il nuovo punto risulta
essere una combinazione convessa di quello di Gauss-Seidel e di
quello precedente. In tal caso si parla di {\em
sottorilassamento}. Invece, per $\om>1$ si parla di {\em
sovrarilassamento}.

Infine, osserviamo che l'iterazione risulta essere definita sotto
la solita condizione che $D$ sia nonsingolare.

Per quanto riguarda la convergenza del metodo SOR, un ragionamento
del tutto simile a quello fatto per gli altri metodi porta
all'equazione dell'errore $$\bfe_{k+1} = (D-\om
L)^{-1}((1-\om)D+\om U)\bfe_k, \qquad k=0,1,2\dots,$$

\no e, quindi, il metodo sar\`a convergente se e solo se la
matrice di iterazione,
\begin{equation}\label{Gw}
G_\om = (D-\om L)^{-1}((1-\om)D+\om U),
\end{equation}

\no \`e convergente, ovvero, se il suo raggio spettrale, che
denoteremo con $\rho_\om$, \`e minore di 1. Anche la convergenza
di questo metodo sar\`a discussa nella prossima sezione. Si
premette solo qualche risultato preliminare, che fornisce delle
limitazioni sul {\em parametro di rilassamento} $\om$.

\begin{teo} $\rho_\om\ge|1-\om|$.\end{teo}

\proof Infatti, se $A\in\RR^{n\times n}$, dalla (\ref{Gw}) si ha,
tenendo conto della triangolarit\`a dei suoi fattori,
\begin{eqnarray*}
|\det(G_\om)| &=& |\det( I-\om D^{-1} L)^{-1}\det((1-\om)I+\om
D^{-1}U)|\\ &=& |1\cdot (1-\om)^n| = |1-\om|^n \equiv
\prod_{\lam\in\sigma(G_\om)} |\lam|\le \rho_\om^n.
\end{eqnarray*}\QED

\begin{cor}\label{w02} $\rho_\om<1 \Rightarrow \om\in(0,2)$.\end{cor}

\begin{oss} Ovviamente, la scelta pi\`u appropriata per il
parametro di rilassamento $\om$ sarebbe quella per la quale
$\rho_\om$ \`e minimizzato. Non approfondiremo oltre questo
aspetto, per il quale si rimanda alla estensiva trattazione in
\cite{Yo}.\end{oss}

\begin{eser} Scrivere un codice che implementi il metodo SOR
applicato al problema (\ref{test}) in Esercizio~\ref{ex1} e
confrontare i valori $\om=0.5,1,1.5$, per diversi valori di $N$.
\end{eser}

\section{Splitting regolari di matrici}\label{splreg}

I metodi visti fino ad ora per risolvere il sistema lineare
(\ref{sist}), sono tutti del tipo
\begin{equation}\label{split1}
M\bfx_{k+1} = N\bfx_k + \bfb, \qquad k=0,1,2,\dots,
\end{equation}

\no dove
\begin{equation}\label{split2}
A=M-N, \qquad \det(M)\ne0.
\end{equation}

\begin{defi} Due matrici $M$ ed $N$ che soddisfino alle
propriet\`a (\ref{split2}), definiscono uno {\em splitting} della
matrice $A$. Se, inoltre,
$$ M^{-1}\ge0, \qquad N\ge0, $$

\no allora lo splitting si dice {\em regolare} per tale matrice.
\end{defi}

Ad esempio:

\begin{itemize}

\item $M_J=D$ e $N_J=L+U$ definisce il metodo di Jacobi;

\item $M_{GS}=D-L$ e $N_{GS}=U$ definisce il metodo di Gauss-Seidel;

\item $M_\om=\om^{-1}(D-\om L)$ e $N_\om=\om^{-1}((1-\om)D+\om U)$
definisce il metodo SOR con parametro di rilassamento $\om$.

\end{itemize}

\begin{eser}
Dimostrare che, nel caso $D,L,U\ge0$,

\begin{itemize}

\item gli splitting di Jacobi e Gauss-Seidel sono regolari;

\item lo splitting di SOR \`e regolare per $0<\om\le1$.

\end{itemize}
\end{eser}

\begin{oss}\label{stopsr} Considerazioni simili a quelle viste nelle
Osservazioni~\ref{stopjac} e \ref{stopgs} per i metodi di Jacobi e
Gauss-Seidel, possono essere fatte per il generico metodo iterativo
(\ref{split1}): se $$\bfr_k=\bfb-A\bfx_k,$$ \`e il residuo al passo $k$,
allora dalla (\ref{split2}) segue che $$M\bfu_k = \bfr_k, \qquad \bfx_{k+1} =
\bfx_k+\bfu_k.$$\end{oss}

Le propriet\`a degli splitting regolari hanno importanti ricadute
sulla convergenza del metodo iterativo (\ref{split1}), che risulta
essere convergente se e solo se
\begin{equation}\label{rosplit}
\rho(M^{-1}N)<1.
\end{equation}

\no Vale, infatti, il seguente risultato.

\begin{teo} Se lo splitting (\ref{split2}) \`e regolare, e
$A^{-1}\ge0$, allora vale (\ref{rosplit}), ovvero il metodo
iterativo (\ref{split1}) \`e convergente.\end{teo}

\proof Cominciamo con l'osservare che, poich\'e lo splitting \`e
regolare,
$$0\le M^{-1}N = (A+N)^{-1}N = (I+A^{-1}N)^{-1} A^{-1}N.$$

\no Ora, se $\tau\in\sigma(A^{-1}N)$ e $\bfv$ \`e il
corrispondente autovettore, allora
\begin{eqnarray*}
A^{-1}N\bfv = \tau\bfv, &\Leftrightarrow& (I+A^{-1}N)^{-1}\bfv =
(1+\tau)^{-1}\bfv\\ &\Leftrightarrow& (I+A^{-1}N)^{-1}A^{-1}N\bfv
=  \tau(1+\tau)^{-1}\bfv,
\end{eqnarray*}

\no ovvero $$\lam = \frac{\tau}{1+\tau}$$ \`e autovalore di
$M^{-1}N$, con autovettore $\bfv$. Osserviamo che, essendo
$M=(A+N)$ nonsingolare per ipotesi, allora $\tau$ non pu\`o
assumere il valore $-1$. Inoltre, dal Teorema~\ref{PF}, segue che
l'autovalore di massimo modulo, sia esso $\lam$, \`e reale e non
negativo, con un autovettore $\bfv\ge0$. Per quanto visto innanzi,
ad esso corrisponder\`a un autovalore $\tau$ di $A^{-1}N$, che \`e
una matrice non negativa, per l'ipotesi $A^{-1}\ge0$. Poich\`e
$\lam\ge0$, segue che $\tau\ge0$ o $\tau<-1$. L'ultima
possibilit\`a \`e per\`o preclusa dal fatto che il corrispondente
autovettore di $M^{-1}N$, $\bfv\ge0$, \`e anche autovettore di
$A^{-1}N$. Infatti, in caso contrario, si avrebbe:
$$0\le A^{-1}N\bfv =\tau\bfv \le -\bfv \le 0,$$

\no dove la prima diseguaglianza segue dal fatto che $A^{-1}N\ge0$
e $\bfv\ge0$. Pertanto, all'autovalore $\lambda$ corrisponder\`a
un autovalore $\tau\ge0$. La tesi si completa osservando che la
funzione
$$f(\tau) = \frac{\tau}{1+\tau}$$

\no \`e, per $\tau\ge0$, crescente e strettamente minore di 1.\QED
\vspace{2mm}

Valgono i seguenti risultati, la cui dimostrazione \`e lasciata
come esercizio.

\begin{cor}\label{Msplit} Se lo splitting (\ref{split2}) \`e regolare, e
$A$ \`e una $M$-matrice, allora il metodo iterativo (\ref{split1})
\`e convergente.\end{cor}

\begin{cor} Se $A$ \`e una $M$-matrice e la matrice $M$ dello
splitting (\ref{split2}) \`e ottenuta ponendo a zero qualche
elemento extra-diagonale, allora il metodo iterativo
(\ref{split1}) \`e convergente.\end{cor}

Vediamo ora alcuni risultati che consentono di comparare la
velocit\`a di convergenza di metodi iterativi definiti da
splitting regolari.

\begin{teo}\label{comsplit} Siano $A=M_1-N_1=M_2-N_2$ due
splitting regolari di $A$, matrice tale che $A^{-1}\ge0$. Se
$N_1\le N_2$, allora
$\rho(M_1^{-1}N_1)\le\rho(M_2^{-1}N_2).$\end{teo}

\proof Per ipotesi, $0\le N_1\le N_2$ e $A^{-1}\ge0$. Ne consegue
che $0\le A^{-1}N_1\le A^{-1}N_2$ e, quindi, dal
Teorema~\ref{monot}, $\rho(A^{-1}N_1)\le\rho(A^{-1}N_2)$. La tesi
discende infine dal fatto che
$$\rho(M_1^{-1}N_1) = \frac{\rho(A^{-1}N_1)}{1+\rho(A^{-1}N_1)}
\le \frac{\rho(A^{-1}N_2)}{1+\rho(A^{-1}N_2)} =
\rho(M_2^{-1}N_2).$$\QED \vspace{2mm}

Da questo risultato si possono trarre alcune significative
conclusioni.

\begin{cor} Se $A$ \`e una $M$-matrice, allora (vedi (\ref{GJ}),
(\ref{Ggs}) e (\ref{Gw})),

\begin{itemize}

\item $\rho(G_{GS})\le\rho(G_J);$

\item se $0\le \om_1 \le \om_2 \le 1: \rho(G_{GS}) \le \rho(G_{\om_2}) \le \rho(G_{\om_1})<1.$

\end{itemize} \end{cor}

\proof Infatti, entrambi gli splitting al primo punto risultano
essere regolari e (vedi (\ref{LDU})) $N_J = L+U\ge U=N_{GS}\ge0$.

Il secondo punto discende, analogamente, dal fatto che $$N_{\om_1}
\equiv \frac{1-\om_1}{\om_1}\, D +U \ge N_{\om_2} \equiv
\frac{1-\om_2}{\om_2}\, D +U \ge N_{GS} \equiv U.$$\QED

\begin{oss}\label{confr}
Da questo ultimo risultato, segue che, nel caso di $M$-matrici, il
metodo di Gauss-Seidel converge pi\`u velocemente, in genere, di
quello di Jacobi, come anche di SOR, se il parametro di
rilassamento \`e pi\`u piccolo di 1. Se ne conclude che,
quest'ultimo, dovrebbe essere scelto (vedi Corollario~\ref{w02})
nell'intervallo [1,2).

Ad esempio (vedi \cite{Va,Yo}), nel caso della matrice
dell'equazione (\ref{test}) in Esercizio~\ref{ex1}, si pu\`o
dimostrare che

\begin{enumerate}

\item $\rho(G_{GS}) = \rho(G_J)^2$;

\item il parametro di rilassamento che minimizza $\rho(G_\om)$ \`e
dato da $$\om^* = 1 + \left(
\frac{\rho(G_J)}{1+\sqrt{1-\rho(G_{GS})}}\right)^2.$$
Evidentemente, $\rho(G_{\om^*})\le \rho(G_{GS})$.
\end{enumerate}

\begin{figure}[t]
\begin{center}
\includegraphics[width=13cm,height=10cm]{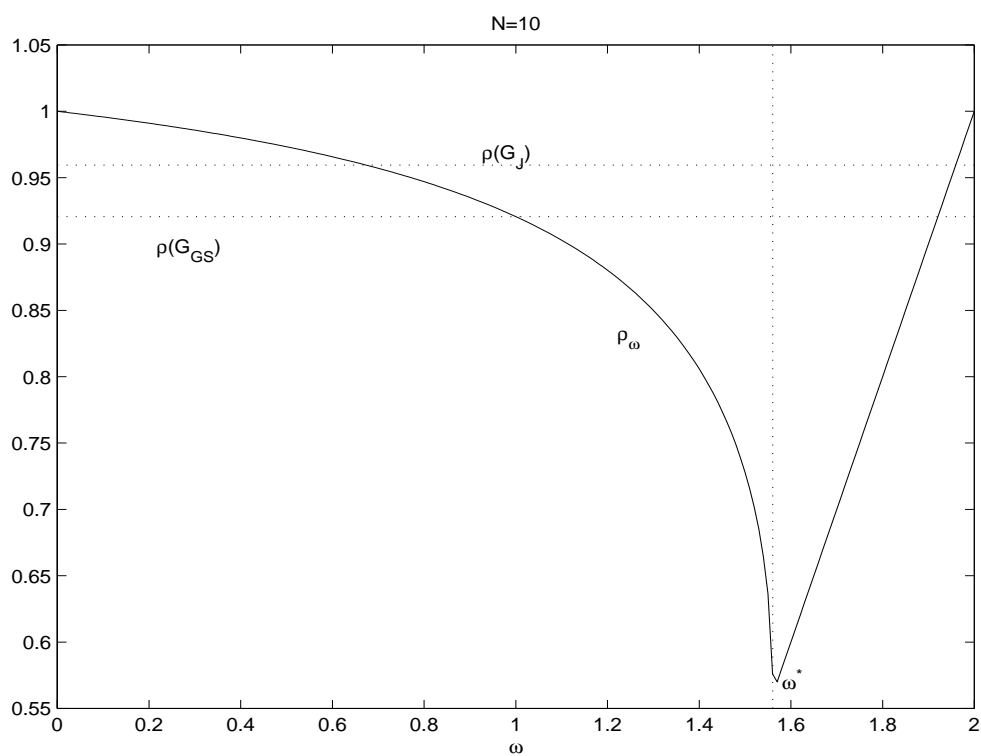}
\caption{\protect\label{fig2} Raggio spettrale di SOR in funzione
del parametro di rilassamento $\om$, relativo al problema
(\ref{test}) dell'Esercizio~\ref{ex1}, $N=10$.}
\end{center}
\end{figure}

\no Come esempio illustrativo, in Figura~\ref{fig2} sono riportati
i grafici dei raggi spettrali dei metodi di Jacobi, Gauss-Seidel,
e di SOR in funzione del parametro di rilassamento $\om$, relativi
al problema su menzionato, per $N=10$. In tal caso, si ottiene
che:
\begin{equation}\label{punt}
\rho(G_J) \approx 0.96, \qquad \rho(G_{GS}) \approx 0.92, \qquad
\om^* \approx 1.56, \qquad \rho(G_{\om^*}) \approx 0.57.
\end{equation}

\no Se ne conclude che il potenziale guadagno nell'utilizzo di SOR
\`e notevole, a patto di conoscere il parametro di rilassamento
ottimale.

\end{oss}

\begin{eser} Costruire una procedura che calcoli i dati per
ottenere la Figura~\ref{fig2}.\end{eser}

\begin{eser} Sempre con riferimento ai metodi descritti in Figura~\ref{fig2},
si supponga di utilizzarli per risolvere il problema (\ref{test})
partendo da un punto ini\-ziale cui corrisponde un errore di norma
uguale ad 1. Tenendo conto del fatto che una iterata del metodo di
Gauss-Seidel ha praticamente lo stesso costo di quella di Jacobi,
che sono approssimativamente i 5/6 del co\-sto dell'ite\-razione
di SOR, qual'\`e il potenziale guadagno (stimato), rispetto al
metodo di Jacobi, dei metodi di Gauss-Seidel e di SOR con
parametro ottimale $\om^*$, se $N=10$ e si desidera abbattere
l'errore di un fattore $10^6$?\end{eser}

\section{Generalizzazione a blocchi dei metodi}\label{blocks}

Il precedente criterio di confronto per splitting regolari, dato
dal risultato del Teorema~\ref{comsplit}, fornisce il supporto
teorico che giustifica la estensione a blocchi dei metodi su
esaminati. Vedremo i dettagli per i metodi di Jacobi e
Gauss-Seidel, in quanto l'estensione ad SOR \`e del tutto simile.

Quanto andiamo a dire si applica in modo naturale a $M$-matrici
con una struttura a blocchi della forma vista in (\ref{A}).
Infatti, il partizionamento (\ref{LDU}) da origine agli splitting
regolari
$$M_J = D, \quad N_J = L+U, \qquad M_{GS} = D-L, \quad N_{GS} =
U,$$

\no cui corrispondono i metodi di Jacobi e Gauss-Seidel,
rispettivamente, che sono, in questo caso, convergenti (vedi
Corollario~\ref{Msplit}). Si pu\`o, tuttavia, considerare anche il
partizionamento a blocchi
$$A = \hat{D} -\hat{L} -\hat{U},$$

\no che, nel caso della matrice (\ref{A}), \`e definito da
$$\hat{D} = \pmatrix{cccc} A_1 \\ &A_2\\&&\ddots \\
&&&A_N\endpmatrix, \qquad \hat{L} = \hat{U}^\top  = \pmatrix{cccc} O\\
I_N &O\\ &\ddots &\ddots \\ &&I_N &O\endpmatrix.$$

\no Questo induce naturalmente gli splitting
\begin{equation}\label{blockJGS}
\hat{M}_J = \hat{D}, \quad \hat{N}_J = \hat{L}+\hat{U}, \qquad
\hat{M}_{GS} = \hat{D}-\hat{L}, \quad \hat{N}_{GS} = \hat{U},
\end{equation}

\no che definiscono le {\em versioni a blocchi} dei metodi di
Jacobi e Gauss-Seidel, rispettivamente.

\begin{eser} Dimostrare che, nel caso in cui $A$ sia una
$M$-matrice, allora gli splitting (\ref{blockJGS}) sono regolari
e, pertanto, i corrispondenti metodi iterativi sono convergenti.
\end{eser}

A questo punto, dal Teorema~\ref{comsplit} segue facilmente che,
nel caso di $M$-matrici,
\begin{equation}\label{block1}
\rho(\hat{G}_J)\le \rho(G_J), \qquad
\rho(\hat{G}_{GS})\le\rho(G_{GS}),
\end{equation}

\no dove $$\hat{G}_J = \hat{M}_J^{-1}\hat{N}_J, \qquad
\hat{G}_{GS} = \hat{M}_{GS}^{-1}\hat{N}_{GS},$$

\no sono le matrici di iterazione dei metodi di Jacobi e
Gauss-Seidel a blocchi, rispettivamente. Come esempio, nel caso
della matrice esaminata nella Osservazione~\ref{confr}, si ottiene
(confrontare con (\ref{punt})) $$\rho(\hat{G}_J) \approx 0.92,
\qquad \rho(\hat{G}_{GS}) \approx 0.85.$$

\begin{eser} Dimostrare le diseguaglianze
(\ref{block1}).\end{eser}

\section{Il metodo semi-iterativo di Chebyshev}\label{s-ceby}

Si supponga di aver generato approssimazioni
\begin{equation}\label{xk}
\bfx_0,\bfx_1,\dots,\bfx_j,
\end{equation}

\no della soluzione, sia essa $\bfx^*$, dell'equazione
(\ref{sist}), mediante il metodo iterativo (\ref{split1}), la cui
matrice di iterazione, $$G=M^{-1}N,$$ supponiamo essere simmetrica
e convergente. Pertanto, per una data matrice ortogonale $V$,
\begin{equation}\label{VLV}
G = V\Lambda V^\top , \qquad \Lambda =
\pmatrix{ccc}\lam_1\\&\ddots\\&&\lam_n\endpmatrix,
\end{equation}

\no dove i suoi autovalori saranno
\begin{equation}\label{albe}
-1<\aa\le\lam_1\le\lam_2\le\dots\le\lam_n\le\bb<1. \end{equation}

\no Per il momento, supporremo che gli scalari $-1<\aa<\bb<1$
siano noti.

Vogliamo ora determinare dei coefficienti $\{c_{kj}\}$ tali che il
vettore
\begin{equation}\label{xceby}
\bfx^{(j)} = \sum_{k=0}^j c_{kj} \bfx_k
\end{equation}

\no sia una approssimazione della soluzione $\bfx^*$ migliore di
quelle iniziali (\ref{xk}). Poich\`e nel caso in cui
$\bfx_k=\bfx^*, k=0,\dots,j,$ \`e lecito aspettarsi
$\bfx^{(j)}=\bfx^*$, dalla (\ref{xceby}) segue che i coefficienti
$\{c_{kj}\}$ devono soddisfare al vincolo
\begin{equation}\label{sum1}
\sum_{k=0}^j c_{kj} = 1.
\end{equation}

\no Fissato questo vincolo, scegliamo i coefficienti $\{c_{kj}\}$
in modo da minimizzare, in qualche senso, l'errore $\bfe^{(j)}
=\bfx^*- \bfx^{(j)}$. Osserviamo che, per la (\ref{sum1}), si ha
$$
\bfe^{(j)} = \bfx^*-\bfx^{(j)} = \bfx^*-\sum_{k=0}^j
c_{kj}\bfx_k =\sum_{k=0}^j c_{kj}(\bfx^* - \bfx_k) \equiv
\sum_{k=0}^j c_{kj}\bfe_k.
$$

\no Inoltre, poich\`e le approssimazioni (\ref{xk}) sono state
ottenute dall'ite\-razione (\ref{split1}), segue che gli errori
soddisfano all'equazione $$\bfe_{k+1} = G\bfe_k \equiv M^{-1}N
\bfe_k,\qquad k=0,1,2,\dots,$$ e, pertanto,
$$ \bfe^{(j)} = \sum_{k=0}^j c_{kj}G^k\bfe_0 \equiv p_j(G)\bfe_0,
$$

\no dove $p_j\in\Pi_j$ \`e un polinomio che deve anche soddisfare
al vincolo (vedi (\ref{sum1}))
$$ p_j(1) = \sum_{k=0}^j c_{kj} = 1. $$

\no Osservando che $$\|\bfe^{(j)}\| = \| p_j(G)\bfe_0 \|\le
\|p_j(G)\|\cdot \|\bfe_0\|,$$

\no al fine di rendere la soluzione indipendente dalla scelta del
vettore iniziale, e quindi da $\bfe_0$, potremmo ricercare $p_j$
come soluzione del seguente problema di minimo, $$p_j =
\arg\min_{p\in\Pi_j,\,p(1)=1} \|p(G)\|.$$

\no Tuttavia, se la norma utilizzata \`e la norma $\|\cdot\|_2$,
dalla (\ref{VLV}) segue che:
\begin{equation}\label{pG}
\|p(G)\| = \|p(V\Lambda V^\top )\| = \|V p(\Lambda) V^\top \| =
\|p(\Lambda)\| = \max_{\lam\in\sigma(G)} | p(\lam) |.
\end{equation}

\no Questo presuppone per\`o la conoscenza dello spettro di $G$,
cosa assai impro\-babile. Quello che noi supporremo di conoscere
\`e un intervallo di inclusione, $[\aa,\bb]$, come definito in
(\ref{albe}). Pertanto, usando la maggiorazione
$$ \max_{\lam\in\sigma(G)} | p(\lam) |\le \max_{\aa\le\lam\le\bb}
|p(\lam)|,$$

\no possiamo scegliere il polinomio $p_j$ come soluzione del
seguente problema di minimassimo vincolato:
$$ p_j = \arg\min_{p\in\Pi_j,\,p(1)=1}\,\max_{\aa\le\lam\le\bb}
|p(\lam)|. $$

\no La soluzione di questo problema (vedi propriet\`a C8 in
Sezione~\ref{cebypol}), \`e data da
\begin{equation}\label{siceb}
p_j(x) = \frac{T_j(\mu(x))}{T_j(\mu(1))}, \qquad \mu(x) =
\frac{2x-\aa-\bb}{\bb-\aa}, \end{equation}

\no dove $T_j$ \`e il $j$-esimo polinomio di Chebyshev di prima
specie. Poich\`e per $x\in[\aa,\bb]$, $\mu(x)\in[-1,1]$, e
$\mu(1)>1$, si conclude facilmente che (vedi propriet\`a C3 e C7
in Sezione~\ref{cebypol})
\begin{eqnarray}\label{pjsol}
\max_{\aa\le\lam\le\bb} |p_j(\lam)| &=& \frac{1}{|T_j(\mu(1))|} =
\frac{1} {T_j(\mu(1))} \\ \nonumber &\le& 2\left(
\mu(1)+\sqrt{\mu(1)^2-1}\right)^{-j}\rightarrow 0, \mbox{\qquad
per~}j\rightarrow\infty.
\end{eqnarray}

\no Osserviamo che
$$\mu(1) = 1 +2\frac{1-\bb}{\bb-\aa},$$

\no con, idealmente, $\aa=\lam_1$, $\bb=\lam_n$. Tuttavia, tali
valori sono in generale non noti, per cui valori di $\aa$ e $\bb$
come in (\ref{albe}) sono in generale utilizzati. Tuttavia, se
$\aa$ e $\bb$ sono vicini ai rispettivi estremi (-1 e 1), allora
$\mu(1)\approx 1$ e, quindi, l'errore si smorza assai lentamente.
Viceversa, se $[\aa,\bb]$ non contiene l'inviluppo convesso di
$\sigma(G)$, allora la (\ref{pjsol}) non \`e pi\`u valida.

In conclusione, esiste il problema della scelta ottimale dei
parametri $\aa$ e $\bb$, del tutto analogo alla scelta ottimale
del parametro di rilassamento $\om$ in SOR.

\begin{eser} Dimostrare la (\ref{pjsol}).\end{eser}

\begin{oss} Il metodo risultante dalla (\ref{siceb}), denominato
{\em metodo semi-iterativo di Chebyshev}, \`e effettivamente
semi-iterativo, perch\`e il numero di iterazioni per abbattere
l'errore di un prefissato fattore pu\`o essere determinato a
priori, per esempio dalla (\ref{pjsol}).\end{oss}

\begin{oss} Nel caso in cui la matrice $G$ non sia simmetrica, ma
sia diagonalizzabile con autovalori reali,
$$G = V\Lambda V^{-1},$$

\no $\Lambda$ come in (\ref{VLV})-(\ref{albe}), il metodo
semi-iterativo di Chebyshev risulta essere ancora efficace.
Infatti, in questo caso (vedi (\ref{pG})), $$\|p(G)\| =
\|p(V\Lambda V^{-1})\| = \|V p(\Lambda) V^{-1}\| \le
\kappa_2(V)\cdot\|p(\Lambda)\|,$$

\no dove $\kappa(V)$ \`e il numero di condizione (in norma-2)
della matrice $V$.\end{oss}

Anche supponendo che il problema della scelta ottimale dei
parametri $\aa$ e $\bb$ sia stato risolto, nondimeno
l'implementazione del metodo semi-iterativo di Chebyshev richiede
una attenta analisi, al fine di renderla efficiente. Infatti, la
formulazione appena vista prevederebbe di dover conservare tutte
le approssimazioni (\ref{xk}) precedentemente generate dal metodo
iterativo di base (\ref{split1}). In realt\`a, \`e possibile
evitare questo, osservando che i polinomi di Chebyshev di prima
specie soddisfano all'equazione alle differenze
$$T_{k+1}(x) = 2xT_k(x)-T_{k-1}(x), \quad k\ge1, \qquad
T_0(x)\equiv 1, \quad T_1(x) = x.$$

\no Definendo, quindi, la successione $\{\cc_k\}$,  dove (vedi
(\ref{siceb})) $\cc_k = T_k(\mu(1))$,
$$\cc_{k+1} = 2\mu(1)\cc_k-\cc_{k-1}, \quad k\ge1, \qquad \cc_0=1,
\quad \cc_1=\mu(1),$$

\no \`e possibile quindi ottenere una relazione di ricorrenza per
i polinomi $\{p_k(x)\}$:
$$p_{k+1}(x) = \frac{T_{k+1}(\mu(x))}{T_{k+1}(\mu(1))} =
\frac{1}{\cc_{k+1}}\left( 2\mu(x)\cc_k p_k(x)
-\cc_{k-1}p_{k-1}(x)\right).$$

\no Considerando che
$$\bfe^{(k)} = \bfx^*-\bfx^{(k)} = p_k(G)\bfe_0,$$

\no si ottiene, ponendo
\begin{equation}\label{muG}
\mu(G) = (\bb-\aa)^{-1}\left(2G-(\aa+\bb)I\right),
\end{equation}
\begin{equation}\label{xk1}
\bfx^{(k+1)} = 2\frac{\cc_k}{\cc_{k+1}}\mu(G)\bfx^{(k)}
-\frac{\cc_{k-1}}{\cc_{k+1}}\bfx^{(k-1)} +
\left(I-2\frac{\cc_k}{\cc_{k+1}}\mu(G)
+\frac{\cc_{k-1}}{\cc_{k+1}}I\right)\bfx^*,
\end{equation}

\no in cui $\bfx^*$ non \`e, ovviamente, nota. Tuttavia, dalla
relazione
$$(I-G)\bfx^* = (I-M^{-1}N)\bfx^* = M^{-1}(M-N)\bfx^* =
M^{-1}A\bfx^* = M^{-1}\bfb,$$

\no si ottiene:
\begin{eqnarray*}
\lefteqn{\left(I-2\frac{\cc_k}{\cc_{k+1}}\mu(G)
+\frac{\cc_{k-1}}{\cc_{k+1}}I\right)\bfx^* =}\\
&=&\frac{1}{\cc_{k+1}}\left( (\cc_{k+1}-2\mu(1)\cc_k+\cc_{k-1})I +
2\cc_k( \mu(1)I-\mu(G)) \right)\bfx^*\\ &=&
2\frac{\cc_k}{\cc_{k+1}}(\mu(1)I-\mu(G))\bfx^* =
\frac{\cc_k}{\cc_{k+1}}\frac{4}{\bb-\aa}(I-G)\bfx^* \\ &=&
\frac{\cc_k}{\cc_{k+1}}\frac{4}{\bb-\aa}M^{-1}\bfb.
\end{eqnarray*}

\no Da questa equazione e dalla (\ref{xk1}) si ottiene, infine,
tenendo conto della (\ref{muG}),
$$\bfx^{(k+1)} = \frac{\cc_k}{\cc_{k+1}}\frac{4}{\bb-\aa}
M^{-1}(N\bfx^{(k)}+\bfb) -
\frac{\cc_k}{\cc_{k+1}}\frac{2(\aa+\bb)}{\bb-\aa}\bfx^{(k)}
-\frac{\cc_{k-1}}{\cc_{k+1}}\bfx^{(k-1)}.$$

Tutti questi argomenti possono essere riassunti nel seguente
algoritmo, che descrive la implementazione efficiente del metodo
semi-iterativo di Chebyshev.

\begin{algo} Metodo semi-iterativo di Chebyshev.\\
\rm \fbox{\parbox{12cm}{
\begin{itemize}

\nulit siano assegnati $\aa,\bb$ ed il vettore iniziale
$\bfx^{(0)}$

\nulit inizializzo $\cc_0=1, \cc_1=\mu(1)$

\nulit  risolvo $M\bfx^{(1)} = N\bfx^{(0)}+\bfb$

\nulit per $k=1,2,\dots:$

\begin{itemize}

\nulit  $\cc_{k+1} = 2\mu(1)\cc_k-\cc_{k-1}$

\nulit risolvo $M\bfz_k = N\bfx^{(k)}+\bfb$

\nulit $\bfx^{(k+1)} = \frac{\cc_k}{\cc_{k+1}}\frac{4}{\bb-\aa}
\bfz_k -
\frac{\cc_k}{\cc_{k+1}}\frac{2(\aa+\bb)}{\bb-\aa}\bfx^{(k)}
-\frac{\cc_{k-1}}{\cc_{k+1}}\bfx^{(k-1)}$

\end{itemize}

\nulit fine per
\end{itemize}
}}\end{algo}

\section{Il metodo SOR simmetrico (SSOR)}\label{ssor}

Abbiamo visto che il metodo semi-iterativo di Chebyshev pu\`o
essere applicato quando la matrice di iterazione del metodo
sottostante (\ref{split1}), $$G=M^{-1}N$$ \`e simmetrica o,
almeno, simile ad una matrice simmetrica. Sarebbe auspicabile che
questo avvenisse, quando la matrice $A$ \`e sdp.

Questo sicuramente avviene per il metodo di Jacobi. Infatti, in
questo caso, se $A=D-L-L^\top $, con $D$ nonsingolare, allora
$$G = D^{-1}(L+L^\top ) \sim D^{-\frac{1}2}(L+L^\top )D^{-\frac{1}2},$$

\no che \`e simmetrica.

Tuttavia, il metodo di Jacobi \`e, tra i metodi iterativi di base,
quello con le caratteristiche di convergenza meno favorevoli.
Sarebbe, pertanto, auspicabile che la stessa propriet\`a valesse
per i metodi di Gauss-Seidel e SOR. Questo non \`e per\`o sempre
garantito, per cui si procede, ora, a definire il {\em metodo SOR
simmetrico (SSOR)}, che include, come caso particolare relativo a
$\om=1$, la simmetrizzazione del metodo di Gauss-Seidel.

Si procede, innanzitutto, a scalare simmetricamente la diagonale
di $A$. Se
$$A\bfx=\bfb, \qquad A = D-L-L^\top ,$$

\no allora questo equivale a risolvere
$$\hat{A}\hat\bfx = \hat\bfb,$$

\no dove $~\hat{A} = D^{-\frac{1}2}AD^{-\frac{1}2},~ \hat\bfb =
D^{-\frac{1}2}\bfb, ~ \hat\bfx = D^{\frac{1}2}\bfx$. Pertanto:
$$\hat{A} = I -\hat{L}-\hat{L}^\top , \qquad \hat{L} =
D^{-\frac{1}2}LD^{-\frac{1}2}.$$

\no Osserviamo che, se
\begin{equation}\label{MwNw}
\hat{M}_\om = \om^{-1}(I-\om \hat{L}),\qquad \hat{N}_\om
=\om^{-1}((1-\om)I+\om \hat{L}^\top ),
\end{equation}

\no allora
$$\hat{A}=\hat{M}_\om-\hat{N}_\om=\hat{M}_\om^\top -\hat{N}_\om^\top $$
\no sono entrambi splitting di $\hat{A}$. Conseguentemente, si
pu\`o suddividere il passo di SOR in 2 passi consecutivi
simmetrici: $$\hat{M}_\om \hat\bfx_{k+\frac{1}2} = \hat{N}_\om
\hat\bfx_k + \hat\bfb, \qquad \hat{M}_\om^\top  \hat\bfx_{k+1} =
\hat{N}_\om^\top  \hat\bfx_{k+\frac{1}2} +\hat\bfb,$$

\no che, cumulativamente, danno l'iterazione
\begin{equation}
\label{ssorit}\hat\bfx_{k+1} = \left(\hat{M}_\om^{-T}
\hat{N}_\om^\top  \hat{M}_\om^{-1} \hat{N}_\om\right) \hat\bfx_k +
\hat{M}_\om^{-T}\left(I+\hat{N}_\om^\top 
\hat{M}_\om^{-1}\right)\hat\bfb.\end{equation}

\no Considerando che (vedi (\ref{MwNw})) $\hat{M}_\om$ e
$\hat{N}_\om^\top $ commutano, si ha che la matrice di iterazione,
\begin{equation}\label{hGw}
\hat{G}_\om = \hat{M}_\om^{-T} \hat{N}_\om^\top  \hat{M}_\om^{-1}
\hat{N}_\om = \hat{M}_\om^{-T} \hat{M}_\om^{-1} \hat{N}_\om^\top 
\hat{N}_\om,\end{equation}

\no risulta essere simile ad una matrice simmetrica. Infatti,
$B_\om \equiv \hat{N}_\om^\top  \hat{N}_\om$ \`e simmetrica e
semi-definita positiva e, quindi, $$\hat{G}_\om \sim
\hat{M}_\om^{-1} B_\om\hat{M}_\om^{-T},$$

\no che \`e simmetrica e semi-definita positiva. Pertanto,
$\hat{G}_\om$ \`e diagonalizzabile, con autovalori reali e non
negativi, ed il metodo semi-iterativo di Chebyshev \`e
applicabile. \`E per questo motivo che il metodo, definito
formalmente dall'iterazione (\ref{ssorit}), \`e denominato SOR
simmetrico (SSOR).

\begin{figure}[t]
\begin{center}
\includegraphics[width=13cm,height=10cm]{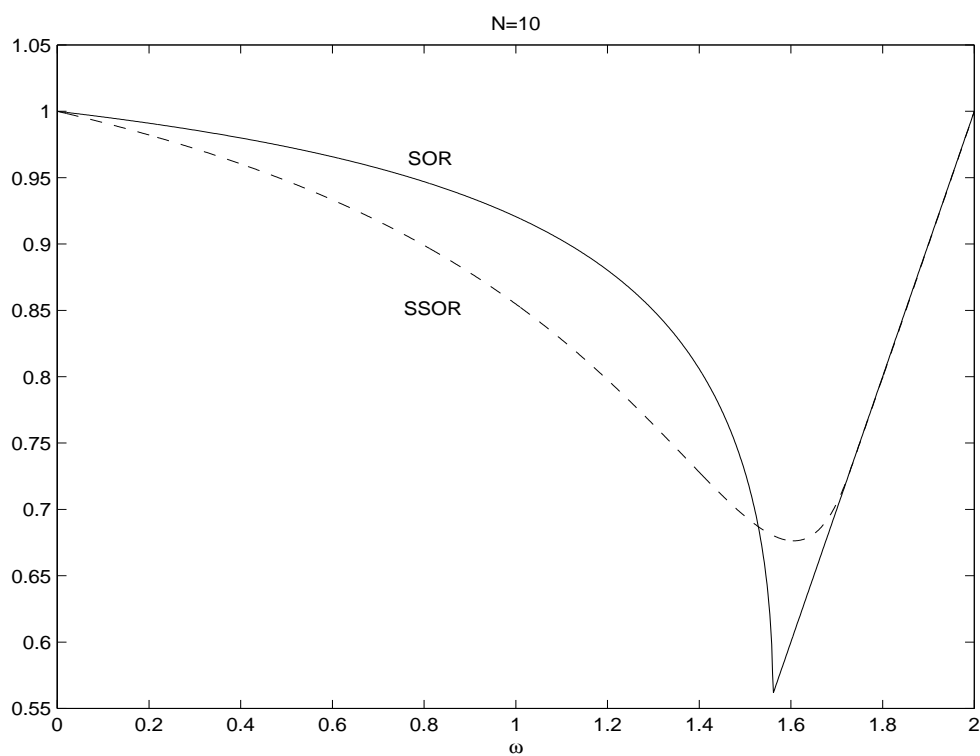}
\caption{\protect\label{fig3} Raggio spettrale di SOR e SSOR in
funzione del parametro di rilassamento $\om$, relativo al problema
(\ref{test}) dell'Esercizio~\ref{ex1}, $N=10$.}
\end{center}
\end{figure}

\begin{oss} Le matrici di iterazione del metodo
SOR, $G_\om$ (vedi (\ref{Gw})), e quella del metodo SSOR,
$\hat{G}_\om$ (vedi (\ref{hGw})), hanno, in ge\-nerale, differenti
raggi spettrali. A riprova di questo, in Figura~\ref{fig3} si
riportano i corrispondenti raggi spettrali, per $\om\in[0,2]$,
relativamente al problema (\ref{test}) dell'Esercizio~\ref{ex1},
$N=10$. Per confrontare adeguatamente le due curve bisogna,
tuttavia, considerare che una iterazione di SSOR ha un costo
computazionale doppio rispetto a quella di SOR.
\end{oss}

\begin{eser} Scrivere un codice che implementi il metodo
semi-iterativo di Chebyshev, applicato a SSOR, per il problema
(\ref{test}) in Esercizio~\ref{ex1}.\end{eser}

\begin{eser} Si supponga che, per un dato metodo iterativo, la
matrice di iterazione sia simmetrica e con autovalori
nell'intervallo $[-0.5,0.5]$. Quante iterazioni del metodo
semi-iterativo di Chebyshev sono richieste, al fine di ridurre
l'errore iniziale di un fattore (almeno) $10^6$? Quante iterazioni
richiede\-rebbe il metodo iterativo di base? (Considerare che
$\log_2(2+\sqrt{3})\approx1.9$).\end{eser}
 
%
%
\chapter{Il caso $A$ sdp}\label{cap4}

{\em
In questa sezione tratteremo il metodo del {\em gradiente
coniugato (GC)}, che \`e convenientemente impiegato nel caso in
cui la matrice $A$ del problema,
\begin{equation}\label{sistsdp}
A\bfx = \bfb,
\end{equation}

\no sia $n\times n$ e sdp.}

\section{Fattorizzazione $AU=UT$}\label{AUeqUT}

Cominciamo con il ricercare la seguente fattorizzazione della
matrice $A$:
\begin{eqnarray}\nonumber
AU&=&UT \equiv (\bfu_1 \dots \bfu_n)\pmatrix{cccc} \cc_1 &\bb_1\\
1 &\ddots &\ddots\\
      &\ddots &\ddots &\bb_{n-1}\\
      &       &1      &\cc_n\endpmatrix,
\\ \label{fattA}\\ \nonumber
      U^\top BU &=& D\equiv\pmatrix{ccc} d_1\\&\ddots\\&&d_n\endpmatrix,
\end{eqnarray}

\no in cui si assume che:

\begin{description}

\item[CG1:~] $B$ sia sdp;

\item[CG2:~] $AB=BA$.

\end{description}

Osserviamo che la seconda condizione in (\ref{fattA}) implica che
i vettori $\{\bfu_i\}$ sono $B$-ortogonali. Dal Teorema~\ref{Aort}
segue che, se essi sono non nulli, allora sono linearmente
indipendenti.

Vediamo in che modo sia possibile determinare la
fattorizzazione (\ref{fattA}), ovvero i vettori $\{\bfu_i\}$ e gli
scalari $\{\cc_i\}$, $\{\bb_i\}$ e $\{d_i\}$. Dalla seconda
equazione, segue immediatamente che
\begin{equation}\label {d_i}
d_i = \bfu_i^\top B\bfu_i \equiv \|\bfu_i\|_B^2, \qquad i=1,\dots,n.
\end{equation}

\no Inoltre, moltiplicando la prima equazione, membro a membro, per
$U^\top B$, si ottiene:
$$U^\top BAU = U^\top BUT = DT,$$

\no da cui segue che
$$d_i\cc_i = E_i^\top  DTE_i = E_i^\top U^\top BAUE_i = (UE_i)^\top  BA (UE_i) =
\bfu_i^\top BA\bfu_i,$$

\no e, quindi,
\begin{equation}\label{c_i}
\cc_i = \frac{\bfu_i^\top BA\bfu_i}{d_i}, \qquad i=1,\dots,n.
\end{equation}

\no Similmente, si ottiene
\begin{eqnarray*} d_{i-1}\bb_{i-1} &=& E_{i-1}^\top DTE_i = E_{i-1}^\top U^\top BAUE_i\\[1mm]
&=&(UE_{i-1})^\top BA(UE_i) = \bfu_{i-1}^\top BA\bfu_i,
\end{eqnarray*}

\no che dalla CG2, e considerando che
\begin{equation}\label{Aui}
A\bfu_i = AUE_i = UTE_i = \bb_{i-1}\bfu_{i-1} +\cc_i\bfu_i
+\bfu_{i+1},\end{equation}

\no permette di ricavare
\begin{eqnarray}\nonumber
\bb_{i-1} &=& \frac{\bfu_{i-1}^\top BA\bfu_i}{d_{i-1}} =
\frac{\bfu_i^\top BA\bfu_{i-1}}{d_{i-1}}\\[1mm] &=&
\frac{\bfu_i^\top B(\bb_{i-2}\bfu_{i-2}+\cc_{i-1}\bfu_{i-1}+\bfu_i)}{d_{i-1}}
= \frac{\bfu_i^\top B\bfu_i}{d_{i-1}} = \frac{d_i}{d_{i-1}}.
\label{b_i}
\end{eqnarray}

\no Per convenzione, nella (\ref{Aui}) poniamo $\bb_0=0$ e
$\bfu_{n+1}=\bfo$. Le relazioni (\ref{d_i})--(\ref{b_i})
permettono di definire l'Algoritmo~\ref{algo2}, che \`e semi-iterativo e
risulta essere definito fintantoch\'e il vettore $\bfu_i$ corrente
\`e non nullo.

\begin{algo}\label{algo2} Fattorizzazione (\ref{fattA}).\\ \rm
\fbox{\parbox{12cm}{
\begin{itemize}

\nulit sia assegnato un vettore iniziale $\bfu_1$

\nulit per i=1,\dots,n:

\begin{itemize}

\nulit $\bfv_i=A\bfu_i$

\nulit $\bfz_i =B\bfu_i$

\nulit $d_i = \bfu_i^\top \bfz_i$

\nulit $\cc_i = \bfv_i^\top \bfz_i/d_i$

\nulit se $i=1, \bb_0=0$, altrimenti, $\bb_{i-1}=d_i/d_{i-1},$
fine se

\nulit $\bfu_{i+1} = \bfv_i-\cc_i\bfu_i-\bb_{i-1}\bfu_{i-1}$

\end{itemize}

\nulit fine per

\end{itemize}
}}\end{algo}

Riguardo alla possibilit\`a che uno dei vettori $\bfu_i$ si
annulli, osserviamo che, dalle (\ref{Aui})-(\ref{b_i}) segue che,
se $\bfu_i=\bfo$, allora $\bfu_{i+1}=\bfo$. Conseguentemente,
sempre dalla (\ref{Aui}), tutti i successivi vettori sono nulli.
Questi argomenti permettono quindi di concludere quanto segue.

\medskip
\begin{teo}\label{cor4.1} $\bfu_i\ne\bfo\Rightarrow \bfu_j\ne\bfo, \, j\le
i$.\end{teo}
\medskip

Verifichiamo che i vettori $\{\bfu_i\}$ sono $B$-ortogonali.
Osserviamo che questa propriet\`a ne implica la lineare
indipendenza, se sono non nulli e, pertanto, la convenzione
$\bfu_{n+1}=\bfo$ rispecchia un dato di fatto.

\medskip
\begin{teo}\label{teo4.1} Se $\bfu_i\ne\bfo,\,(i\le n)$, allora i vettori
$\{\bfu_1,\dots,\bfu_i\}$ sono $B$-ortogonali.\end{teo}

\proof \`E sufficiente far vedere che, se $\bfu_{i+1}\ne\bfo$, allora
$$\bfu_j^\top B\bfu_{i+1}=0,\qquad j\le i,$$
cosa che dimostriamo per induzione su $i$. Per $i=1$, la tesi segue
dalle corrispondenti (\ref{c_i})-(\ref{Aui}), tenendo conto che $\bb_0=0$.
Supponiamo ora vera la tesi per $i-1$ e dimostriamo per $i$. Per $j=i$ e $j=i-1$,
la tesi segue, rispettivamente, da (\ref{c_i}) e (\ref{b_i}). Per
$j\le i-2$, tenendo conto di (\ref{Aui}), di CG2, e dell'ipotesi di
induzione, si ottiene infine:
\begin{eqnarray*}
\bfu_j^\top B\bfu_{i+1} &=&
\bfu_j^\top B(A\bfu_i-\cc_i\bfu_i-\bb_{i-1}\bfu_{i-1}) = \bfu_j^\top BA\bfu_i\\[1mm]
&=& \bfu_j^\top AB\bfu_i =
(\bfu_{j+1}+\cc_j\bfu_j+\bb_{j-1}\bfu_{j-1})^\top B\bfu_i = 0.\QED
\end{eqnarray*}

Nella prossima sezione vedremo come sia possibile utilizzare la
fattorizzazione (\ref{fattA}) per risolvere il sistema lineare
(\ref{sistsdp}).

\section{Il metodo dei gradienti coniugati (GC)}\label{CG}

Sia assegnato un vettore, $\bfx_0$, quale approssimazione iniziale
della soluzione $\bfx^*$ di (\ref{sistsdp}), e sia, in generale,
$$S_i = \bfx_0 + [\bfu_1,\dots,\bfu_i], \qquad i=1,2,\dots,$$

\no la variet\`a affine generata dai vettori
$\{\bfu_1,\dots,\bfu_i\}$ prodotti dall'Algoritmo~\ref{algo2}. Al
passo $i$-esimo dello stesso, ricerchiamo una nuova
approssimazione, $\bfx_i$, che sia la approssimazione ai minimi
quadrati, nella $B$-norma, di $\bfx^*$ nella variet\`a $S_i$. Per
quanto visto in Sezione~\ref{minquad}, tale soluzione \`e
caratterizza dalla $B$-ortogonalit\`a dell'errore,
$$\bfe_i =\bfx^*-\bfx_i,$$

\no ai vettori $\bfu_j,j\le i$. Ovvero,
\begin{equation}\label{eBort}
\bfu_j^\top B\bfe_i = 0, \qquad j\le i.
\end{equation}

Osserviamo che questa caratterizzazione della soluzione non \`e
immediatamente utilizzabile, in quanto la soluzione $\bfx^*$ non
\`e nota. Tuttavia, scegliendo
\begin{equation}\label{AeqB}
B=A,
\end{equation}

\no si ottiene che i requisiti CG1 e CG2 sono evidentemente
soddisfatti e, inoltre,
$$B\bfe_i = A\bfx^* -A\bfx_i = \bfb-A\bfx_i \equiv \bfr_i$$

\no \`e il {\em residuo} corrispondente ad $\bfx_i$.

\begin{oss} Poich\'e $A$ \`e nonsingolare,
$\bfr_i=\bfo \Leftrightarrow \bfx_i=\bfx^*$.\end{oss}

Quindi, la condizione di $B$-ortogonalit\`a dell'errore (\ref{eBort}),
a causa della scelta (\ref{AeqB}), risulta essere equivalente alla
condizione di ortogonalit\`a sul residuo,
\begin{equation}\label{rort}
\bfu_j^\top \bfr_i = 0, \qquad j\le i.
\end{equation}

Al fine di ottenere un procedimento iterativo, ricerchiamo la
approssimazione al passo $i$-esimo nella forma
\begin{equation}\label{x_i}
\bfx_i = \bfx_{i-1}+\lam_i\bfu_i.
\end{equation}

\no Ne consegue che
\begin{equation}\label{r_i}
\bfr_i = \bfr_{i-1}-\lam_iA\bfu_i.
\end{equation}

\no Imponendo la (\ref{rort}) per $j=i$, si ottiene
\begin{equation}\label{lami}
\lam_i = \frac{\bfu_i^\top \bfr_{i-1}}{\bfu_i^\top A\bfu_i} \equiv
\frac{\bfu_i^\top \bfr_{i-1}}{d_i}.
\end{equation}

Andiamo a verificare che tale scelta \`e sufficiente a garantire tutte
le condizioni di ortogonalit\`a  (\ref{rort}). Vale, infatti, il
seguente risultato.

\medskip
\begin{teo}\label{teo4.2}
Se i vettori $\{\bfu_1,\dots,\bfu_i\}$ sono non nulli e
$A$-ortogonali, e le approssimazioni sono definite da
(\ref{x_i})-(\ref{lami}), allora (\ref{rort}) \`e
verificata.\end{teo}

\proof Per induzione su $i$. Per $i=1$, la tesi \`e ovvia.
Supponiamo, quindi vero che
$$\bfu_j^\top \bfr_{i-1} = 0, \qquad j\le i-1,$$

\no e dimostriamo la (\ref{rort}). Per $j=i$ la tesi segue dalla
(\ref{lami}). Per $j\le i-1$, dalla (\ref{r_i}) segue
$$\bfu_j^\top \bfr_i = \bfu_j^\top \bfr_{i-1}-\lam_i\bfu_j^\top A\bfu_i = 0,$$

\no per l'ipotesi di induzione e la $A$-ortogonalit\`a dei vettori
$\{\bfu_j\}$.\QED

\medskip
\begin{cor}\label{cor4.2}
Se i vettori $\{\bfu_1,\dots,\bfu_n\}$ sono non nulli e
$A$-ortogonali, allora $\bfx_i=\bfx^*$, per qualche $i\le
n$.\end{cor}

\proof \`E sufficiente osservare che il residuo, al passo $n$, non
pu\`o che essere il vettore nullo, dalla (\ref{rort}), essendo i
vettori $\{\bfu_i\}$ una base per $\RR^n$.\QED\bigskip

Rimane da vedere cosa accade se qualcuno dei vettori $\bfu_i$ si
annulla, per $i\le n$. A questo fine, effettueremo la scelta del
vettore $\bfu_1$, che definisce l'intera successione, mediante
l'Algoritmo~\ref{algo2} con $B=A$, come
\begin{equation}\label{u1}
\bfu_1 = \bfr_0 \equiv \bfb-A\bfx_0.
\end{equation}

\no Questa scelta definisce il metodo dei {\em gradienti coniugati}.
Si premette il seguente risultato.

\medskip
\begin{lem}\label{Krylov} Se $\bfx_i\ne\bfx^*$, allora:
$$[\bfu_1,\dots,\bfu_{i+1}] =
[\bfu_1,A\bfu_1,\dots,A^i\bfu_1] = [\bfr_0,\dots,\bfr_i].
$$\end{lem}

\proof Dimostriamo la tesi per induzione su $i$. Per $i=0$, essa
discende dalla scelta (\ref{u1}) del vettore $\bfu_1$.
Supponiamola vera per $i-1$ e dimostriamo per $i$.

\begin{description}
\item[\mbox{$[\bfu_1,\dots,\bfu_{i+1}]=[\bfu_1,\dots,A^i\bfu_1]$.~}] Infatti, dalla
(\ref{Aui}) e dall'ipotesi di induzione, segue che
$$\bfu_{i+1}=A\bfu_i-\cc_i\bfu_i-\bb_{i-1}\bfu_{i-1},$$
con
\begin{eqnarray*}
\bfu_{i-1},\bfu_i&\in&[\bfu_1,\dots,A^{i-1}\bfu_1]\subseteq[\bfu_1,\dots,A^i\bfu_1],\\[1mm]
A\bfu_i&\in&[A\bfu_1,\dots,A^i\bfu_1]\subseteq[\bfu_1,A\bfu_1,\dots,A^i\bfu_1].
\end{eqnarray*}
Pertanto $\bfu_{i+1}\in[\bfu_1,\dots,A^i\bfu_1]$. Viceversa,
sempre dalla (\ref{Aui}), segue che $$A^i\bfu_1 =
A(A^{i-1}\bfu_1)\in[A\bfu_1,\dots,A\bfu_i]\subseteq[\bfu_1,\dots,\bfu_{i+1}].$$

\item[\mbox{$[\bfu_1,\dots,\bfu_{i+1}]=[\bfr_0,\dots,\bfr_i]$.~}]
Infatti, dall'ipotesi di induzione e dalla (\ref{r_i}), segue, in virt\`u della (\ref{Aui}), che
$$\bfr_i = \bfr_{i-1}-\lam_i
A\bfu_i\in[\bfu_1,\dots,\bfu_i,A\bfu_i]\subseteq[\bfu_1,\dots,\bfu_{i+1}].$$
Viceversa, sempre dalla (\ref{r_i}), segue che, se
$\bfx_i\ne\bfx^*$, allora il parametro $\lam_i$ dato dalla
(\ref{lami}), deve essere non nullo. Infatti, se cos\`\i\, non
fosse, dalla (\ref{rort}) seguirebbe che $$\bfr_i=\bfr_{i-1}\in[\bfu_1,\dots,\bfu_i], \qquad
\bfu_j^\top \bfr_i=0,\quad j\le i,$$ per cui dovrebbe aversi
$\bfr_i=\bfo$, contro l'ipotesi fatta. Pertanto, $$A\bfu_i =
\frac{\bfr_{i-1}-\bfr_i}{\lam_i} =
\bfu_{i+1}+\cc_i\bfu_i+\bb_{i-1}\bfu_{i-1},$$ da cui, per
l'ipotesi di induzione, segue infine che $$\bfu_{i+1} =
\frac{\bfr_{i-1}-\bfr_i}{\lam_i}
-\cc_i\bfu_i-\bb_{i-1}\bfu_{i-1}\in[\bfr_0,\dots,\bfr_i].~\QED$$
\end{description}

\medskip
\begin{cor}\label{cor4.3} I residui prodotti dal metodo dei gradienti coniugati
sono tra di loro ortogonali.
\end{cor}

\proof Segue immediatamente dal precedente Lemma~\ref{Krylov} e dalla
(\ref{rort}).\QED

\begin{algo}\label{algo3} Metodo dei gradienti coniugati.\\ \rm
\fbox{\parbox{12cm}{
\begin{itemize}

\nulit sia assegnata una approssimazione iniziale $\bfx_0$

\nulit inizializzo $\bfu_1=\bfr_0 \equiv \bfb-A\bfx_0$

\nulit per i=1,\dots,n:

\begin{itemize}

\nulit se $\bfr_{i-1}=\bfo, \bfx^*=\bfx_{i-1}$, esci, fine se

\nulit $\bfv_i=A\bfu_i$

\nulit $d_i = \bfu_i^\top \bfv_i$

\nulit $\lam_i = \bfu_i^\top \bfr_{i-1}/d_i$

\nulit $\bfx_i = \bfx_{i-1}+\lam_i\bfu_i$

\nulit $\bfr_i = \bfr_{i-1}-\lam_i\bfv_i$

\nulit $\cc_i = \bfv_i^\top \bfv_i/d_i$

\nulit se $i=1, \bb_0=0$, altrimenti, $\bb_{i-1}=d_i/d_{i-1},$
fine se

\nulit $\bfu_{i+1} = \bfv_i-\cc_i\bfu_i-\bb_{i-1}\bfu_{i-1}$

\end{itemize}

\nulit fine per

\nulit se $i=n, \bfx^*=\bfx_n$, fine se
\end{itemize}

}}\end{algo}

\medskip
\begin{cor}\label{cor4.4} $\bfu_i\ne\bfo \mbox{~e~} \bfu_{i+1}=\bfo \Rightarrow
\bfx_i = \bfx^*$.\end{cor}

\proof Infatti, se cos\`\i\, non fosse, si avrebbe
$$[\bfu_1,\dots,\bfu_{i+1}]\equiv[\bfu_1,\dots,\bfu_i]\subset[\bfr_0,\dots,\bfr_i],$$
contraddicendo, in virt\`u del Corollario~\ref{cor4.3}, il risultato del Lemma~\ref{Krylov}.\QED

\medskip
\begin{cor}\label{cor4.5} In aritmetica esatta, il metodo dei GC converge in al pi\`u
$n$ passi.\end{cor}
\vspace{2mm}

I precedenti risultati permettono di ottenere
l'Algo\-ritmo~\ref{algo3} per il metodo dei gradienti coniugati.

\section{Implementazione efficiente del metodo dei GC}\label{CGeff}

Esaminiamo brevemente il costo computazionale del precedente
Algoritmo~\ref{algo3} per il metodo dei gradienti coniugati. Tale
costo pu\`o essere quantificato in termini di occupazione di
memoria e numero di operazioni. Di entrambe le cose,
considereremo solo le quantit\`a principali.

\begin{description}

\item[Occupazione di memoria:~] avremo sicuramente bisogno dello
spazio ne\-cessario per allocare la matrice $A$, sebbene questo
richieda, nel caso sparso, generalmente $O(n)$ celle di memoria
(vedi Appendice~\ref{appA}).
Delle approssimazioni della soluzione, solo quella corrente \`e
necessaria per produrre la successiva, pertanto un unico vettore
$\bfx$ \`e sufficiente. Lo stesso discorso vale per i residui, dei
quali \`e necessario solo quello corrente. Due vettori sono invece
necessari per le direzioni generate, pi\`u un vettore di appoggio
che \`e inizialmente utilizzato per accumulare il prodotto
matrice-vettore ({\em matvec}) e, successivamente, per accumulare
la nuova direzione.

\item[Operazioni:~] suddivideremo queste ultime in base ad una
semplice classificazione. Abbiamo gi\`a introdotto il {\em matvec}, di
cui ne necessita 1 per iterazione. Similmente, si vede che sono
necessari 3 prodotti scalari ({\em scal}\,) e 4 {\em axpy} (che \`e
l'acronimo per ``alpha x plus y''), ovvero operazioni del tipo

\centerline{vettore = scalare$\times$vettore+vettore.}
\end{description}

\no Riassumendo, il costo computazionale
dell'Algoritmo~\ref{algo3} \`e dato da 1 matrice e 5 vettori,
riguardo alla occupazione di memoria, e da 1 matvec, 3 scal e 4
axpy per iterazione.

\begin{eser} Riscrivere l'Algoritmo~\ref{algo3} tenendo conto
delle sole strutture di memoria effettivamente richieste nella sua
implementazione.\end{eser}

Vediamo adesso come riformulare il metodo dei gradienti
coniugati, in modo da renderne pi\`u efficiente l'implementazione.
A tal fine, vale il seguente risultato.

\begin{teo}\label{teo4.3} Se $\bfx_i\ne\bfx^*$, allora
$\bfu_{i+1}\in[\bfu_i,\bfr_i]$.\end{teo}

\proof Distinguiamo i due casi $i=1$ e $i>1$. Per $i=1$, dalla
(\ref{u1}) e dalla (\ref{lami}) segue che $\bfr_1\ne\bfo$ e
$\lam_1\ne0$. Pertanto, dalla (\ref{Aui}) e dalla (\ref{r_i}) si
ottiene
$$\bfu_2 = A\bfu_1 -\cc_1\bfu_1 =
\frac{\bfr_0-\bfr_1}{\lam_1}-\cc_1\bfu_1,$$

\no da cui, ricordando che $\bfu_1=\bfr_0$, si conclude che
$\bfu_2\in[\bfu_1,\bfr_1]$. Per quanto riguarda il caso $i>1$,
nella dimostrazione del Lemma~\ref{Krylov} abbiamo gi\`a visto che, se
$\bfx_i\ne\bfx^*$, allora $\lam_i$, come definito da (\ref{lami}),
deve essere non nullo. Ne consegue che
\begin{eqnarray*}
\bfu_{i+1} &=& A\bfu_i -\cc_i\bfu_i-\bb_{i-1}\bfu_{i-1} =
\frac{\bfr_{i-1}-\bfr_i}{\lam_i}-\cc_i\bfu_i-\bb_{i-1}\bfu_{i-1}\\[1mm]
&=&-\frac{\bfr_i}{\lam_i}-\cc_i\bfu_i+\left(\frac{\bfr_{i-1}}{\lam_i}-\bb_{i-1}\bfu_{i-1}\right)
\equiv-\frac{\bfr_i}{\lam_i}-\cc_i\bfu_i+\bfw_i.
\end{eqnarray*}

\no Pertanto, considerando che dal Lemma~\ref{Krylov} segue che
$\bfw_i\in[\bfu_1,\dots,\bfu_i]$, la tesi risulta dimostrata se proviamo che
$\bfw_i\in[\bfu_i]$ o, equivalentemente, che
$\bfw_i^\top A\bfu_j=0,\,j<i$. Solo i seguenti due casi sono
possibili:

\begin{description}

\item[$j<i-1$:~]in tal caso,
$$\bfw_i^\top A\bfu_j =
\frac{\bfr_{i-1}^\top A\bfu_j}{\lam_i}-\bb_{i-1}\bfu_{i-1}^\top A\bfu_j =
\frac{\bfr_{i-1}^\top (\bfu_{j+1}+\cc_j\bfu_j+\bb_{j-1}\bfu_{j-1})}{\lam_i}=0,$$

\no in virt\`u della (\ref{rort}) e della $A$-ortogonalit\`a dei
vettori $\{\bfu_j\}$;

\item[$j=i-1$:~] dalle (\ref{b_i}), (\ref{rort}) e (\ref{lami}) segue che
\begin{eqnarray*}
\bfw_i^\top A\bfu_{i-1} &=&
\left(\frac{\bfr_{i-1}}{\lam_i}-\bb_{i-1}\bfu_{i-1}\right)^\top A\bfu_{i-1}\\[1mm]
&=&\frac{\bfr_{i-1}^\top (\bfu_i+\cc_{i-1}\bfu_{i-1}+\bb_{i-2}\bfu_{i-2})}{\lam_i}-\bb_{i-1}\bfu_{i-1}^\top A\bfu_{i-1}\\[1mm]
&=&\frac{\bfr_{i-1}^\top \bfu_i}{\lam_i}-d_i = d_i-d_i = 0.~\QED
\end{eqnarray*}
\end{description}

Il precedente risultato ci dice che possiamo ricercare la nuova
direzione, $\bfu_{i+1}$, come combinazione del residuo e della
direzione correnti, ovvero,
$$\bfu_{i+1} = c_1^{(i)} \bfr_i + c_2^{(i)}\bfu_i,$$

\no in cui i due coefficienti della combinazione devono essere
scelti in modo da soddisfare solo la condizione di ortogonalit\`a
\begin{equation}\label{lastone}
\bfu_{i+1}^\top A\bfu_i=0,\end{equation}

\no in quanto si verifica facilmente che la
condizione $\bfu_{i+1}^\top A\bfu_j=0$, $j<i$, risulta essere verificata
in ogni caso.  Osserviamo che le coppie che
soddisfano la (\ref{lastone}) sono definite a meno di una costante
moltiplicativa non nulla. Inoltre, le approssimazioni $\{\bfx_i\}$ (e quindi i
residui $\{\bfr_i\}$) sono invarianti, rispetto alla norma dei
vettori $\{\bfu_i\}$. Di conseguenza, definiamo dei nuovi vettori,
\begin{equation}\label{pi}
\bfp_i = \bfu_i\sigma_i, \quad i=1,2,\dots,\qquad \sigma_1=1,
\end{equation}

\no per i quali la costante $c_1^{(i)}$ sia normalizzata a 1.
Ponendo $c_2^{(i)}=\mu_i$, ricer\-chiamo, quindi
\begin{equation}\label{newpi}
\bfp_{i+1} = \bfr_i +\mu_i\bfp_i, \qquad i = 1,2,\dots,
\end{equation}

\no dove $\mu_i$ si ottiene imponendo la condizione di
ortogonalit\`a $\bfp_i^\top A\bfp_{i+1}=0$. Al fine di ricavare una
espressione conveniente per tale quantit\`a, osserviamo che, in
analogia con (\ref{d_i}) e (\ref{x_i})--(\ref{lami}), ora si
avranno, rispettivamente,
\begin{eqnarray*}
\hat{d}_i &=& \bfp_i^\top A\bfp_i,\\[1mm]
\bfx_i &=& \bfx_{i-1}+\hat\lam_i\bfp_i, \\[1mm]
\bfr_i &=& \bfr_{i-1}-\hat\lam_iA\bfp_i,\qquad i=1,2,\dots,
\end{eqnarray*}

\no dove
$$
\hat\lam_i = \frac{\bfp_i^\top \bfr_{i-1}}{\bfp_i^\top A\bfp_i} =
\frac{(\bfr_{i-1}+\mu_{i-1}\bfp_{i-1})^\top \bfr_{i-1}}{\hat{d}_i}
= \frac{\bfr_{i-1}^\top \bfr_{i-1}}{\hat{d}_i}.
$$

\no Nell'ultima equazione, si \`e sfruttato il fatto che condizioni
di ortogonalit\`a simili alla (\ref{rort}) valgono anche per i vettori
$\{\bfp_i\}$. Questo permette di ottenere, infine,
$$0 = \bfp_i^\top A\bfr_i +\mu_i \bfp_i^\top A\bfp_i =
\frac{(\bfr_{i-1}-\bfr_i)^\top }{\hat\lam_i}\bfr_i +\mu_i \hat{d}_i
  = -\frac{\bfr_i^\top \bfr_i}{\bfr_{i-1}^\top \bfr_{i-1}}\hat{d}_i
  +\mu_i\hat{d_i},$$

\no ovvero,
$$
\mu_i = \frac{\bfr_i^\top \bfr_i}{\bfr_{i-1}^\top \bfr_{i-1}}, \qquad
i=1,2,\dots.
$$

\no Questo ci permette di formulare l'Algoritmo~\ref{algo4}
per il metodo dei gradien\-ti coniugati, che risulta essere pi\`u
efficiente dell'Algoritmo~\ref{algo3} precedentemente proposto.
Si verifica facilmente, infatti, che il primo necessita di un vettore in
meno (per quanto concerne l'occupazione di memoria) e di una scal
ed una axpy in meno per iterazione, rispetto al secondo.

\begin{figure}[t]
\begin{algo}\label{algo4} Metodo dei gradienti coniugati.\\ \rm
\fbox{\parbox{12cm}{
\begin{itemize}

\nulit sia assegnata una approssimazione iniziale $\bfx_0$

\nulit inizializzo $\bfp_1=\bfr_0 \equiv \bfb-A\bfx_0$ e
 $\eta_0=\bfr_0^\top \bfr_0$

\nulit per i=1,\dots,n:

\begin{itemize}

\nulit se $\eta_{i-1}=0, \bfx^*=\bfx_{i-1}$, esci, fine se
\qquad\qquad\qquad{\bf (*)}

\nulit $\bfv_i=A\bfp_i$

\nulit $\hat{d}_i = \bfp_i^\top \bfv_i$

\nulit $\hat\lam_i = \eta_{i-1}/\hat{d}_i$

\nulit $\bfx_i = \bfx_{i-1}+\hat\lam_i\bfp_i$

\nulit $\bfr_i = \bfr_{i-1}-\hat\lam_i\bfv_i$

\nulit $\eta_i = \bfr_i^\top \bfr_i$

\nulit $\mu_i = \eta_i/\eta_{i-1}$

\nulit $\bfp_{i+1} = \bfr_i+\mu_i\bfp_i$

\end{itemize}

\nulit fine per

\nulit se $i=n, \bfx^*=\bfx_n$, fine se
\end{itemize}

}}\end{algo}
\end{figure}

\begin{eser} Riscrivere l'Algoritmo~\ref{algo4} tenendo conto
delle sole strutture di memoria effettivamente richieste nella sua
implementazione.\end{eser}

\begin{eser}\label{ex2} Supponendo che l'Algoritmo~\ref{algo4}
sia iterato fino ad $i=n$, verificare
che esso definisce la seguente fattorizzazione, alternativa alla
(\ref{fattA}) con $B=A$:
\begin{equation}\label{APPT}
AP = P\hat{T}, \qquad P^\top AP = \hat{D} \equiv \pmatrix{ccc}
\hat{d}_1\\ &\ddots \\ &&\hat{d}_n\endpmatrix,
\end{equation}

\no dove (vedi (\ref{pi}))
$$P \equiv (\bfp_1\dots\bfp_n) = U\Sigma \equiv U\pmatrix{ccc} \sigma_1 \\
&\ddots\\&&\sigma_n\endpmatrix, \qquad \hat{D}=D\Sigma^2,$$

\no e la matrice $\hat{T}=\Sigma^{-1}T\Sigma$ \`e ancora
tridiagonale. Ricavare, quindi, la relazione di ricorrenza, simile
alla (\ref{Aui}), cui soddisfano i vettori $\{\bfp_i\}$.
\end{eser}

\begin{eser}\label{ex3} Tenendo conto della (\ref{newpi}), e
definendo le matrici
$$R=(\bfr_0,\dots,\bfr_{n-1}), \qquad B = \pmatrix{cccc}
1 &-\mu_1 \\ &\ddots &\ddots\\ &&\ddots
&-\mu_{n-1}\\&&&1\endpmatrix,$$

\no dimostrare che:

\begin{enumerate}

\item $R = PB;$

\item $\Delta^2 = R^\top R$ \`e diagonale;

\item $Z=R\Delta^{-1}$ \`e ortogonale.

\end{enumerate}

\no Verificare, inoltre, che
\begin{equation}\label{barT}
A\sim Z^\top AZ \equiv \bar{T} = \pmatrix{cccc}
a_1 &-b_1\\ -b_1 &\ddots &\ddots\\ &\ddots &\ddots
&-b_{n-1}\\&&-b_{n-1} & a_n\endpmatrix,
\end{equation}

\no dove
\begin{eqnarray*}
a_1 &=& \frac{\hat{d}_1}{\|\bfr_0\|^2}, \\
a_{i+1} &=& \frac{\hat{d}_{i+1}+\hat{d}_i\mu_i^2}{\|\bfr_i\|^2},\\
b_i &=& \frac{\hat{d}_i\mu_i}{\|\bfr_{i-1}\|\,\|\bfr_i\|}, \qquad
i=1,\dots,n-1.
\end{eqnarray*}

\no Concludere che l'Algoritmo~\ref{algo4} equivale alla
fattorizzazione $LU$ della matrice tridiagonale $\bar{T}$.
(Suggerimento: utilizzare la fattorizzazione (\ref{APPT}) e la
relazione $R=PB$.)

\end{eser}

\begin{eser}\label{ex4.5} Con la notazione degli Esercizi~\ref{ex2} e \ref{ex3},
dimostrare che
\begin{equation}\label{Tcappello}
\hat{T} = B L \Lambda^{-1} \equiv \pmatrix{ccc} \frac{1+\mu_1}{\hat\lam_1}
&-\frac{\mu_1}{\hat\lam_2}\\ -\frac{1}{\hat\lam_1} &\ddots &\ddots \\
              &\ddots\endpmatrix,\qquad (\mu_n=0),\end{equation}

\no con $$L = \pmatrix{rrrr} 1\\-1&\ddots\\&\ddots&\ddots\\&&-1&1\endpmatrix,
\qquad \Lambda = \pmatrix{ccc} \hat\lam_1\\ &\ddots \\
&&\hat\lam_n\endpmatrix.$$\end{eser}

\medskip
\begin{oss}\label{cglanc} Osserviamo che la equazione (\ref{barT})
definisce, implicitamente, la fattorizzazione
\begin{equation}\label{AZZT}
AZ = Z\bar{T}, \qquad Z^\top Z = I,
\end{equation}

\no e che, dalla (\ref{rort}), dal Lemma~\ref{Krylov} e dal Corollario~\ref{cor4.3},
segue che, al passo $i$-esimo, il metodo dei gradienti coniugati ricerca la
soluzione $\bfx_i$ in modo tale che, se $$ Z = (\bfz_1\dots\bfz_n), \qquad
\bfz_i =\frac{\bfr_{i-1}}{\|\bfr_{i-1}\|},$$

\no allora il residuo $\bfr_i=\bfb-A\bfx_i$ soddisfa le condizioni
di ortogonalit\`a
\begin{equation}\label{rortz}
\bfz_j^\top \bfr_i = 0, \qquad j\le i.
\end{equation}
\end{oss}

\section{Criteri di arresto}\label{stopcg}

Abbiamo visto che il metodo del gradiente coniugato converge, in
aritmetica esatta, in al pi\`u $n$ iterazioni, se $n$ \`e la
dimensione del sistema lineare da risolvere. Va sottolineato che
spesso, a causa degli errori dovuti all'utilizzo dell'aritmetica
finita, questa propriet\`a pu\`o non essere effettiva.

D'altronde, talora \`e possibile raggiungere una soluzione
soddisfacente con un numero di iterazioni molto pi\`u esiguo.
Nasce, quindi, l'esigenza di riguardare il metodo pi\`u come un
metodo iterativo, che semi-iterativo, cui associare un opportuno
criterio di arresto. Nell'Algoritmo~\ref{algo4} questo \`e
rap\-presentato dal controllo sul residuo, che abbiamo etichettato
con (*). Tuttavia, il controllo in questa forma \`e stato inserito
pi\`u per un motivo concettuale, che pratico. In effetti, la
condizione $\bfr_i=\bfo$ potrebbe non essere mai soddisfatta in
aritmetica finita. Un controllo, molto pi\`u ragionevole, potrebbe
essere
$$\|\bfr_i\|^2\equiv\eta_i \le \mbox{\em tol},$$

\no essendo {\em tol} una opportuna tolleranza.
Una variante di questo controllo \`e costituita da
$$\eta_i \le \mbox{\em tol}\cdot\eta_0,$$

\no ovvero la norma (al quadrato) del residuo corrente \`e stata
abbattuto di un fattore {\em tol} rispetto a quella del residuo iniziale.
Sebbene questi due criteri siano molto semplici e, di fatto, assai
utilizzati in pratica, essi non sono in grado di assicurare {\em a
priori} una prefissata accuratezza della soluzione approssimata.
Infatti, se $\bfe_i = \bfx^*-\bfx_i$ \`e l'errore al passo
$i$-esimo, allora esso \`e soluzione del sistema lineare
$$A\bfe_i = \bfr_i,$$

\no e, pertanto,
$$\bfe_i = A^{-1}\bfr_i.$$

\no Da questo si ricava la diseguaglianza
\begin{equation}\label{stop1}
\|\bfe_i\|\le\|A^{-1}\|\cdot\|\bfr_i\|,
\end{equation}

\no la quale ci dice chiaramente che, se la matrice $A$ \`e molto
malcondizionata, allora $\|\bfe_i\|$ pu\`o essere grande, anche se
$\|\bfr_i\|$ \`e piccola.

\bigskip
Ad esempio, prendendo come matrice $A=(a_{ij})$ la matrice di Hilbert di
dimesione $n$, i cui elementi sono dati da $a_{ij} = (i+j-1)^{-1},
i,j=1,\dots,n$, si ha $\|A^{-1}\|\approx10^{13}$, per $n=10$. Ne consegue
che la norma del residuo pu\`o essere amplificata di un fattore $10^{13}$.
Infatti, se $\bfb$ \`e scelto in modo tale che la soluzione sia
$\bfx^*=(1\dots1)^\top $, e si considera l'errore (arrotondato alla
quinta cifra significativa) dato da
$$\bfe = \pmatrix{r}
  1.6740e-003\\
 -1.4557e-001\\
  3.1159e+000\\
 -2.8438e+001\\
  1.3607e+002\\
 -3.7503e+002\\
  6.1665e+002\\
 -5.9702e+002\\
  3.1392e+002\\
 -6.9130e+001\endpmatrix\qquad \Rightarrow \qquad\|\bfe\|=10^3,$$

\no si verifica che ad esso corrisponde un residuo
$$\bfr = \pmatrix{r}
 -5.3291e-015\\
 -1.3767e-014\\
  3.3862e-013\\
 -3.1082e-012\\
  1.4869e-011\\
 -4.1000e-011\\
  6.7417e-011\\
 -6.5260e-011\\
  3.4316e-011\\
 -7.5624e-012\endpmatrix\qquad \Rightarrow\qquad  \|\bfr\|\approx
 1.0932e-010,$$

\no da cui si evince che il fattore di amplificazione \`e
dell'ordine di $10^{13}$.

\begin{eser}\label{c_stop} Dimostrare che, se si utilizza il criterio di arresto
dato da $$\|\bfr_i\|\le\eps \cdot\|\bfb\|,$$ allora l'errore
soddisfa la diseguaglianza
$$\|\bfe_i\|\le\eps\cdot\kappa(A)\cdot\|\bfx^*\|,$$ dove
$\kappa(A)=\|A\|\cdot\|A^{-1}\|$ \`e il numero di condizionamento
di $A$ e $\bfx^*$ \`e, al solito, la soluzione del
problema.\end{eser}

\bigskip
Tuttavia, dalla (\ref{stop1}) si ottiene che, poich\'e
$$\|A^{-1}\| = \lam_{min}^{-1} \equiv \left(\min_{\lam\in\sigma(A)}
\lam\right)^{-1},$$

\no allora \`e possibile utilizzare il controllo sull'errore dato da
\begin{equation}\label{stop2}
\|\bfe_i\|\le\frac{\|\bfr_i\|}{\lam_{min}},
\end{equation}

\no a patto che sia noto l'autovalore minimo di $A$, o una sua opportuna
approssimazione. A questo riguardo, vale, per matrici simmetriche, la seguente
propriet\`a di {\em interlacciamento} degli autovalori.

\begin{teo} Sia $A\in\RR^{n\times n}$ una matrice simmetrica. Si
denoti con $A_r$ la sua sottomatrice principale di
ordine $r$, e con $$\lam_1(A_r)\le\cdots\le\lam_r(A_r)$$ i
corrispondenti autovalori. Allora, per $r=1,2,\dots,n-1$:
\begin{eqnarray*}
\lefteqn{\lam_1(A)\equiv\lam_1(A_n)\le\lam_1(A_{r+1})\le\lam_1(A_r)
\le\lam_2(A_{r+1})\le}\\[1mm]
&&\lam_2(A_r)\le
\cdots\le\lam_r(A_r)\le\lam_{r+1}(A_{r+1})
\le\lam_n(A_n)\equiv\lam_n(A).\end{eqnarray*}
\end{teo}

Pertanto, essendo la matrice $A$ simile alla matrice $\bar{T}$ in
(\ref{barT}), di cui al passo $i$-esimo possiamo costruire la sottomatrice
principale $\bar{T}_i$, possiamo approssimarne gli autovalori estremi con
quelli estremi di tale sottomatrice. Alternativamente, si pu\`o ricorrere alla
matrice $\hat{T}$ (vedi (\ref{Tcappello})), che si vede facilmente essere
simile ad una matrice tridiagonale simmetrica. Va altres\`\i\, sottolineato
che \`e generalmente richiesto un esiguo numero di ite\-razioni per pervenire
a soddisfacenti approssimazioni degli autovalori estremi (qualche decina, in
genere).

\begin{figure}[t]
\begin{center}
\includegraphics[width=13cm,height=10cm]{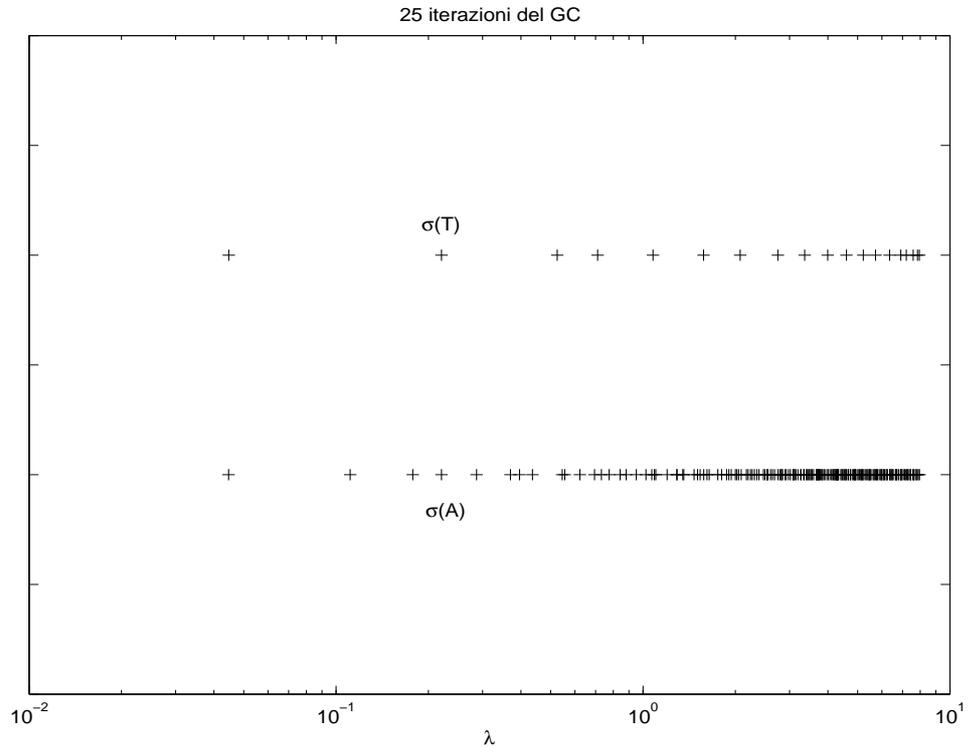}
\caption{\protect\label{fig4} Spettro della matrice $A$ dell'equazione
(\ref{test}), $N=20$, e della sua approssimazione dopo 25 iterazioni del
metodo dei GC.}
\end{center}
\end{figure}

Ad esempio, per la matrice $A$ dell'equazione (\ref{test}) in
Esercizio~\ref{ex1} si ha, per $N=20$ (quindi,
$A\in\RR^{400\times400}$), $$\lam_{min}(A) \approx 4.4677e-2,
\qquad \lam_{max}(A)\approx 7.9553.$$

\no Dopo 25 iterazioni del metodo del gradiente coniugato, si ottengono le
seguenti approssimazioni,
$$\lam_{min}(\bar{T}_{25})\approx 4.4682e-2, \qquad
\lam_{max}(\bar{T}_{25})\approx 7.8636,$$

\no che sono delle approssimazioni decisamente soddisfacenti, in questo
contesto. Per completezza, in Figura~\ref{fig4} si riporta l'intero spettro di
$A$ (indicato, nella figura, con $\sigma(A)$) e quello di $\bar{T}_{25}$
(indicato con $\sigma(T)$).

\section{Velocit\`a di convergenza del metodo dei GC}

In Sezione~\ref{CG} abbiamo visto che il metodo dei gradienti
coniugati minimizza, al passo $i$-esimo, la $A$-norma dell'errore,
$$\|\bfe_i\|_A = \|\bfx^*-\bfx_i\|_A = \sqrt{
(\bfx^*-\bfx_i)^\top A(\bfx^*-\bfx_i)},$$

\no nella variet\`a affine
$$S_i = \bfx_0 + [\bfu_1,\dots,\bfu_i], \qquad i=1,2,\dots.$$

\no Sottolineando ancora una volta che questa propriet\`a potrebbe
non essere pi\`u verificata in aritmetica finita, abbiamo che essa
pu\`o essere scritta, equivalentemente, come
$$E(\bfx_i) \equiv (\bfx^*-\bfx_i)^\top A(\bfx^*-\bfx_i) =
\min_{\bfx\in S_i}\, (\bfx^*-\bfx)^\top A(\bfx^*-\bfx).$$

\no Tuttavia, dalla (\ref{x_i}) e dal Lemma~\ref{Krylov} segue che:

\begin{itemize}

\item $\bfe_i = \bfe_{i-1} -\lam_i \bfu_i = \dots = \bfe_0
-\sum_{j=1}^i \lam_j \bfu_j$;

\item $[\bfu_1,\dots,\bfu_i] =
[\bfr_0,A\bfr_0,\dots,A^{i-1}\bfr_0]\equiv[A\bfe_0,\dots,A^i\bfe_0]$.

\end{itemize}

\no Pertanto, si conclude che
$$\bfe_i = \bfe_0 + \sum_{j=1}^i c_j A^j \bfe_0 \equiv
p_i(A)\bfe_0, \qquad\mbox{con}\qquad p_i\in\Pi_i,~p_i(0)=1.$$

\no Osserviamo che (vedi Sezione~\ref{polmat}), se
$$A=Q\Lambda Q^\top  \equiv Q\pmatrix{ccc} \lam_1(A)
\\&\ddots\\ &&\lam_n(A)\endpmatrix Q^\top ,$$

\no con $Q$ matrice ortogonale e $$0<\lam_{min}\equiv\lam_1(A)
\le\dots\le\lam_n(A)\equiv \lam_{max},$$

\no allora $p_i(A) = Qp_i(\Lambda) Q^\top $~ e ~$E(\bfx_0)=\bfe_0^\top Q\Lambda
Q^\top \bfe_0$. Segue quindi che

\begin{eqnarray*}
E(\bfx_i) &=& \bfe_0^\top  p_i(A)^2A \bfe_0 = (Q^\top \bfe_0)^\top 
p_i(\Lambda)^2\Lambda(Q^\top \bfe_0)\\[1mm]
&\equiv& \min_{p\in\Pi_i,\,p(0)=1} (Q^\top \bfe_0)^\top p(\Lambda)^2 \Lambda(Q^\top \bfe_0) \\[1mm]
& \le& (Q^\top \bfe_0)^\top \Lambda(Q^\top \bfe_0)\min_{p\in\Pi_i,\,p(0)=1}\max_{~\lam\sigma(A)}p(\lam)^2\\[1mm]
&\equiv& E(\bfx_0) \min_{p\in\Pi_i,\,p(0)=1}~ \max_{\lam\in\sigma(A)}
p(\lam)^2
\le E(\bfx_0) \min_{p\in\Pi_i,\,p(0)=1}~ \max_{\lam\in[\lam_{min},\lam_{max}]}
p(\lam)^2 \\[1mm]
&=& E(\bfx_0) \left(\min_{p\in\Pi_i,\,p(0)=1}~ \max_{\lam\in[\lam_{min},\lam_{max}]}
|p(\lam)|\right)^2.
\end{eqnarray*}

\no La soluzione dell'ultimo problema di minimassimo \`e noto
essere data da (vedi propriet\`a C8 in Sezione~\ref{cebypol})
$$p(z) = \frac{T_i(\mu(z))}{T_i(\mu(0))}, \qquad \mu(z) =
\frac{\lam_{max}+\lam_{min}-2z}{\lam_{max}-\lam_{min}},$$

\no e, pertanto, tenendo conto anche della propriet\`a C3, si
ottiene $$E(\bfx_i)\le T_i(\mu(0))\/^{-2} E(\bfx_0).$$

\no Considerando che, detto $\kappa(A) = \lam_{max}/\lam_{min}\equiv \rho^2$ il numero di
condizionamento della matrice $A$,
$$\bar{z}\equiv\mu(0) = \frac{\lam_ {max}+\lam_{min}}{\lam_{max}-\lam_{min}} = \frac{\rho^2+1}{\rho^2-1}>1,$$

\no e tenendo conto della propriet\`a C7 in Sezione~\ref{cebypol}, si ottiene, infine,
\begin{eqnarray*}
E(\bfx_i) &\le& 4 \left( \bar{z}+\sqrt{\bar{z}^2-1} \right)^{-2i}E(\bfx_0)
= 4 \left( \bar{z}-\sqrt{\bar{z}^2-1} \right)^{2i}E(\bfx_0)\\[1mm]
&=& 4\left( \frac{\rho^2+1-2\rho}{\rho^2-1} \right)^{2i} E(\bfx_0)
= 4\left( \frac{\rho-1}{\rho+1} \right)^{2i} E(\bfx_0) \\[1mm] &\equiv&
4\left( \frac{\sqrt{\kappa(A)}-1}{\sqrt{\kappa(A)}+1} \right)^{2i}
E(\bfx_0).
\end{eqnarray*}

\begin{oss}
 Da questa espressione si evince che la convergenza del metodo
dei gradienti coniugati diviene via, via pi\`u lenta, al crescere
del numero di condizionamento della matrice $A$. Osserviamo che un
grande numero di condizionamento significa anche avere un criterio di
arresto pi\`u restrittivo se, ad esempio, si utilizza la strategia
(\ref{stop2}) (vedi anche l'Esercizio~\ref{c_stop}).
Si conclude che, se $A$ fosse ben condizionata,
ovvero $\kappa(A)$ fosse relativamente ``piccolo'', si avrebbe un
vantaggio, in termini di costo computazionale del metodo.
Questi argomenti costituiscono la motivazione principale per
quanto andremo a trattare nel prossimo capitolo.
\end{oss}

\begin{figure}[t]
\begin{center}
\includegraphics[width=13cm,height=10cm]{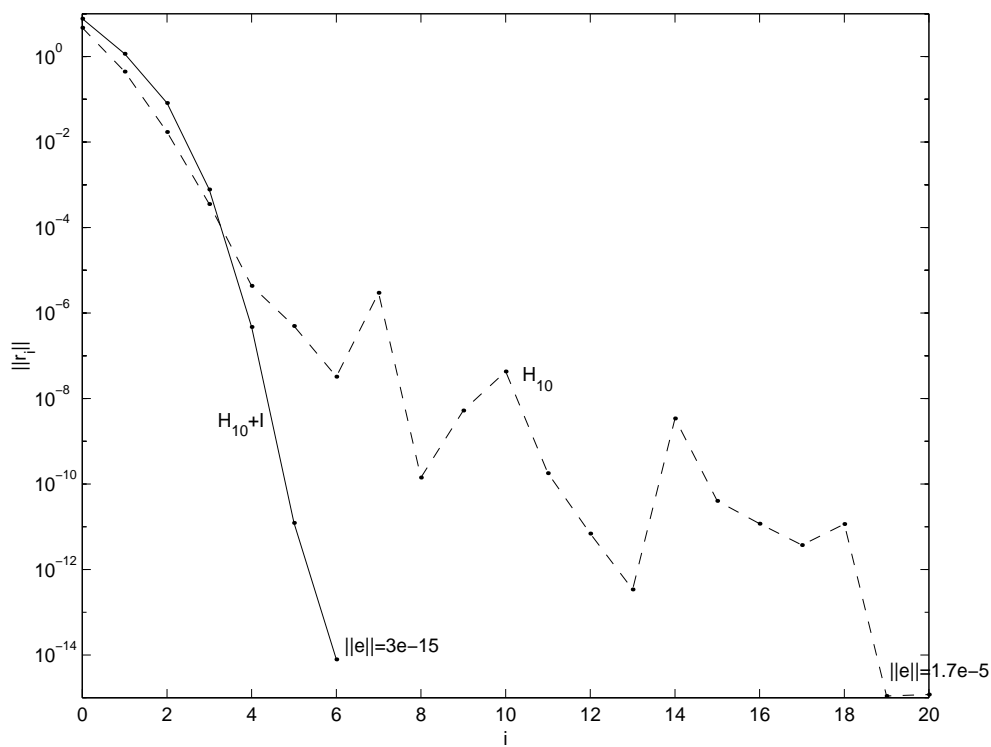}
\caption{\protect\label{fig5} Metodo del gradiente coniugato
applicato al problema (\protect\ref{sistsdp}) nei casi
$A=H_{10}$ e $A=H_{10}+I$.}
\end{center}
\end{figure}

Come esempio, si consideri l'applicazione del metodo del gradiente
coniu\-gato al sistema lineare (\ref{sistsdp}), di dimensione $n=10$,
e in cui il termine noto \`e definito in modo da avere come soluzione
il vettore $\bfx^*=(1,\dots,1)^\top $, nei seguenti due casi:
\begin{enumerate}
\item $A=H_{10}\equiv\left((i+j-1)^{-1}\right)$, la matrice di Hilbert
di dimensione 10, gi\`a introdotta nella precedente Sezione~\ref{stopcg};

\item $A=H_{10}+I$.
\end{enumerate}

\no Nel caso 1, il numero di condizione di $A$ \`e dell'ordine di
$10^{13}$, mentre nel caso 2, esso vale circa 2.8. Ne consegue
che, nel caso 2, il metodo dei gradienti coniugati converger\`a
pi\`u rapidamente che nel caso 1. Inoltre, a parit\`a di norma del
residuo, l'errore sulla soluzione sar\`a minore nel caso 2,
rispetto al caso 1. Queste affermazioni trovano piena conferma
sperimentale, come si evince dalla Figura~\ref{fig5} in cui sono
riportati i grafici delle norme dei residui rispetto all'indice di
iterazione. Nel caso 2, dopo 6 iterazioni il residuo ha norma pari
a circa $8\cdot10^{-15}$ e il corrispondente errore \`e dello
stesso ordine. Nel caso 1 sono richieste 19 iterazioni per avere
un residuo di norma comparabile, ma l'errore risulta essere assai
pi\`u grande. Infatti, anche se alla iterazione 20 il residuo ha
norma di poco superiore a $10^{-15}$, il corrispondente errore ha
norma dell'ordine di $10^{-5}$. In entrambi i casi, il vettore
iniziale utilizzato \`e il vettore nullo.

\`E da sottolineare come, nel caso 1, la propriet\`a di
terminazione finita del metodo dei GC (vedi Corollario~\ref{cor4.5})
 non sia pi\`u valida, a causa dell'utilizzo dell'aritmetica finita.

%
%
\chapter{Precondizionamento}\label{cap5}

{\em
Abbiamo visto che la convergenza del metodo del gradiente
coniugato dipende crucialmente dal numero di condizione della
matrice $A$. Quando questo \`e molto grande, la convergenza
potrebbe essere troppo lenta o, a causa degli errori di macchina,
addirittura  non esserci. In questa evenienza, \`e raccomandabile
una riformulazione del problema mediante {\em precondizionamento}
che, spesso, \`e in grado di ricondurre il costo computazionale
entro limiti accettabili. In questo capitolo tratteremo il
precondizionamento di sistemi lineari con matrice sdp, sebbene
alcune delle tecniche possano essere estese a casi pi\`u generali.
Non analizzeremo tutte le implementazioni di questo approccio, ma
solo alcune tra quelle pi\`u significative e generali.}

\section{Fattorizzazione $CAU=UT$}\label{CAUeqUT}

Per ovviare ai problemi derivanti dall'avere $\kappa(A)\gg1$,
invece di risolvere il problema (\ref{sistsdp}), si risolve
l'equivalente {\em sistema precondizionato}
\begin{equation}\label{CA}
CA\bfx = C\bfb, \qquad C\approx A^{-1}.
\end{equation}

\no In questo modo, ci si aspetta di avere
$\kappa(CA)\ll\kappa(A)$. \`E possibile definire il sistema
precondizionato anche nella forma pi\`u generale
$$C_1AC_2\,\, C_2^{-1} \bfx = C_1\bfb, \qquad C_2C_1\approx
A^{-1},$$

\no tuttavia nel seguito esamineremo, per brevit\`a,
solo la forma (\ref{CA}). Assume\-remo, inoltre, che la matrice
$C$ sia, al pari di $A$, sdp.
Osserviamo che la scelta $C=A^{-1}$ darebbe $CA=I$, che \`e
perfettamente condizionata. Tuttavia, tale scelta non \`e
proponibile, in quanto equivalente alla risoluzione di
(\ref{sistsdp}) con un metodo diretto.
Con queste premesse, distinguiamo i seguenti due casi.

\begin{description}

\item[$CA$ sdp:~] in tal caso, il sistema (\ref{CA}) \`e
formalmente equivalente a (\ref{sistsdp}). Possiamo quindi utilizzare la
matrice $B=CA$ ed applicare formalmente tutto quello che si \`e
visto nel precedente Capitolo~\ref{cap4}, sostituendo la matrice
$CA$ alla matrice $A$;

\item[$CA$ non sdp:~] in questo caso \`e richiesta una
ulteriore analisi, che andiamo a fare nel seguito.

\end{description}

Nel caso in cui $CA$ non sia sdp (bench\`e lo siano $C$ ed $A$),
si ricerca la seguente fattorizzazione, alternativa alla
(\ref{fattA}):
\begin{equation}\label{fattCA}
CAU=UT, \qquad U^\top BU = D,
\end{equation}

\no in cui le matrici $U,T$ e $D$ hanno la stessa struttura vista
in (\ref{fattA}). Tuttavia, per tale fattorizzazione, il requisito
CG2 visto in Sezione~\ref{AUeqUT} viene rimpiazzato da quello seguente
(che assumeremo soddisfatto nel seguito):

\bigskip
\no{\bf PCG:~} $ACB = BCA$.
\bigskip

Si verifica facilmente che la fattorizzazione
(\ref{fattCA}) risulta essere definita ancora dalle espressioni
(\ref{d_i})--(\ref{b_i}), previa la sostituzione
\begin{equation}\label{AeqCA}
A\leftarrow CA,
\end{equation}

\no Sempre con la (\ref{AeqCA}), continuano a valere i risultati
dei Teoremi~\ref{cor4.1} e \ref{teo4.1}.

Andando avanti, la scelta (\ref{AeqB}) soddisfa il requisito PCG,
consentendo di ricercare la nuova approssimazione ancora nella
forma (\ref{x_i}), in modo da soddisfare le propriet\`a di
ortogonalit\`a (\ref{rort}). Questo porta alla stessa relazione di
ricorrenza (\ref{r_i})-(\ref{lami}) per i vettori residuo ed i
risultati del Teorema~\ref{teo4.2} e del Corollario~\ref{cor4.2}
continuano a valere.

A questo punto, si introducono i {\em residui precondizionati}
\begin{equation}\label{s_i}
\bfs_i = C\bfr_i, \qquad i=0,1,\dots,
\end{equation}

\no e la scelta (\ref{u1}) del vettore iniziale viene rimpiazzata
da
$$
\bfu_1 = \bfs_0 \equiv C\bfr_0\equiv C(\bfb-A\bfx_0),
$$

\no che definisce il metodo dei {\em gradienti coniugati
precondizionati}. Tenendo conto della (\ref{AeqCA}) e sostituendo
i vettori residuo $\{\bfr_i\}$ con i corrispondenti residui
precondizionati (\ref{s_i}), i risultati del Lemma~\ref{Krylov} (vedi
Esercizio~\ref{ex5.*}) e dei Corollari~\ref{cor4.4} e \ref{cor4.5} si
generalizzano al presente caso. Il Corollario~\ref{cor4.3} si riformula,
invece, come segue.

\begin{cor} I vettori dei residui $\{\bfr_i\}$ e quelli dei
residui precondizionati $\{\bfs_i\}$ sono tra loro bi-ortogonali,
ovvero:
\begin{equation}\label{siri}
\bfs_j^\top \bfr_i = 0, \qquad i\ne j. \end{equation}\end{cor}

\proof Per la (\ref{s_i}), ed essendo $C$ sdp, basta dimostrare
la tesi per $i>j$. Questa segue, quindi, dalla (\ref{rort}),
ricordando che, in questo caso,
$[\bfu_1,\dots,\bfu_i]\equiv[\bfs_0,\dots,\bfs_{i-1}]$.\QED

\begin{oss} La (\ref{siri}) pu\`o essere riscritta, in virt\`u
della (\ref{s_i}), come $$\bfr_j^\top C\bfr_i = 0, \qquad i\ne j.$$
Ovvero i residui $\{\bfr_i\}$ sono tra loro $C$-ortogonali. Questo
conferma la terminazione (in aritmetica esatta) in al pi\`u $n$
passi.
\end{oss}

Per ultimo, continua a valere, sempre considerando la
(\ref{AeqCA}) ed i residui precondizionati al posto di quelli non
precondizionati, una estensione del Teorema~\ref{teo4.3} (vedi
Esercizio~\ref{ex5.**}), il quale permette, alfine, di derivare
l'Algoritmo~\ref{algo5} per il metodo dei gradienti coniugati precondizionati.

\begin{figure}[t]
\begin{algo}\label{algo5} Metodo dei gradienti coniugati precondizionati.\\ \rm
\fbox{\parbox{12cm}{
\begin{itemize}

\nulit sia assegnata una approssimazione iniziale $\bfx_0$

\nulit inizializzo $\bfr_0 = \bfb-A\bfx_0, \bfp_1 = \bfs_0\equiv
C\bfr_0$ e $\eta_0=\bfr_0^\top \bfs_0$

\nulit per i=1,\dots,n:

\begin{itemize}

\nulit se $\eta_{i-1}=0, \bfx^*=\bfx_{i-1}$, esci, fine se

\nulit $\bfv_i=A\bfp_i$

\nulit $\hat{d}_i = \bfp_i^\top \bfv_i$

\nulit $\hat\lam_i = \eta_{i-1}/\hat{d}_i$

\nulit $\bfx_i = \bfx_{i-1}+\hat\lam_i\bfp_i$

\nulit $\bfr_i = \bfr_{i-1}-\hat\lam_i\bfv_i$

\nulit $\bfs_i = C\bfr_i$

\nulit $\eta_i = \bfr_i^\top \bfs_i$

\nulit $\mu_i = \eta_i/\eta_{i-1}$

\nulit $\bfp_{i+1} = \bfs_i+\mu_i\bfp_i$

\end{itemize}

\nulit fine per

\nulit se $i=n, \bfx^*=\bfx_n$, fine se
\end{itemize}

}}\end{algo}
\end{figure}

\begin{oss} Rispetto all'Algoritmo~\ref{algo4}, l'Algoritmo~\ref{algo5}
richiede un costo computazionale aggiuntivo dato da:
\begin{itemize}

\item un vettore per il residuo precondizionato corrente e lo
spazio per $C$, in termini di occupazione di memoria (osserviamo, tuttavia,
che il vettore $\bfs_i$ potrebbe essere riscritto sul vettore $\bfv_i$);

\item il calcolo del residuo precondizionato,
in termini di operazioni per ite\-razione.

\end{itemize}
Questo costo aggiuntivo in memoria ed operazioni per iterata \`e
tuttavia giu\-stificato dalla aspettativa di eseguire molte meno
iterazioni, come vedremo, rispetto al metodo non precondizionato.
\end{oss}

\begin{oss} Osserviamo che il criterio di arresto
nell'Algoritmo~\ref{algo5} controlla la $C$-norma del residuo,
laddove, nell'Algoritmo~\ref{algo4}, se ne controllava la norma
euclidea. Anche ora (vedi Sezione~\ref{stopcg}) il criterio di arresto pi\`u
appropriato non \`e $\eta_{i-1}=0$, ma $\eta_{i-1}\le tol$, per una prefissata
tolleranza $tol$.\end{oss}

\begin{oss} In genere, come vedremo nel seguito, il residuo
precondizionato \`e calcolato risolvendo il sistema lineare
\begin{equation}\label{M}
M\bfs_i = \bfr_i, \qquad M = C^{-1}.
\end{equation}

\no invece che eseguendo il prodotto (\ref{s_i}). Pertanto, il
{\em precondizionatore} $M$ deve essere caratterizzato da una
semplice ed economica risoluzione del sistema lineare (\ref{M}).
\end{oss}

\begin{eser}\label{ex5.*} Dimostrare che, per il metodo dei
gradienti coniugati precondizionati, se $\bfx_i\ne\bfx^*$
allora $$[\bfu_1,\dots,\bfu_{i+1}] =
[\bfu_1,CA\bfu_1,\dots,(CA)^i\bfu_1] = [\bfs_0,\dots,\bfs_i].$$
\end{eser}

\begin{eser}\label{ex5.**} Dimostrare che, per il metodo dei
gradienti coniugati precondizionati, se $\bfx_i\ne\bfx^*$
allora $\bfu_{i+1}\in[\bfu_i,\bfs_i]$.
\end{eser}

\begin{eser} Riscrivere l'Algoritmo~\ref{algo5} tenendo conto
della (\ref{M}) e delle sole strutture di memoria effettivamente
richieste nella sua implementazione. Riformularlo, inoltre, come
metodo iterativo (e non semi-iterativo), inserendo un idoneo
criterio di arresto.\end{eser}

Nelle prossime sezioni, vedremo alcune tecniche per derivare
precondizio\-natori efficienti.

\section{Fattorizzazioni incomplete}\label{IC}

\`E noto che, se $A$ \`e sdp, allora essa \`e fattorizzabile nella
forma \begin{equation}\label{ldl}A = LDL^\top ,\end{equation}

\no (vedi Sezione~\ref{fattLUmat}) dove $$D=\pmatrix{ccc} d_1\\ &\ddots \\&&d_n\endpmatrix,
\qquad d_i>0,\quad i=1,\dots,n,$$

\no e $$L = \pmatrix{ccc} l_{11}\\ \vdots &\ddots\\
l_{n1}&\dots&l_{nn}\endpmatrix, \qquad l_{ii}=1,\quad
i=1,\dots,n.$$

\no Tuttavia, anche se $A$ \`e una matrice sparsa, non \`e detto
che il suo fattore $L$ sia sparso a sua volta. Ad esempio, la
matrice (\ref{test}) in Esercizio~\ref{ex1} \`e una matrice a
banda con ampiezza di banda $N$, ma con solo 5 diagonali
significative. Il suo fattore $L$, bench\`e preservi l'ampiezza di
banda $N$, ha tutte le $N$ sottodiagonali non nulle. Pertanto,
quando $N$ \`e grande, l'occupazione di memoria diviene cospicua.
Per questo motivo, si ricercano fattorizzazioni approssimate della
matrice $A$, che possano essere utilizzate come precondizionatori.
Ovvero, si ricerca una decomposizione del tipo
\begin{equation}\label{ldlr}
A = \tL\tD\tL^\top -R \equiv M -R,
\end{equation}

\no dove $R$ \`e la matrice di errore, $M$ \`e il
precondizionatore, $$\tD=\pmatrix{ccc}\td_1\\
&\ddots\\&&\td_n\endpmatrix,$$

\no e $\tL$ \`e triangolare inferiore con elementi diagonali uguali ad
1, ma tipicamente molto pi\`u sparsa del corrispondente fattore
$L$ in (\ref{ldl}). Poich\`e siamo inte\-ressati ad avere $M=C^{-1}$
sdp, ricerchiamo la fattorizzazione approssimata in modo tale da
avere $\td_i>0,i=1,\dots,n$. Vediamo sotto quali condizioni questo
si verifica.

Vi sono vari modi per definire la fattorizzazione approssimata
(\ref{ldlr}). Ad esempio, si possono memorizzare solo gli elementi extra-diagonali
del fattore $\tL$ che hanno modulo pi\`u grande di una prefissata tolleranza.
Un'alternativa, tra le pi\`u utilizzate, \`e quella di fissare {\em a priori}
la struttura di sparsit\`a del fattore $\tL$. Formalmente, se
$$G_n = \left\{ (i,j)\in\{1,\dots,n\}^2 : i >j\right\},$$

\no si considera un opportuno $\tG_n\subseteq G_n$. Il sottoinsieme di indici $\tG_n$
definisce la fattorizzazione approssimata (\ref{ldlr}) tale che:
\begin{equation}\label{tL}
\tL = \pmatrix{ccc} \ell_{11} \\ \vdots &\ddots\\ \ell_{n1} &\dots
&\ell_{nn}\endpmatrix, \qquad\mbox{dove:}\quad
\left\{\begin{array}{ll} \ell_{ii}=1, & i=1,\dots,n,\\ \\
\ell_{ij}=0, &(i,j)\not\in\tG_n.\end{array}\right.
\end{equation}

\begin{oss} Se $\tG_n= G_n$, ad esso corrisponde
la fattorizzazione esatta (\ref{ldl}). D'altronde, al sottoinsieme
banale, $\tG_n=\emptyset$, corrisponde $\tL=I$, mentre $\tD\equiv M$
contiene gli elementi diagonali di $A$. Evidentemente, in
entrambi i casi il corrispondente precondizionatore $M$ \`e sdp.
\end{oss}

Studieremo le condizioni per ottenere $M$ sdp nel caso, abbastanza
importante, in cui $A$ \`e una matrice di Stieltjes, ovvero una
$M$-matrice simmetrica (e, quindi, definita positiva). Si premette
il seguente risultato.

\begin{lem}\label{ldlmmat} Se $A$ \`e una matrice di Stieltjes, allora
il fattore $L$ in (\ref{ldl}) \`e una $M$-matrice.\end{lem}

\proof \`E sufficiente dimostrare che $l_{ij}\le0$, se $i>j$.
Dimostriamo per induzione sulla dimensione, $n$, di A. Per $n=1$,
$L_1=(1)$ e, pertanto, la tesi \`e vera. Supponiamola vera per
$n-1$ e dimostriamo per $n$. Si ha:
\begin{eqnarray}\nonumber
A&\equiv& A_n =\pmatrix{cc} A_{n-1} &\bfa_n\\ \bfa_n^\top 
&a_{nn}\endpmatrix = \pmatrix{cc} L_{n-1}\\ \bfl_n^\top 
&1\endpmatrix\pmatrix{cc} D_{n-1}\\&d_n\endpmatrix\pmatrix{cc}
L_{n-1}^\top  &\bfl_n\\ &1\endpmatrix\\[1mm] &=& \pmatrix{cc}
L_{n-1}D_{n-1}L_{n-1}^\top  & L_{n-1}D_{n-1}\bfl_n\\
\bfl_n^\top D_{n-1}L_{n-1}^\top  &d_n+\bfl_n^\top D_{n-1}\bfl_n\endpmatrix.
\label{dl}
\end{eqnarray}

\no Pertanto, $$L_n =\pmatrix{cc} L_{n-1}\\ \bfl_n^\top 
&1\endpmatrix.$$

\no Dall'ipotesi di induzione, segue che $L_{n-1}$ \`e una
$M$-matrice. Rimane quindi da dimostrare che $\bfl_n\le0$. Dalla
(\ref{dl}) si ottiene, infatti, $$0\ge\bfa_n \equiv
L_{n-1}D_{n-1}\bfl_n.$$

\no essendo $D_{n-1}$ con elementi diagonali positivi, e $L_{n-1}$
una $M$-matrice, segue quindi che $$\bfl_n =
D_{n-1}^{-1}L_{n-1}^{-1}\bfa_n\le0.~\QED$$

\medskip
\begin{eser} Dall'equazione (\ref{dl}) ricavare le seguenti
espressioni per gli elementi della fattorizzazione (\ref{ldl}):
\begin{eqnarray}\nonumber
d_j &=& a_{jj} -\sum_{k=1}^{j-1} l_{jk}^2d_k, \qquad\qquad j = 1,\dots,n,\\
l_{jj} &=& 1,
\label{chol}\\ \nonumber
l_{ij} &=& \frac{a_{ij}-\sum_{k=1}^{j-1}l_{ik}l_{jk}d_k}{d_j},
\qquad i=j+1,\dots,n.
\end{eqnarray}
\end{eser}

\medskip
\begin{eser} Generalizzare l'equazione (\ref{dl}) al caso della
fattorizzazione approssimata (\ref{ldlr})-(\ref{tL}). Verificare
che gli elementi dei fattori sono formalmente dati da:
\begin{eqnarray}\nonumber
\td_j &=& a_{jj} -\sum_{k=1}^{j-1} \ell_{jk}^2\td_k, \qquad\qquad
j = 1,\dots,n,\\ \ell_{jj} &=& 1, \label{cholapp}\\ \nonumber
\ell_{ij} &=&\left\{\begin{array}{cl}
\frac{\displaystyle a_{ij}-\sum_{k=1}^{j-1}\ell_{ik}\ell_{jk}\td_k}{\displaystyle \td_j},
&\quad (i,j)\in\tG_n,\\ \\ 0,
&\quad (i,j)\not\in\tG_n. \end{array}\right.\nonumber
\end{eqnarray}
\end{eser}

Possiamo ora enunciare il risultato di stabilit\`a riguardo alle
fattorizzazioni approssimate di matrici di Stieltjes.

\begin{teo} Se $A$ \`e una matrice di Stieltjes e $\tG_n\subseteq G_n$,
allora (vedi (\ref{ldl})-(\ref{ldlr})) ~$\tD\ge\tL\tD\ge LD$~ e
~$\tD\ge D$.\end{teo}

\proof Dimostriamo per induzione sulla dimensione, $n$, di $A$.
Per $n=1$ la tesi \`e banalmente vera. Supponiamola vera per $n-1$
e dimostriamo per $n$. Infatti, la fattorizzazione esatta sar\`a
data da (\ref{dl}), mentre per i fattori approssimati si avr\`a
$$\tL_n = \pmatrix{cc} \tL_{n-1}\\ \bfel_n^\top  &1\endpmatrix,\qquad
  \tD_n = \pmatrix{cc} \tD_{n-1} \\ &\td_n\endpmatrix,$$

\no con ~$\tD_{n-1}\ge\tL_{n-1}\tD_{n-1}\ge L_{n-1}D_{n-1}$~ e
~$\tD_{n-1}\ge D_{n-1}(\ge0)$, per ipotesi di induzione. Gli
elementi delle due fattorizzazioni sono dati da (\ref{chol}) e
(\ref{cholapp}), rispettivamente. Rimane, quindi, da dimostrare
che
\begin{eqnarray}\nonumber
0&\ge&\bfel_n^\top \tD_{n-1}\equiv(\ell_{n1}\td_1,\dots,\ell_{n,n-1}\td_{n-1})\\[2mm]&\ge&
(l_{n1}d_1,\dots,l_{n,n-1}d_{n-1})\equiv\bfl_n^\top D_{n-1},\label{ellgel}
\end{eqnarray}

\no in quanto l'altra diseguaglianza discende dal fatto che la
precedente implica
\begin{eqnarray*}
0&\le&\bfel_n^\top \tD_{n-1}\bfel_n =
(-\bfel_n^\top \tD_{n-1})\tD_{n-1}^{-1}(-\bfel_n^\top \tD_{n-1})^\top \\[2mm]
&\le& (-\bfl_n^\top D_{n-1})D_{n-1}^{-1}(-\bfl_n^\top D_{n-1})^\top  =
\bfl_n^\top D_{n-1}\bfl_n,
\end{eqnarray*}

\no da cui segue che $$\td_n = a_{nn}-\bfel_n^\top \tD_{n-1}\bfel_n \ge
a_{nn}-\bfl_n^\top D_{n-1}\bfl_n=d_n.$$

\no Dimostriamo la (\ref{ellgel}), ovvero che $\ell_{nj}\td_j\ge l_{nj}d_j$,
 per induzione su $j$. Per $j=1$ ho due possibilit\`a:

\begin{itemize}

\item $(n,1)\not\in\tG_n$: pertanto $\ell_{n1}\td_1 = 0 \ge
a_{n1} \equiv l_{n1}d_1$;

\item $(n,1)\in\tG_n$: in questo caso $0\ge\ell_{n1}\td_1 \equiv
a_{n1}\equiv l_{n1}d_1$.

\end{itemize}

\no Supponiamo ora vera la tesi per $k\le j-1$, e dimostramo per
$k=j$. Ancora una volta, posso avere i seguenti due casi:

\begin{itemize}

\item $(n,j)\not\in\tG_n$: pertanto $\ell_{nj}\td_j = 0 \ge
l_{nj}d_j$, dal Lemma~\ref{ldlmmat};

\item $(n,j)\in\tG_n$: in questo caso dall'ipotesi di induzione (sia quella
su $n$ che quella su $j$), e tenendo conto delle
(\ref{chol})-(\ref{cholapp}), segue che
\begin{eqnarray*}
0&\ge& \td_j\ell_{nj} = a_{nj}-\sum_{k=1}^{j-1} \ell_{nk}\ell_{jk}\td_k
= a_{nj}-\sum_{k=1}^{j-1} (-\ell_{nk}\td_k)(-\ell_{jk}\td_k)\td_k^{-1} \\
&\ge& a_{nj}-\sum_{k=1}^{j-1} (-l_{nk}d_k)(-l_{jk}d_k)d_k^{-1} =
 a_{nj}-\sum_{k=1}^{j-1} l_{nk}l_{jk}d_k =
d_j l_{nj}.\end{eqnarray*}

\no Considerando che $\td_j\ge d_j>0$, segue quindi che $0\ge
\ell_{nj}\ge l_{nj}$. \QED
\end{itemize}

\begin{cor} $I\ge\tL\ge L$ e, pertanto, $\tL$ \`e una $M$-matrice.\end{cor}

\begin{cor} La fattorizzazione approssimata di $A$ definita da un
qualunque sottoinsieme di $G_n$ \`e stabile e il corrispondente
precondizionatore \`e sdp.\end{cor}

\begin{cor} Se $M\equiv C^{-1}$ \`e definito come in (\ref{ldlr}),
allora $$A^{-1}\ge M^{-1}=C\ge 0.$$\end{cor}

\section{Fattorizzazione incompleta di Cholesky}\label{IC1}

Se $A=(a_{ij})\in\RR^{n\times n}$, un modo abbastanza ovvio di definire il
sottoinsieme $\tG_n$ \`e il seguente: $$\tG_n = \left\{ (i,j): i>j ~\wedge~ a_{ij}\ne0
\right\},$$

\no ovvero, imponiamo che $\tL$ abbia la stessa struttura di
sparsit\`a della parte triangolare inferiore di $A$. Questa scelta
d\`a origine alla {\em fattorizzazione incompleta di Cholesky},
spesso denotata con {\em IC (incomplete Cholesky)}. Deri\-viamo la
sua espressione nel caso in cui $A$ sia una matrice di Stieltjes
con una struttura di sparsit\`a assimilabile a quella della
matrice del problema (\ref{test}) dell'Esercizio~\ref{ex1}. In tal
caso, la matrice $A$ \`e pentadiagonale con ampiezza di banda $N$:
\begin{equation}\label{penta}
A = \pmatrix{cccccc}
a_1 &b_2    &      &c_{N+1}\\
b_2 &\ddots &\ddots&       &\ddots\\
    &\ddots &\ddots &\ddots&      &c_n\\
c_{N+1} &   &\ddots &\ddots&\ddots\\
       &\ddots&     &\ddots &\ddots&b_n\\
       &      &c_n  &       &b_n   &a_n\endpmatrix, \qquad a_i>0,~b_i,c_i\le0.
       \end{equation}

\no Pertanto, $\tG_n = \{(i,i-1), (i,i-N)\}$ e, conseguentemente, dalle
(\ref{cholapp}) si ottengono le seguenti espressioni (ponendo, per convenzione,
$b_1=c_1=\dots=c_N=0$):
\begin{eqnarray}\nonumber
\ell_{i,i-1} &=& \frac{ b_i -\sum_{k=1}^{i-2}
\ell_{ik}\ell_{i-1,k}\td_k}{\td_{i-1}} \equiv
\frac{b_i}{\td_{i-1}},\quad i = 2,\dots,n,
\\[1mm] \label{ellij}
\ell_{i,i-N} &=& \frac{ c_i -\sum_{k=1}^{i-N-1}
\ell_{ik}\ell_{i-N,k}\td_k}{\td_{i-N}} \equiv \frac{c_i}{\td_{i-N}},
\quad i = N+1,\dots,n,~~~~~~\\[1mm] \nonumber
\td_i &=& a_i -\sum_{k=1}^{i-1}\ell_{ik}^2\td_k \equiv
a_i-\frac{b_i^2}{\td_{i-1}}-\frac{c_i^2}{\td_{i-N}}, \quad
i=1,\dots,n.
\end{eqnarray}

\medskip
\begin{oss} Vista la particolare forma degli elementi (\ref{ellij}) della
fattorizzazione incompleta, nel caso della matrice (\ref{penta})
risulta essere conveniente riscrivere la fattorizzazione
$\tL\tD\tL^\top $ nella forma equivalente
$$(\tL\tD) \tD^{-1} (\tL\tD)^\top ,$$

\no in quanto si verifica facilmente che $$(\tL\tD)
=\pmatrix{cccccc} \td_1 \\ b_2 &\td_2\\
    &\ddots &\ddots \\
c_{N+1} &   &\ddots &\ddots\\
       &\ddots&     &\ddots &\ddots\\
       &      &c_n  &       &b_n   &\td_n\endpmatrix,$$

\no la cui parte strettamente triangolare coincide con quella di
$A$. Gli  elementi diagonali sono, invece, quelli omologhi della
matrice $\tD$. In tal modo, l'occupazione di memoria richiesta per
il precondizionatore si riduce ad un solo vettore di lunghezza
$n$, in cui memorizzare i $\{\td_i\}$.
\end{oss}

\medskip
\begin{eser} Calcolare il costo computazionale per risolvere il
sistema $$M\bfs =\bfr$$ nel caso della fattorizzazione IC della
matrice (\ref{penta}), e quello del corrispondente matvec
$A\bfp$.\end{eser}

\medskip
\begin{eser}\label{ex5.4} Verificare che, nel caso della fattorizzazione
incompleta (\ref{ellij}) della matrice (\ref{penta}), la
matrice di errore $R$ in (\ref{ldlr}) \`e data da
\begin{equation}\label{R1}
R=\pmatrix{cccc}
   &r_N\\
   &    &\ddots\\
r_N&    &      &r_n\\
   &\ddots\\
   &     &r_n\endpmatrix,
\end{equation}
con
\begin{equation}\label{R2}
r_N=0,\qquad r_i = \frac{b_{i-N+1}c_i}{\td_{i-N}},\quad
i=N+1,\dots,n.
\end{equation}
Osserviamo che $r_i\ge0$, se $A$ \`e, come abbiamo assunto, una
matrice di Stieltjes.
\end{eser}

\section{Fattorizzazione incompleta modificata}\label{MIC}

La bont\`a del precondizionatore $M$ (o, equivalentemente, di
$C=M^{-1}$) determina quanto la {\em matrice precondizionata},
$$CA \equiv M^{-1}A$$

\no ``assomiglia'' alla matrice identit\`a $I$. Un modo per
misurare questo, \`e ana\-lizzarne lo spettro: pi\`u questo avr\`a
gli autovalori addensati intorno a $\lam=1$, tanto migliore sar\`a
il precondizionatore. Dalla (\ref{ldlr}), segue che $$M^{-1}A = I
- M^{-1}R.$$

\no Inoltre, osserviamo che, essendo $M^{-1}$ sdp, allora sar\`a
fattorizzabile nella forma $M^{-1} = G^\top G$, con $G=(\tL\tD^{\frac{1}2})^{-1}$
(fattorizzazione di Cholesky). Pertanto, $$M^{-1}R = (G^\top G)R \sim GRG^\top .$$

\no Tenendo conto del fatto che $R$ \`e simmetrica, si
conclude che gli autovalori di $M^{-1}R$ sono reali ed aventi gli
stessi segni di quelli di $R$.

Osserviamo che, essendo $M^{-1}A \sim GAG^\top $ che \`e definita
positiva, allora il numero di condizionamento di quest'ultima
matrice sar\`a dato da $$\frac{1-\mu_{min}}{1-\mu_{max}},$$ dove
$\mu_{min}$ e $\mu_{max}$ sono, rispettivamente, il pi\`u piccolo
ed il pi\`u grande autovalore di $M^{-1}R$. Ricordiamo che
$1>\mu_{max}$, sebbene quest'ultimo potrebbe essere, in alcuni
casi, assai vicino a $1$, rendendo la matrice malcondizionata.
Tuttavia, se riuscissimo a modificare la fattorizzazione
incompleta in modo da garantire che
\begin{equation}\label{neg}\sigma(R)\subset\RR^-,
\end{equation}

\no allora seguirebbe che $\mu_{max}\le0$. Pertanto,
l'autovalore pi\`u piccolo della matrice precondizionata sarebbe
non minore di $\lam=1$, migliorandone gene\-ralmente il
condizionamento. Questa modifica, nel caso della fattorizzazione
IC, ne origina la versione {\em modificata}, denominata {\em MIC
(modified IC)}.

Ricaviamo la espressione della fattorizzazione MIC nel caso della
matrice (\ref{penta}) che, ricordiamo, assumiamo essere una
matrice di Stieltjes. Si premette il seguente risultato, riguardo
alla localizzazione degli autovalori di una matrice.

\begin{teo}{\bf (Gershg\"orin)} Se $A=(a_{ij})\in\RR^{n\times n}$,
allora $$\sigma(A)\subset \bigcup_{i=1}^n {\cal D}_i,$$ dove
$${\cal D}_i = \left\{z\in\CC: |z-a_{ii}|\le \sum_{j\ne
i}|a_{ij}|\right\},\qquad i=1,\dots,n.$$\end{teo}

\medskip
\begin{eser} Dimostrare il Teorema di Gershg\"orin.\end{eser}

Nel caso in esame, quindi, considerando che gli elementi
$\{r_i\}$ di $R$ (vedi (\ref{R1})-(\ref{R2})) sono non negativi,
un modo semplice per garantire la (\ref{neg}) \`e
fare in modo che la matrice $R$ abbia, sulla posizione
$i$-esima della dia\-gonale principale, $-(r_i+r_{i+N-1})$ (per convenzione,
$r_j=0$, se $j<N$ o $j>n$). In tal modo, dal Teorema di Gershg\"orin
segue che gli autovalori di $R$ dovranno essere non positivi. Questo
si ottiene semplicemente utilizzando la seguente definizione degli
elementi $\{\td_i\}$ della fattorizzazione, al posto di quella vista
in (\ref{ellij}):
\begin{equation}\label{mic}
\td_i =
a_i-\frac{b_i(b_i+c_{i+N-1})}{\td_{i-1}}-\frac{c_i(c_i+b_{i-N+1})}{\td_{i-N}},
\qquad i=1\dots,n,
\end{equation}

\no dove, al solito, $b_i=0$, per $i\le1$, e $c_i=0$, per $i\le N$ o $i>
n$.

\medskip
\begin{eser} Derivare la equazione di ricorrenza (\ref{mic}) per la
fattorizzazione modificata e verificare che essa soddisfa la
(\ref{neg}).\end{eser}

\begin{eser} Dimostrare che, nel caso della fattorizzazione
(\ref{mic}), l'autovalore pi\`u piccolo della matrice
precondizionata \`e $\lam=1$, cui corrisponde un autovettore con
componenti tutte uguali.\end{eser}

\begin{figure}[hp]
\begin{center}
\includegraphics[width=13cm,height=8.5cm]{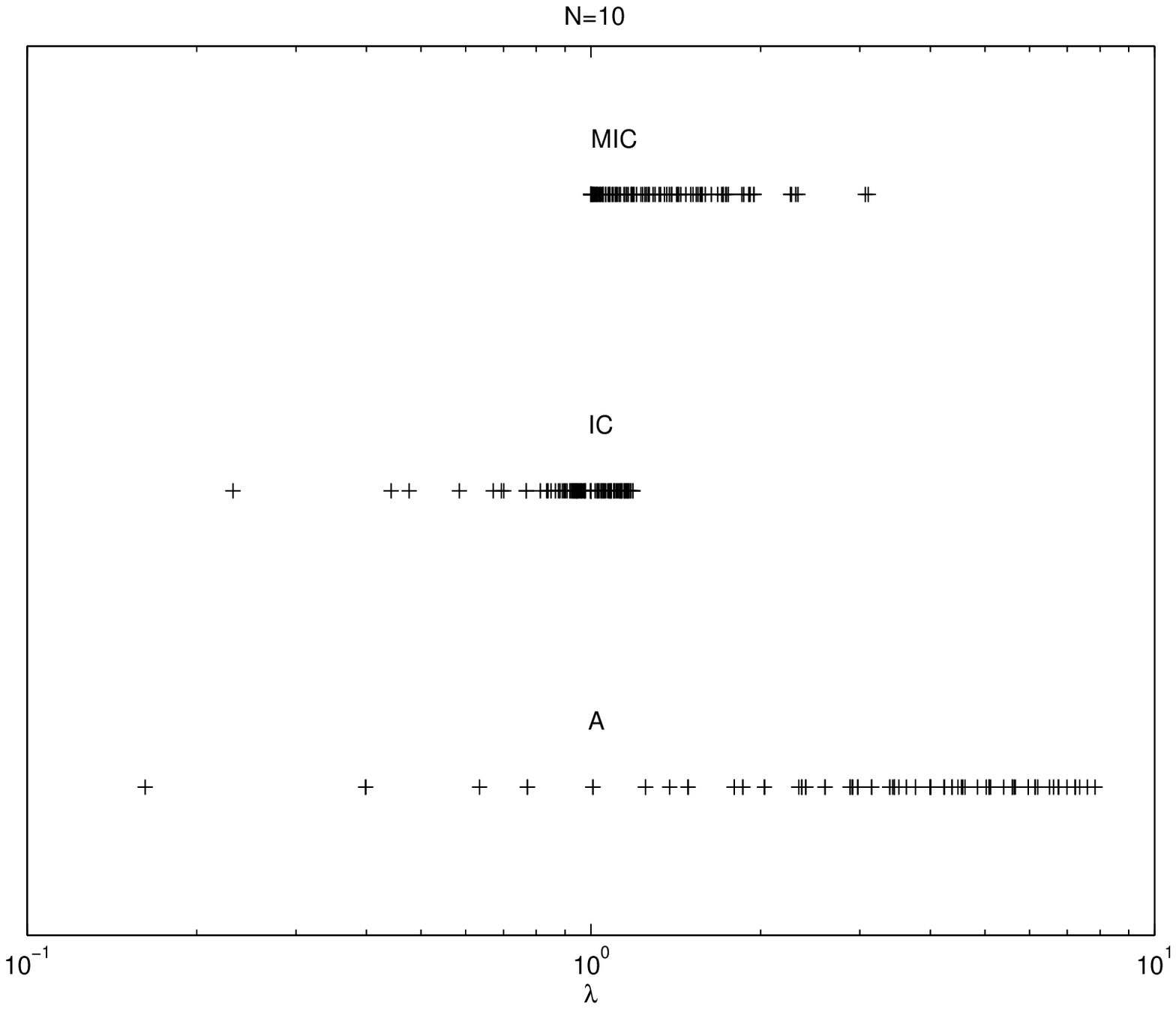}
\caption{\protect\label{fig6} Spettro della matrice $A$ in (\ref{test}), e
delle corrispondenti versioni precondizionate.}
\end{center}
\begin{center}
\includegraphics[width=13cm,height=8.5cm]{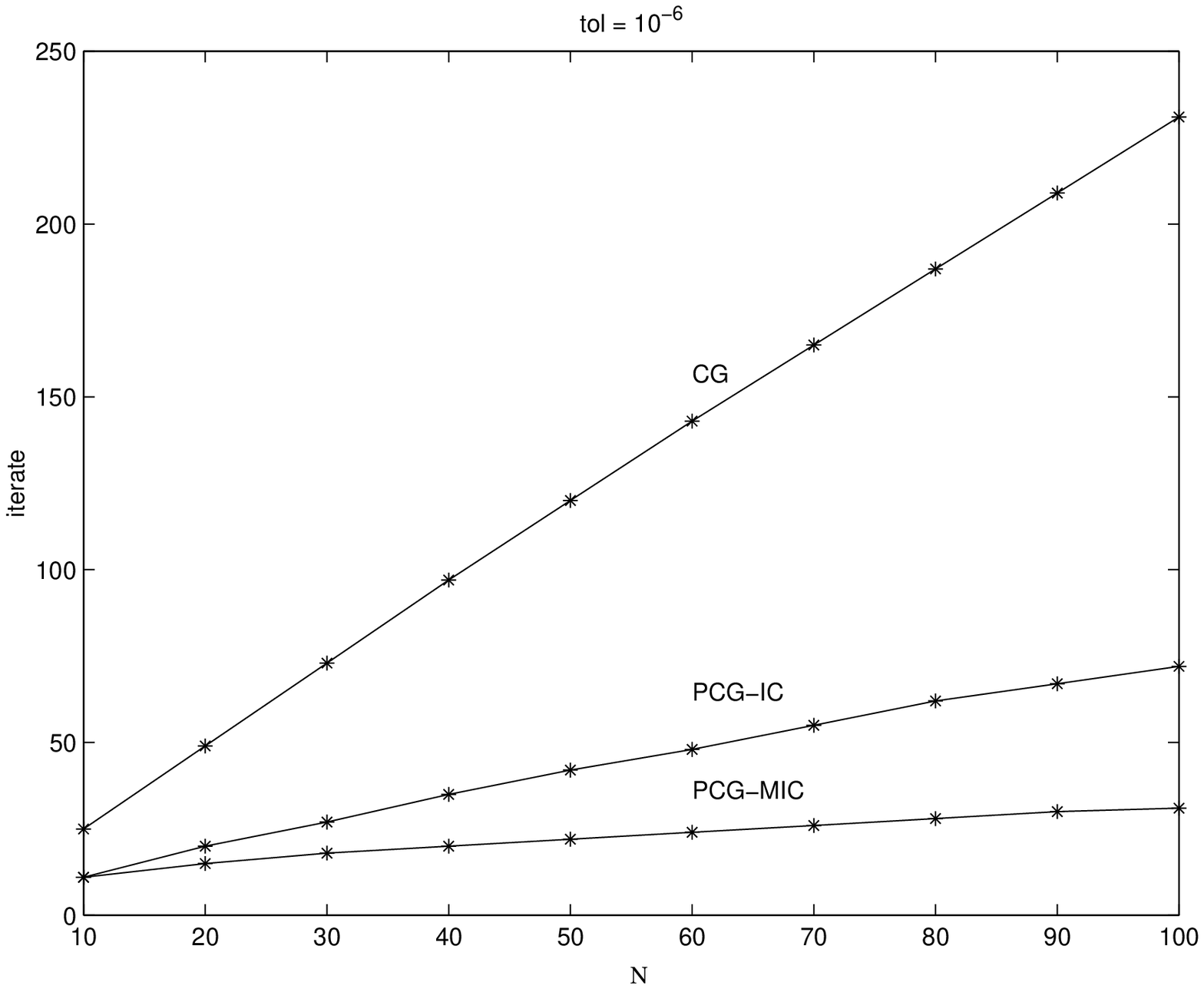}
\caption{\protect\label{fig7} Problema (\ref{test}), $N=10,\dots,100$;
confronto tra CG, IC e MIC.}
\end{center}
\end{figure}

A titolo di esempio, in Figura~\ref{fig6} riportiamo lo spettro
della matrice $A$ in (\ref{test}), $N=10$ (vedi
Esercizio~\ref{ex1}), assieme a quello della matrice
precondizionata mediante IC e mediante la versione modificata MIC.
Nell'ultimo caso, \`e evidente come il pi\`u piccolo autovalore della matrice
precondizionata sia $\lam=1$.

Concludiamo questa sezione comparando il metodo dei gradienti
coniugati non precondizionati (CG) con quello precondizionato con
IC (PCG-IC) e MIC (PCG-MIC), applicati sempre al problema
(\ref{test}) per valori di $N$ compresi tra 10 e 100 (le
dimensioni delle matrici corrispondenti vanno, quindi, da $10^2$ a
$10^4$). In Figura~\ref{fig7} si riportano le iterazioni richieste
dai vari metodi per soddisfare lo stesso criterio di arresto
basato sulla norma del residuo ($\|\bfr_i\|<10^{-6}$). In tutti i
casi, la soluzione iniziale considerata \`e stata il vettore
nullo. Come si evince dalla figura, le versioni precondizionate
risultano essere particolarmente vantaggiose, rispetto al metodo
senza precondizionatore. Tra le versioni precondizionate, MIC
risulta essere decisamente migliore (avendo, tra l'altro,
esattamente lo stesso costo per iterazione della IC).

\section{Generalizzazione a blocchi}\label{blockC}

La matrice del problema (\ref{test}), come anche la matrice
(\ref{A}), \`e una matrice di Stieltjes pentadiagonale. Tuttavia,
se la riguardiamo partizionata a blocchi (di dimensione $N\times
N$), essa risulta essere tridiagonale. Ovvero della forma
\begin{equation}\label{tridA}
A = \pmatrix{cccc} A_1 &B_2^\top \\ B_2 &\ddots &\ddots\\ &\ddots
&\ddots &B_N^\top \\ &&B_N &A_N\endpmatrix.\end{equation}

\no In particolare, nel caso delle due matrici su menzionate,
$B_i=-I_N, i=2,\dots,N$, ed i blocchi $A_i$ sono tridiagonali. Se
la matrice \`e sdp, allora sar\`a fattorizzabile nella forma
\begin{equation}\label{ldlb}A = L D^{-1}L^\top ,\end{equation}

\no dove \begin{equation} \label{ldlb1} L = \pmatrix{cccc}
D_1\\B_2&\ddots\\&\ddots &\ddots\\&&B_N &D_N\endpmatrix, \qquad D
= \pmatrix{ccc} D_1\\&\ddots\\&&D_N\endpmatrix,\end{equation}

\no e i blocchi diagonali sono definiti, ponendo per convenzione
$B_1=O$, da \begin{equation}\label{ldlb2} D_i =A_i -B_i
D_{i-1}^{-1} B_i^\top , \qquad i=1,\dots,N.\end{equation}

\medskip
\begin{eser} Verificare la fattorizzazione
(\ref{ldlb})--(\ref{ldlb2}).\end{eser}
\medskip

Osserviamo che il costo computazionale della fattorizzazione
(\ref{ldlb})--(\ref{ldlb2}) risulta essere di $O(N^4)$ operazioni
e $N^3$ posizioni di memoria. Risulta, pertanto, oneroso, quando
$N$ \`e grande, soprattutto quando i blocchi $A_i$ e $B_i$ sono
sparsi. \`E possibile, tuttavia, definire una versione
approssimata della (\ref{ldlb}), che definisce un {\em
precondizionatore a blocchi} per la matrice $A$, da utilizzarsi
con il metodo dei gradienti coniugati precondizionati. In maggior
dettaglio, in analogia con la (\ref{ldlr}), definiamo
\begin{equation}\label{blockM} A = \tL \tD^{-1} \tL -R \equiv
M-R,\end{equation}

\no dove $M$ \`e il precondizionatore, $R$ \`e la matrice di
errore, e $$\tL = \pmatrix{cccc} \tD_1\\B_2&\ddots\\&\ddots
&\ddots\\&&B_N &\tD_N\endpmatrix, \qquad \tD = \pmatrix{ccc}
\tD_1\\&\ddots\\&&\tD_N\endpmatrix.$$

\no I blocchi diagonali $\{\tD_i\}$ sono definiti da $$\tD_i =A_i
-B_i \Sigma_{i-1} B_i^\top , \qquad i=1,\dots,N,$$

\no essendo $\Sigma_{i-1}$ una {\em approssimazione sparsa} di
$\tD_{i-1}^{-1}$.

\medskip
\begin{eser} Verificare che, nella (\ref{blockM}),
$$R = \pmatrix{ccc} R_1\\ &\ddots \\&&R_N\endpmatrix,\qquad R_i =
B_i(\Sigma_{i-1}-\tD_{i-1}^{-1})B_i^\top .$$\end{eser}

\medskip
\begin{oss} In modo analogo a quanto visto nella
precedente sezione per la IC, \`e possibile definire una
corrispondente versione modificata del precondizionatore $M$ in
(\ref{blockM}). Tuttavia, questo pu\`o non essere sempre
conveniente, a causa del costo per ottenerla.\end{oss}

\section{Precondizionatori polinomiali}\label{polyC}

Fino ad ora abbiamo derivato, in diversi modi, la matrice di
precondizio\-namento $M\equiv C^{-1} \approx A$. Pertanto $M$
risulta essere una approssimazione di $A$, la cui inversa sia
facilmente (ed efficientemente) calcolabile. In que\-sta sezione
cercheremo di derivare, per la prima volta, una espressione per
$C\approx A^{-1}$. Questa approssimazione, vedremo, commuta con
$A$, cosicch\'e se $C$ \`e sdp, allora tale \`e anche $CA$.
Pertanto, la scelta $B=CA$ per i gradienti coniugati risulta
essere la pi\`u naturale.

\begin{oss} Abbiamo visto che il metodo dei gradienti coniugati
mi\-nimizza, ad ogni passo, la $B$-norma dell'errore in un opportuno
sottospazio. Quindi, se si utilizza la norma $B=CA$, con $C\approx
A^{-1}$, segue che $B\approx I$, ovvero la norma dell'errore che
viene minimizzata \`e approssimativamente quella
euclidea.\end{oss}

Ricercheremo il precondizionatore $C$ nella forma $$C = C(A),
\qquad C(z)\in\Pi_{m-1},$$

\no per un certo $m\in\NN$ assegnato. Osserviamo che $C(A)$ \`e
simmetrica e, inol\-tre, commuta con $A$. Idealmente, si
desidererebbe $C(A)\equiv A^{-1}$, ma questo implicherebbe che
$C(\lam)$ dovrebbe essere il polinomio interpolante la funzione
$f(\lam)=\lam^{-1}$ sui punti dello spettro di $A$, che indichiamo con
$$\sigma(A)=\{0<\lam_{min}\equiv\lam_1\le\lam_2\le\dots\le\lam_n\equiv\lam_{max}\}.$$

\no Questo implica che il relativo grado \`e potenzialmente
elevato, se $A$ ha molti autovalori distinti e, peraltro, $\sigma(A)$
\`e generalmente non noto. In realt\`a, noi
ricerchiamo solo una approssimazione dell'inversa, per cui,
fissato $m$, cercheremo il polinomio di grado al pi\`u $m-1$ che
minimizza la differenza $$|\lam^{-1}-C(\lam)|,\qquad
\lam\in[\lam_{min},\lam_{max}].$$

\no Se invece di minimizzare l'errore assoluto si minimizza
quello relativo (il che \`e pi\`u pertinente, nel nostro caso),
si ottiene, infine, il seguente problema di minimassimo vincolato:
\begin{eqnarray}\nonumber
\lefteqn{
\min_{C\in\Pi_{m-1}}\,\max_{\lam\in[\lam_{min},\lam_{max}]}\,|1-C(\lam)\lam|}\\[1mm]
\nonumber &\qquad&\equiv \min_{p\in\Pi_m, p(0)=0}\,
\max_{\lam\in[\lam_{min},\lam_{max}]}\,|1-p(\lam)|\\[1mm] &&\equiv
\min_{e\in\Pi_m, e(0)=1}\,
\max_{\lam\in[\lam_{min},\lam_{max}]}\,|e(\lam)|.\label{minmaxe}
\end{eqnarray}

\no Poich\`e il punto di interpolazione, $\lam=0$, \`e esterno
all'intervallo $[\lam_{min},\lam_{max}]$, la soluzione del
problema (\ref{minmaxe}) \`e noto essere data da (vedi
Sezione~\ref{cebypol})
\begin{equation}\label{em}
e_m(\lam) = \frac{T_m(\mu(\lam))}{T_m(\mu(0))}, \qquad \mu(\lam) =
\frac{\lam_{max}+\lam_{min}-2\lam}{\lam_{max}-\lam_{min}},
\end{equation}

\no dove $T_m(x)$ \`e l'$m$-esimo polinomio di Chebyshev di prima
specie. Ne consegue che la corrispondente matrice
precondizionata sar\`a
\begin{equation}\label{pm}
C(A)A \equiv C_{m-1}(A)A = p_m(A) = I - e_m(A)\equiv
I-\frac{T_m(\mu(A))}{T_m(\mu(0))},
\end{equation}

\no dove, evidentemente,
\begin{equation}\label{muA}
\mu(A) =
\frac{(\lam_{max}+\lam_{min})I-2A}{\lam_{max}-\lam_{min}}.
\end{equation}

\no Osserviamo che gli autovalori di $p_m(A)$ sono dati da
$$1-e_m(\lam), \qquad \lam\in\sigma(A),$$

\no per cui, se
\begin{equation}\label{epsm}
\eps_m \equiv \max_{\lam\in[\lam_{min},\lam_{max}]} |e_m(\lam)|
<1,\end{equation}

\no allora $p_m(A)$ \`e sdp ed \`e possibile utilizzare la norma
indotta da $B=p_m(A)$. Inoltre, in tal caso si verifica facilmente
che \begin{equation}\label{kpmA}\kappa(p_m(A)) \le
\frac{1+\eps_m}{1-\eps_m}.\end{equation}

\no Osserviamo che l'ottimalit\`a di (\ref{em}) implica che la
maggiorazione di $\kappa(p_m(A))$ data dalla (\ref{kpmA}) \`e
minimizzata. Valgono i seguenti risultati.

\begin{teo}\label{emmin1} Per ogni ~$m\ge1$~ si ha ~$\eps_m<1$.\end{teo}

\proof Infatti, dalle (\ref{em})--(\ref{epsm}) segue che $$\eps_m
= T_m(\mu(0))^{-1}, \qquad
\mu(0)=\frac{\lam_{max}+\lam_{min}}{\lam_{max}-\lam_{min}}>1.$$

\no Infine (vedi Sezione~\ref{cebypol}) la successione
$\{T_m(\mu(0))\}$ \`e strettamente monotona crescente e
divergente, avendosi $$T_m(\mu(0)) = \cosh\left(
m\log\left(\mu(0)+\sqrt{\mu(0)^2-1}\right)\right),$$

\no e, inoltre, $T_1(\mu(0))=\mu(0)>1.$\QED
\bigskip

Pertanto, la successione $\{\eps_m\}$ \`e infinitesima. Da questo
discende immediatamente, considerando la (\ref{kpmA}), il seguente corollario.

\begin{cor} $\kappa(p_m(A)) \rightarrow 1$, per
$m\rightarrow\infty$.\end{cor}

\begin{figure}[hp]
\begin{center}
\includegraphics[width=13cm,height=8.5cm]{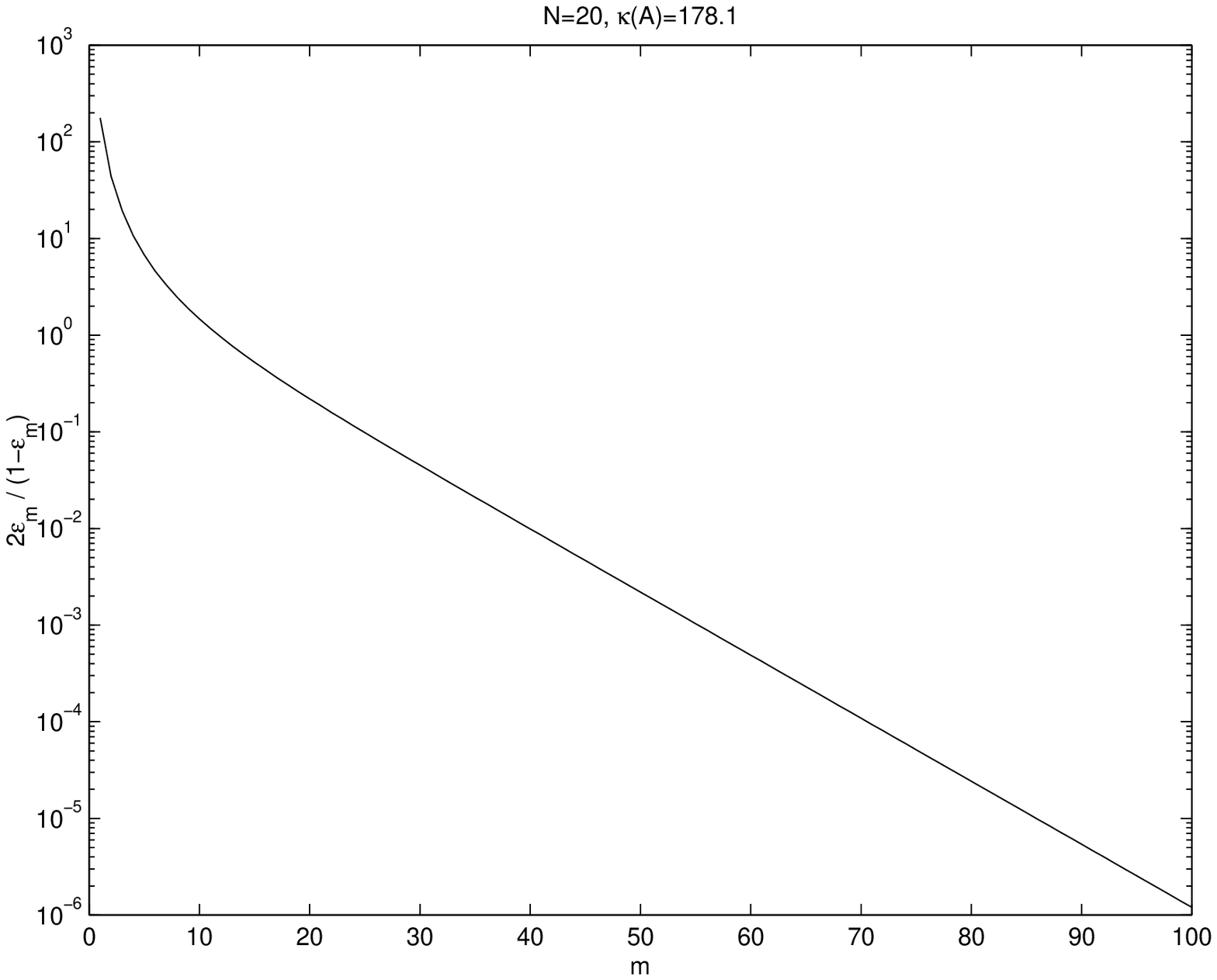}
\caption{\protect\label{fig8} Maggiorazione di $\kappa(p_m(A))-1$
per la matrice in (\ref{test}), $N=20$.}
\end{center}
\begin{center}
\includegraphics[width=13cm,height=8.5cm]{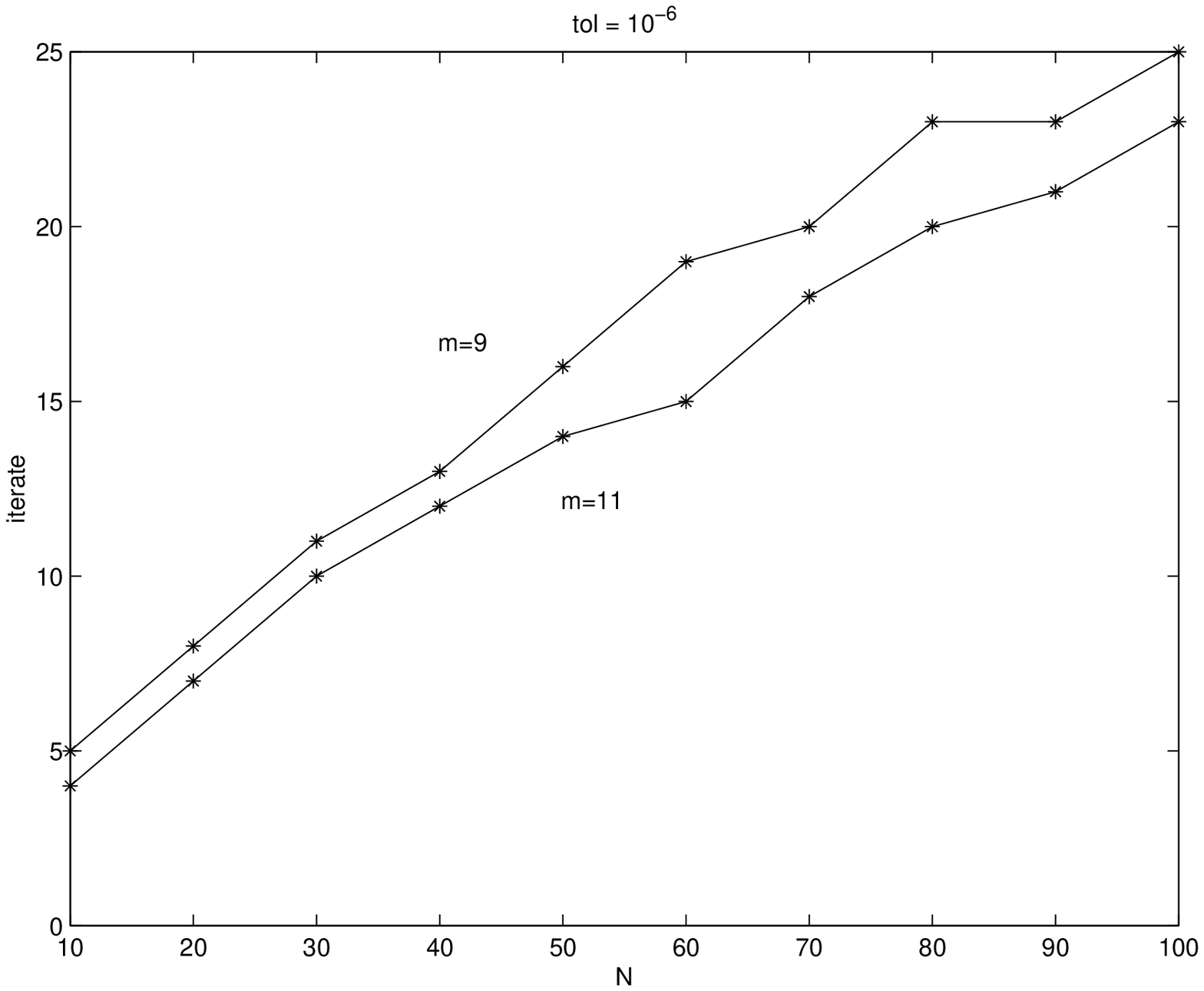}
\caption{\protect\label{fig9} Problema (\ref{test}),
$N=10,\dots,100$; confronto tra i precondiziona\-tori polinomiali di
grado $m=9$ e $m=11$.}
\end{center}
\end{figure}

In Figura~\ref{fig8} riportiamo, a titolo di esempio, la
quantit\`a (vedi (\ref{kpmA})) $$\frac{2\eps_m}{1-\eps_m}\ge \kappa(p_m(A))-1,$$
rispetto ad $m$, per la matrice $A$ nell'equazione (\ref{test}), $N=20$.

Ancora dal Teorema~\ref{emmin1} discende il seguente risultato
che, essenzialmente, ci dice che la matrice precondizionata ha un
numero di condizionamento minore di quello di $A$, per quasi ogni
valore di $m$.

\begin{cor} $\kappa(p_m(A))<\kappa(A)$, per ogni $m\ge2$.\end{cor}

\proof Dalla (\ref{kpmA}) segue infatti che la tesi \`e vera se
$$\frac{1+\eps_m}{1-\eps_m}<\kappa(A) \equiv
\frac{\lam_{max}}{\lam_{min}},$$

\no ovvero quando $$\eps_m <
\frac{\lam_{max}-\lam_{min}}{\lam_{max}+\lam_{min}}\equiv
\mu(0)^{-1} \equiv \eps_1.~\QED$$

\medskip
\begin{oss} Evidentemente, $\kappa(p_1(A))\le \kappa(A)$.\end{oss}

In Figura~\ref{fig9} riportiamo, come esempio, il numero di
iterazioni per soddisfare lo stesso criterio di arresto visto per
i metodi in Figura~\ref{fig7}, usando i precondizionatori di grado
$m=9$ e $m=11$.

\medskip
\begin{oss} Osserviamo che, per confrontare correttamente i
risultati riportati in Figura~\ref{fig9} con quelli di
Figura~\ref{fig7}, bisogna considerare che ogni iterazione con il
precondizionatore polinomiale ha un costo molto superiore a quello
dei precondizionatori IC e MIC. Inoltre, il costo aumenta al crescere del grado
$m$. Questo aspetto sar\`a approfondito nella prossima
sezione.\end{oss}

\section{Dettagli implementativi}

Abbiamo detto che il precondizionatore polinomiale \`e
convenientemente implementato nel metodo precondizionato dei
gradienti coniugati con la scelta $B=C_{m-1}(A)A\equiv p_m(A)$. Questo
implica che \`e formalmente possibile utilizzare
l'Algoritmo~\ref{algo4} con le seguenti sostituzioni formali (vedi
(\ref{pm})):

\begin{itemize}

\item $A\leftarrow p_m(A);$

\item $\bfb\leftarrow C_{m-1}(A)\bfb$.

\end{itemize}

\no Ne consegue che abbiamo bisogno di algoritmi efficienti per
calcolare la matrice $p_m(A)$ e, inizialmente, il vettore
$C_{m-1}(A)\bfb$.

Preliminarmente, osserviamo che abbiamo bisogno anche di conoscere
gli autovalori estremi di $A, \lam_{min}$ e $\lam_{max}$, per
poter definire la (\ref{muA}). Abbia\-mo gi\`a visto in
Sezione~\ref{stopcg} come ottenere buone approssimazioni di queste
quantit\`a, quando non fossero note, eseguendo alcune iterazioni
del metodo dei gradienti coniugati che, tra l'altro, fornirebbero
anche un punto iniziale pi\`u favorevole. Va tuttavia sottolineato
che le propriet\`a di convergenza esa\-minate nella precedente
sezione dipendono in modo cruciale da tali stime, per cui se
queste non fossero sufficientemente accurate, potrebbe aversi una
degradazione significativa della ``performance'' del precondizionatore.

Dopo questa premessa, passiamo ad esaminare i due punti su
evidenziati.
\medskip

Per il calcolo della matrice precondizionata, osserviamo che
nell'Algoritmo~\ref{algo4} \`e solo richiesto di poter
effettuare il prodotto (vedi (\ref{pm}))
\begin{equation}\label{bfv}\bfv = p_m(A)\bfp \equiv \bfp
-T_m(\mu(0))^{-1}\left(\,T_m(\mu(A))\bfp\,\right).\end{equation}

\no Pertanto, il problema \`e risolto ricordando che $T_m(\mu(0))$
si calcola ricorsivamente come
\begin{eqnarray*}
T_0(\mu(0))&=&1,\\[1mm] T_1(\mu(0))&=&\mu(0),\\[1mm]
T_{i+1}(\mu(0))&=&2\mu(0)T_i(\mu(0))-T_{i-1}(\mu(0)), \qquad
i=1,\dots,m-1,\end{eqnarray*}

\no ed individuando una procedura efficiente per calcolare
$T_m(\mu(A))\bfp$. A questo fine, osserviamo che risulta essere
conveniente assemblare la matrice (vedi (\ref{muA}))
\begin{equation}\label{G} G \equiv
2\mu(A),\end{equation} che richiede un costo lineare in $n$, se $A$ \`e
sparsa. Osserviamo che, fisicamente, $G$ pu\`o essere riscritta su
$A$: pertanto non \`e necessario un costo addizionale, in termini
di memoria. Quindi, si verifica facilmente che $\bfy_m\equiv
T_m(\mu(A))\bfp$, dove
\begin{eqnarray*}
\bfy_0 &=& \bfp,\\[1mm] \bfy_1 &=& \frac{1}2( G\bfp ),\\[1mm] \bfy_{i+1} &=&
G\bfy_i -\bfy_{i-1}, \qquad i=1,\dots,m-1.
\end{eqnarray*}

\no Se ne conclude che il calcolo di $\bfv$ richiede (vedi
(\ref{bfv})) $m$ matvec ed $m$ axpy.

\begin{eser} Scrivere un codice che implementi il calcolo del
vettore $\bfv$ in (\ref{bfv}).\end{eser}

Esaminiamo, ora, il calcolo del vettore $C_{m-1}(A)\bfb$,
richiesto per il calcolo del solo residuo iniziale. Teoricamente,
essendo \begin{equation}\label{Cm1}C_{m-1}(\lam) =
\frac{p_m(\lam)}\lam = \frac{1-e_m(\lam)}\lam \equiv
\frac{T_m(\mu(0))-T_m(\mu(\lam))}{\lam\,T_m(\mu(0))},
\end{equation}

\no potremmo utilizzare la base delle potenze per rappresentare
$p_m(\lam)$ ed effettuare la divisione formale per $\lam$, in
quanto, ricordiamo, $p_m(0)=0$. Tuttavia, questa procedura non \`e
raccomandabile, in quanto, anche per modesti valori di $m$, i
coefficienti di $p_m(\lam)$, in questa base, risultano essere
assai grandi e di segno discorde, rendendo inaffidabile il calcolo
in aritmetica finita. Ad esempio, se \begin{equation}\label{exm11}
m=11, \qquad \lam_{min}=10^{-3}, \qquad \lam_{max}=1.001,\end{equation}

\no si ottiene che $\mu(\lam) = 1.002 -2\lam$ e $p_{11}(\lam) =
\sum_{i=0}^{11} c_i \lam^{11-i},$ con i coefficienti della
rappresentazione dati, approssimativamente, da:
\begin{equation}\label{ci}\pmatrix{c} c_0\\ \vdots\\
c_{11}\endpmatrix =
\pmatrix{r}
  1.6752660329527138e+006\\
   -9.2323911076024063e+006\\
   2.1975394326889034e+007\\
   -2.9566871007963315e+007\\
   2.4710334661370583e+007\\
   -1.3273250212561980e+007\\
   4.5826267022750815e+006\\
   -9.8775156863499968e+005\\
   1.2453502849456538e+005\\
   -8.1003898374282262e+003\\
   2.0914774360373008e+002\\
   0.0000000000000000e+000\endpmatrix.
 \end{equation}

\no Per ottenere una rappresentazione computazionalmente pi\`u
efficiente non conviene, quindi, ricorrere alla usuale base delle potenze.
\`E invece preferibile utilizzare, come polinomi di base, dei polinomi
che soddisfino alla stessa relazione di ricorrenza (\ref{ceby}) dei
polinomi di Chebyshev. Dalla (\ref{Cm1}) segue che
\begin{equation}\label{Sm}
C_{m-1}(\lam)
= \frac{S_m(\lam)}{\lam\,T_m(\mu(0))}, \qquad S_m(\lam) =
T_m(\mu(0))-T_m(\mu(\lam)).\end{equation}

\no Si verifica che la successione $\{S_i(\lam)\}$ soddisfa alla
equazione alle differenze
\begin{eqnarray}\nonumber
S_0(\lam) &\equiv& 0,\\ \label{Si}
S_1(\lam) &=& \frac{2\lam}{\lam_{max}-\lam_{min}},
\\ \nonumber
S_{i+1}(\lam)&=&2\mu(\lam)S_i(\lam) -S_{i-1}(\lam) +
\frac{4\lam}{\lam_{max}-\lam_{min}}T_i(\mu(0)), \quad i=1,2,\dots,
\end{eqnarray}

\no la cui soluzione si verifica essere
\begin{equation}\label{solSm}
S_m(\lam) = \frac{2\lam}{\lam_{max}-\lam_{min}}\left(
2\sum_{k=0}^{m-1}
T_k(\mu(0))U_{m-1-k}(\mu(\lam))-U_{m-1}(\mu(\lam))\right),
\end{equation}
\vspace{2mm}

\no in cui gli $\{U_k\}$ sono i polinomi di Chebyshev di seconda
specie. Dalla (\ref{Sm}) segue quindi che
\begin{equation}\label{CmUk}
C_{m-1}(\lam) = \sum_{k=0}^{m-1} \cc_k U_{m-1-k}(\mu(\lam)),
\end{equation}

\no dove
\begin{eqnarray*}
\cc_0 &=& \frac{2}{(\lam_{max}-\lam_{min})T_m(\mu(0))}, \\[2mm]
\cc_k &=& \frac{4T_k(\mu(0))}{(\lam_{max}-\lam_{min})T_m(\mu(0))},
\qquad k=1,\dots,m-1.
\end{eqnarray*}

\no Considerando che la successione $\{T_k(\mu(0))\}$ \`e
crescente e ad elementi positivi, se ne conclude che tutti i
coefficienti appartengono all'intervallo
$$\left(0,4(\lam_{max}-\lam_{min})^{-1}\right).$$

\no Nel caso dell'esempio (\ref{exm11}) si ottiene l'intervallo
$(0,4)$, da confrontare con i coefficienti riportati nella (\ref{ci}).

\medskip
\begin{eser} Verificare che la successione $\{S_i(\lam)\}$,
definita dalla (\ref{Sm}), soddisfa alla equazione (\ref{Si}).
\end{eser}

\medskip
\begin{eser} Verificare che la (\ref{solSm}) \`e la soluzione del
problema (\ref{Si}).\end{eser}

Vediamo adesso di utilizzare la (\ref{CmUk}) per calcolare in modo
efficiente $C_{m-1}$ in un punto $\lam$ assegnato. A questo fine,
utilizziamo il seguente {\em algoritmo di Clenshaw}:
\begin{eqnarray}\nonumber
y_{-1} &=&0,\\ \label{clen}
y_0    &=&\cc_0,\\ \nonumber
y_k    &=& 2\mu(\lam)y_{k-1} -y_{k-2} +\cc_k, \qquad
k=1,\dots,m-1.
\end{eqnarray}

\no Per la successione $\{y_k\}$ generata dalla (\ref{clen}) vale il
seguente risultato.

\begin{teo} $y_{m-1} = C_{m-1}(\lam).$ \end{teo}

\proof Definiamo il vettore $$\bfcc = (\cc_{m-1},\dots,\cc_0)^\top $$
contenente i coefficienti della rappresentazione (\ref{CmUk}) di
$C_{m-1}$ e, fissato $\lam$, il vettore $$\bfu =
(U_0(\mu(\lam)),\dots,U_{m-1}(\mu(\lam)))^\top .$$ Evidentemente,
dalla (\ref{CmUk}) segue che
$$C_{m-1}(\lam)= \bfcc^\top \bfu.$$ Si verifica, inoltre, che il
vettore $\bfu$ \`e soluzione del sistema lineare $$A_\lam \bfu =
E_1,$$ dove
$$A_\lam = \pmatrix{ccccc}
1\\
-2\mu(\lam) &1\\
1 &-2\mu(\lam) &1\\
  &\ddots &\ddots &\ddots\\
  &       &1 &-2\mu(\lam) &1\endpmatrix_{m\times m}.$$
Introducendo la matrice antidiagonale $$J=\pmatrix{cccc}
&&&1\\&&\cdot\\&\cdot\\1\endpmatrix_{m\times m},$$ si verificano
facilmente le seguenti propriet\`a:
\begin{itemize}
\item $J=J^\top , J^2 = I$;
\item $JA_\lam J = A_\lam^\top $,
\item $JE_1 = E_m$;
\item $\hat{\bfcc} \equiv J\bfcc = (\cc_0,\dots,\cc_{m-1})^\top $.
\end{itemize}
Pertanto, segue che $$C_{m-1}(\lam) = \bfcc^\top  A_\lam^{-1}E_1 =
\bfcc^\top  J^2A_\lam^{-1}J^2 E_1 = \hat{\bfcc}^\top A_\lam^{-T}E_m
= E_m^\top  A_\lam^{-1}\hat{\bfcc}.$$
La tesi si completa osservando che, dalla (\ref{clen}), il vettore $$\bfy =
(y_0,\dots,y_{m-1})^\top ,$$ risulta essere soluzione del sistema lineare
$$A_\lam\bfy = \hat{\bfcc},$$ da cui si ottiene $$y_{m-1} =
E_m^\top A_\lam^{-1}\hat{\bfcc}.~\QED$$

Questo risultato pu\`o essere utilmente impiegato per il calcolo
del vettore $C_{m-1}(A)\bfb$. Infatti, ricordando la definizione
della matrice $G$ in (\ref{G}), dalla (\ref{clen}) si ottiene che la
successione
\begin{eqnarray}\nonumber
\bfy_{-1} &=&\bfo,\\ \label{clen1}
\bfy_0    &=&\cc_0\bfb,\\ \nonumber
\bfy_k    &=& G\bfy_{k-1} -\bfy_{k-2} +\cc_k\bfb, \qquad
k=1,\dots,m-1,
\end{eqnarray}

\no \`e tale che $\bfy_{m-1} = C_{m-1}(A)\bfb$, con un costo,
quindi, di $m-1$ matvec e $2m-2$ axpy.

\medskip
\begin{eser} Scrivere un codice che implementi efficientemente l'algoritmo
(\ref{clen1}) per il calcolo di $C_{m-1}(A)\bfb$.\end{eser}

\medskip
\begin{eser} Si supponga che gli autovalori estremi di una matrice
sdp $A$ siano $\lam_{min}=0.5$, ~e~ $\lam_{max} = 10000.5$.
Quale \`e il grado minimo del polinomio $p(z)$ tale che:
\begin{itemize}
\item $p(0)=0$,
\item $p(A)$ \`e la migliore approssimazione della matrice identit\`a,
\item $\kappa(p(A))\le100$?
\end{itemize}
(Considerare che $(\log_2(1+10^{-2}\sqrt{2}))^{-1} \approx 49.4$).
\end{eser}

%
%
\chapter{Il caso $A$ simmetrica}\label{cap6}

{\em In questo capitolo trattiamo il caso della risoluzione di sistemi lineari
la cui matrice dei coefficienti sia una matrice simmetrica nonsingolare.}

\section{Introduzione}
Nel caso in cui la matrice del sistema lineare
\begin{equation}\label{sistsim}
A\bfx = \bfb
\end{equation}

\no non sia pi\`u sdp, allora i metodi analizzati precedentemente
possono non essere pi\`u ben definiti, in quanto, essenzialmente,
la matrice $A$ non induce pi\`u un prodotto scalare ed una norma
corrispondente. Si deve, quindi, modificare l'approccio per la
risoluzione del problema. Tuttavia, laddove queste esistono, tutte
le propriet\`a della matrice $A$ devono essere convenientemente
utilizzate, al fine di ottenere metodi di risoluzione efficienti.
Per questo motivo, in questo capitolo e nel successivo,
``impoveriremo'' progressivamente la matrice $A$. Il motivo di
ci\`o risiede nel fatto che il caso in cui $A$ sia, oltre che
nonsingolare, anche simmetrica, merita una analisi separata dal
caso in cui $A$ sia una generica matrice nonsingolare. \`E questo
il caso che andiamo a trattare in questo capitolo, dove assumeremo
costantemente la propriet\`a di simmetria della matrice $A$ del
sistema lineare (\ref{sistsim}).

Concludiamo questa introduzione osservando che la via di risolvere il problema
(\ref{sistsim}) mediante la risoluzione del sistema equivalente ({\em
equazioni normali})
\begin{equation}\label{normaleq}
A^\top A\bfx = A^\top \bfb,
\end{equation}

\no la cui matrice \`e sdp, se $A$ \`e nonsingolare, non \`e una prassi
raccomandabile. Questo \`e dovuto al fatto che, generalmente, $\kappa(A^\top A)\gg
\kappa(A)$. Ad esempio, se si utilizza la norma 2, $$\kappa(A^\top A) =
\kappa(A)^2.$$

\section{Il metodo di Lanczos}\label{lanc}

Nel caso in cui $A$ sia una matrice simmetrica, se ne ricerca la
seguente fattorizzazione, alternativa alla (\ref{fattA}) esaminata
nel caso sdp,
\begin{equation}\label{fattlanc}
AU=UT \equiv (\bfu_1 \dots \bfu_n)\pmatrix{cccc} \cc_1 &\bb_1\\
\bb_1 &\ddots &\ddots\\
      &\ddots &\ddots &\bb_{n-1}\\
      &       &\bb_{n-1}   &\cc_n\endpmatrix,
\qquad U^\top U=I.
\end{equation}

\no Questa fattorizzazione, se esiste, definisce la riduzione a
forma tridiagonale della matrice $A$ mediante una trasformazione
ortogonale. La corrispondente procedura di tridiagonalizzazione
\`e nota come {\em metodo di Lanczos}.

\begin{oss} Osserviamo che, rispetto alla (\ref{fattA}), la
matrice $B$, che definisce il prodotto scalare, \`e scelta uguale
alla matrice identit\`a $I$. \end{oss}

Vediamo di derivare l'algoritmo corrispondente alla
(\ref{fattlanc}). Osserviamo che, essendo $U$ una matrice
ortogonale, $$U^\top AU = U^\top UT = T.$$

\no Questa equazione:

\begin{itemize}

\item conferma che la matrice $T$ deve essere simmetrica;

\item ci dice che le matrici $A$ e $T$ sono simili;

\item permette di derivare le seguenti espressioni per i
coefficienti di $T$:

\begin{equation}\label{T}
\bfu_i^\top A\bfu_i = \cc_i, \qquad \bfu_i^\top A\bfu_{i+1} = \bb_i,
\qquad i=1,\dots,n,
\end{equation}

\no dove, come nel caso della (\ref{fattA}), assumiamo
$\bfu_{n+1}=\bfo$.

\end{itemize}

\no Inoltre, considerando le colonne $i$-esime dei due membri
della prima equazione in (\ref{fattlanc}) e ponendo per
convenzione $\bb_0=0$, si ottiene la seguente equazione di
ricorrenza,

\begin{equation}\label{Auilanc}
A\bfu_i = \bb_i\bfu_{i+1} +\cc_i\bfu_i +\bb_{i-1}\bfu_{i-1},
\qquad i=1,\dots,n.
\end{equation}

\begin{algo}\label{algo6} Fattorizzazione (\ref{fattlanc}).\\ \rm
\fbox{\parbox{12cm}{
\begin{itemize}

\nulit sia assegnato un vettore iniziale $\bfu_1$, $\|\bfu_1\|=1$

\nulit inizializzo $\bb_0=0$

\nulit per i=1,\dots,n:

\begin{itemize}

\nulit $\bfv_i=A\bfu_i$

\nulit $\cc_i = \bfu_i^\top \bfv_i$

\nulit $\bfv_i \leftarrow \bfv_i -\cc_i\bfu_i-\bb_{i-1}\bfu_{i-1}$

\nulit $\bb_i = \|\bfv_i\|$

\nulit se $\bb_i=0$, esci, fine se

\nulit $\bfu_{i+1} = \bfv_i/\bb_i$

\end{itemize}

\nulit fine per

\end{itemize}
}}\end{algo}

\no Questa relazione, tenendo conto della simmetria della matrice
$T$, induce il seguente modo alternativo per calcolare $\bb_i$ ed
il nuovo vettore $\bfu_{i+1}$:

\begin{eqnarray}\label{vilanc}
\bfv_i &=& A\bfu_i -\cc_i\bfu_i -\bb_{i-1}\bfu_{i-1},\\ \bb_i &=&
\|\bfv_i\|, \qquad \bfu_{i+1} = \bfv_i/\bb_i, \qquad
i=1,\dots,n-1.\nonumber
\end{eqnarray}

\no Possiamo riassumere le precedenti espressioni
(\ref{T})--(\ref{vilanc}) nell'\,Algoritmo~\ref{algo6}.

\begin{oss}\label{invar} Osserviamo che, qualora l'Algoritmo~\ref{algo6}
termini ad un passo $i<n$ (cio\'e si verifichi $\bb_i=0$), allora,
dette
\begin{equation} \label{invar1}U_i = (\bfu_1\dots\bfu_i), \qquad
T_i = \pmatrix{cccc} \cc_1 &\bb_1\\ \bb_1 &\ddots &\ddots\\
      &\ddots &\ddots &\bb_{i-1}\\
      &       &\bb_{i-1}   &\cc_i\endpmatrix,\end{equation}

\no avremo determinato \begin{equation}\label{invar2}
AU_i = U_iT_i,\qquad U_i^\top U_i=I_i,\end{equation}
ovvero una base ortonormale per un sottospazio $i$-dimensionale
invariante di $A$. Ne consegue che gli autovalori di $T_i$ sono
$i$ degli autovalori di $A$.

Viceversa, se $\bb_i\ne0$, allora dalla (\ref{fattlanc}) si
ottiene che
\begin{equation}\label{AUi}AU_i =
U_{i+1}\hat{T}_i \equiv U_{i+1}\pmatrix{cccc} \cc_1 &\bb_1\\ \bb_1
&\ddots &\ddots\\
      &\ddots &\ddots &\bb_{i-1}\\
      &       &\bb_{i-1} &\cc_i\\ \hline
      &       &          &\bb_i\endpmatrix
      \equiv U_{i+1} \pmatrix{c} T_i\\ \hline
      \bb_i(E_i^{(i)})^\top \endpmatrix.
      \end{equation}

\end{oss}

Per verificare che il precedente Algoritmo~\ref{algo6} determina
proprio la fattorizzazione (\ref{fattlanc}), rimane da dimostrare
il seguente risultato.

\begin{teo} I vettori $\{\bfu_i\}$ sono ortogonali.\end{teo}

\proof Dimostriamo, per induzione su $k\le n$, che
$\{\bfu_1\dots,\bfu_k\}$ sono orto\-gonali. Per $k=2$ la tesi
segue facilmente dalle (\ref{T}) e (\ref{Auilanc}) per
$i=1$. Supponiamo ora vera la tesi per $k$ e dimostriamo che
$\bfu_j^\top \bfu_{k+1}=0, j\le k$, assumendo che $\bfu_{k+1}\ne\bfo$.
Per $j\le k-2$ la tesi segue dall'ipotesi di induzione osservando
che $\bb_k\ne0$ e, dalla (\ref{Auilanc}),
\begin{eqnarray*}
\bb_k\bfu_j^\top \bfu_{k+1}&=&\bfu_j^\top (A\bfu_k
-\cc_k\bfu_k-\bb_{k-1}\bfu_{k-1}) = \bfu_j^\top A\bfu_k\\
&=&(\bb_j\bfu_{j+1}+\cc_j\bfu_j+\bb_{j-1}\bfu_{j-1})^\top \bfu_k
= 0.
\end{eqnarray*}

\no Per $j=k-1$ e $j=k$, la tesi discende da simili argomenti, in
virt\`u della (\ref{T}), con $i=k-1$ e $i=k$, rispettivamente.\QED

\begin{eser} Scrivere un codice che implementi efficientemente
l'Algoritmo~\ref{algo6}.\end{eser}

\section{MINRES}\label{minres}

Vediamo come la fattorizzazione (\ref{fattlanc}) possa essere
utilizzata per definire un metodo iterativo di risoluzione per
(\ref{sistsim}). In questo caso, data una approssimazione iniziale
$\bfx_0$ della soluzione, scegliamo il vettore $\bfu_1$ come
\begin{equation}\label{u1lanc}
\bfu_1 = \frac{\bfr_0}{\|\bfr_0\|}, \qquad \bfr_0 = \bfb-A\bfx_0.
\end{equation}

\no Osserviamo che, a meno di un fattore di scalamento, si tratta
della stessa scelta fatta nel caso del metodo dei gradienti
coniugati (vedi (\ref{u1})), nel caso $A$ sdp. Se procedessimo con
questa analogia, potremmo ricercare
\begin{equation}\label{xilanc} \bfx_i\in S_i =
\bfx_0+[\bfu_1,\dots,\bfu_i],\end{equation}

\no imponendo la seguente condizione di ortogonalit\`a sul residuo
corrente, $\bfr_i=\bfb-A\bfx_i$,
\begin{equation}\label{rortlanc}
\bfu_j^\top \bfr_i = 0, \qquad j\le i,
\end{equation}

\no che \`e formalmente uguale alla (\ref{rort}).

\begin{oss}\label{cgfromlanc}
Se $A$ fosse sdp, e imponessimo la (\ref{rortlanc}) per ottenere
la nuova soluzione, riotterremmo, in pratica, il metodo dei
gradienti coniugati. Questo segue osservando che:

\begin{itemize}
\item la matrice tridiagonale $T$ e la matrice ortogonale $U$ della
fattorizzazione (\ref{fattlanc}) sono univocamente determinate (a
meno dei segni degli elementi extra-diagonali di $T$ e delle
colonne di $U$) dalla scelta del vettore iniziale $\bfu_1$ (vedi
Algoritmo~\ref{algo6});

\item il cambiamento di segno di un generico elemento extra-diagonale
$\bb_i$ equivale (dimostrarlo per esercizio) a cambiare di segno
ai vettori $$\bfu_{i+1},\dots,\bfu_n.$$ Evidentemente, questo non
altera la propriet\`a ~$U^\top U=I$, come anche la condizione di
ortogonalit\`a (\ref{rortlanc});

\item nel caso sdp (vedi Osservazione~\ref{cglanc}), il metodo dei
gradienti coniugati definisce implicitamente la fattorizzazione
(\ref{AZZT}), con le condizioni di ortogonalit\`a riscrivibili
nella forma (\ref{rortz}), che sono rispettivamente equivalenti
alle (\ref{fattlanc}) e (\ref{rortlanc}), dal momento che
$\bfu_1=\bfz_1$.

\end{itemize}
\end{oss}

Tuttavia, la condizione (\ref{rortlanc}) comporta un serio
inconveniente. Essa, infatti, risulta essere equivalente a
richiedere che $\bfx_i = \bfx_0+U_i\bfy_i$, con il vettore
$\bfy_i\in\RR^i$ tale che (vedi (\ref{invar1})-(\ref{AUi}), e posto
$\bb=\|\bfr_0\|$),
\begin{eqnarray}\nonumber
\bfo &=& U_i^\top (\bfb-A(\bfx_0+U_i\bfy_i)) = U_i^\top (\bfr_0
-U_{i+1}\hat{T}_i\bfy_i)\\[1mm] \label{Tising} &=& \bb
E_1^{(i)}-(I_i~\bfo)\hat{T}_i\bfy_i = \bb E_1^{(i)} -T_i\bfy_i.
\end{eqnarray}

\no Osserviamo, tuttavia che, se $A$ \`e solo simmetrica, allora
la matrice $T_i$ potrebbe essere singolare. Pertanto, risulta
essere pi\`u appropriato ricercare la nuova soluzione
(\ref{xilanc}) imponendo che la norma del corrispondente residuo
sia minimizzata. Il metodo cos\`\i\, ottenuto \`e denominato {\em
MINRES (MINimal RESidual method)}. Dalle (\ref{invar1}) e (\ref{xilanc}) segue,
quindi, che
\begin{equation}\label{xil} \bfx_i = \bfx_0 +U_i\bfy_i, \qquad \bfy_i\in\RR^i,\quad
i=1,2,\dots,\end{equation}

\no dove \begin{equation}\label {vil} \bfy_i
=\arg\min_{\bfy\in\RR^i} \|\bfb-A(\bfx_0+U_i\bfy)\|.\end{equation}

\no Considerando la (\ref{AUi}), e ponendo al solito
$\bb=\|\bfr_0\|$, dalla (\ref{u1lanc}) e dalle
(\ref{xil})--(\ref{vil}) segue che
\begin{eqnarray}\nonumber \min_{\bfy\in\RR^i}
\|\bfb-A(\bfx_0+U_i\bfy)\| &=& \min_{\bfy\in\RR^i}
\|\bfr_0-AU_i\bfy\| = \min_{\bfy\in\RR^i}
\|\bb\bfu_1-U_{i+1}\hat{T}_i\bfy\|\\ \label{bfvi} &=&
\min_{\bfy\in\RR^i} \|\bb
E_1^{(i+1)}-\hat{T}_i\bfy\|.\end{eqnarray}

\no A questo punto, osserviamo che $\hat{T}_i$ ha rango massimo,
ovvero $i$. Ne consegue che \`e fattorizzabile $QR$ nella forma
\begin{equation}\label{qrti}
\hat{T}_i = Q_i\hat{R}_i, \qquad \hat{R}_i = \pmatrix{c} R_i\\
\bfo^\top \endpmatrix_{i+1\times i},
\end{equation}

\no con $Q_i$ matrice ortogonale e $R_i$ triangolare superiore e
nonsingolare. Ne consegue che, denotato con
\begin{equation}\label{gi}\bfg_i =
\pmatrix{c}\bfg_{i1}\\g_{i2}\endpmatrix \equiv \bb
Q_i^\top E_1^{(i+1)},\end{equation}

\no la soluzione risulta essere
\begin{eqnarray}\nonumber
\bfy_i &\equiv& R_i^{-1}\bfg_{i1} =\arg\min_{\bfy\in\RR^i} \|\bb
E_1^{(i+1)}-Q_i\hat{R}_i\bfy\|\\ \label{vil2} &=&
\arg\min_{\bfy\in\RR^i} \sqrt{\|\bfg_{i1}-R_i\bfy\|^2+|g_{i2}|^2}.
\end{eqnarray}

\no Inoltre, con tale scelta si ottiene che
~$\|\bfr_i\|=|g_{i2}|$.~
Questo implica che possiamo calcolare la norma del residuo
senza calcolare esplicitamente n\'e $\bfr_i$ n\'e $\bfx_i$.
Pertanto, potremmo iterare fino ad avere la norma del residuo pi\`u
piccola di una prefissata tolleranza e, quindi, calcolare la
soluzione alla fine.

Esaminiamo, brevemente, il caso in cui, al generico passo $i$
dell'Algoritmo~\ref{algo6}, si abbia $\bfu_{i+1}=\bfo$, ovvero
$\bb_i=0$. In tal caso, come visto nella Osservazione~\ref{invar},
avremo determinato un sottospazio invariante di $A$, per cui la
(\ref{invar2}) vale in luogo della (\ref{AUi}). Si ottengono, di
conseguenza, i seguenti risultati, del tutto analoghi,
rispettivamente, ai Corollari~\ref{cor4.4} e \ref{cor4.5} visti
per il metodo dei gradienti coniugati.

\begin{teo}\label{teo6.2} $\bfu_{i+1}=\bfo \Rightarrow \bfx_i=\bfx^*$. \end{teo}

\proof Dalle (\ref{invar1})-(\ref{invar2}) segue che la matrice $T_i$ deve essere
nonsingolare, infatti
$$i = \rank(AU_i) \le\min\{ \rank(U_i),\rank(T_i)\}\le i.$$
Pertanto, vedi (\ref{bfvi}),
\begin{eqnarray*}
\min_{\bfy\in\RR^i} \|\bfb-A(\bfx_0+U_i\bfy)\| &=&
\min_{\bfy\in\RR^i} \|\bfr_0-AU_i\bfy\| = \min_{\bfy\in\RR^i}
\|\bb\bfu_1-U_iT_i\bfy\|\\ &=& \min_{\bfy\in\RR^i} \|\bb
E_1^{(i)}-T_i\bfy\|=0,\end{eqnarray*} scegliendo
$\bfy=\bb\,T_i^{-1}E_1^{(i)}$.\QED

\begin{cor}\label{cor6.1} $\bfx_i=\bfx^*$, per qualche $i\le n$.\end{cor}

\proof \`E sufficiente osservare che, se si arrivasse al passo
$n$, (vedi (\ref{xilanc})) $S_n\equiv\RR^n$.\QED

\begin{oss} Va tuttavia sottolineato che, come nel caso del
metodo dei gradienti coniugati, la propriet\`a di terminazione
del Corollario~\ref{cor6.1} potrebbe non valere pi\`u in aritmetica
finita.\end{oss}

\begin{oss}
Osserviamo che, nel caso in cui la matrice $A$ sia solo
simmetrica, l'Algoritmo~\ref{algo4} visto per il metodo dei
gradienti coniugati potrebbe essere non definito anche se la
matrice $T_i$ in (\ref{Tising}) fosse nonsingolare. Infatti (vedi
Esercizio~\ref{ex3}) esso risulta essere equivalente alla
fatto\-rizzazione $LU$ di tale matrice, che potrebbe non essere
definita, anche se $T_i$ fosse nonsingolare.
\end{oss}

\section{Implementazione efficiente di MINRES}\label{effminres}

La procedura definita dalle (\ref{xil})--(\ref{vil2}) ha, tuttavia,
un serio inconveniente, legato al fatto che, per calcolare $\bfx_i$ (in
corrispondenza dell'indice di terminazione) dalla (\ref{xil}),
\`e richiesta la memorizzazione di tutti i vettori $\{\bfu_j\}$ generati
(vedi (\ref{invar1})).
Questo fatto implica che il costo computazionale, in termini di
occupazione di memoria, potrebbe diventare elevato, quando $n$
\`e molto grande. Ci\`o \`e ancor pi\`u vero se si pensa che per
calcolare gli elementi della fattorizzazione, mediante
l'Algoritmo~\ref{algo6} sono necessari solo gli ultimi 2 vettori
generati.

Per ovviare a questo problema, vediamo di esaminare le
relazioni tra l'approssimazione al passo $i$-esimo, $\bfx_i$, e
quella al passo precedente, $\bfx_{i-1}$. Supponiamo, quindi, di
avere gi\`a calcolato, al passo $i-1$, quanto segue:

\begin{itemize}

\item $Q_{i-1}\hat{R}_{i-1}\equiv Q_{i-1}\pmatrix{c}R_{i-1}\\
\bfo^\top \endpmatrix = \hat{T}_{i-1}$;

\item $\bfg_{i-1} \equiv
\pmatrix{c}\bfg_{i-1,1}\\g_{i-1,2}\endpmatrix = \bb\,
Q_{i-1}^\top E_1^{(i)}$;

\item $\bfx_{i-1}\equiv \bfx_0 +U_{i-1}R_{i-1}^{-1}\bfg_{i-1,1}$.

\end{itemize}

\no Osserviamo anche che $Q_{i-1}$ sar\`a data dal prodotto di
$i-1$ matrici elementari di Givens di dimensione $i$:
\begin{equation}\label{Gim1}
Q_{i-1} = \left(G_i^{(i-1)}\cdots G_2^{(i-1)}\right)^\top .
\end{equation}

\no Inoltre, $R_{i-1}$ \`e triangolare superiore con al pi\`u 3
diagonali (vedi Osservazione~\ref{QRT}): $$R_{i-1} = \pmatrix{ccccc}
r_{11} &r_{12} &r_{13}\\
       &\ddots &\ddots &\ddots\\
       &       & \ddots &\ddots &r_{i-3,i-1}\\
       &       &        &\ddots &r_{i-2,i-1}\\
       &       &        &       &r_{i-1,i-1}\endpmatrix.$$

\no Al passo $i$-esimo, d'altronde, abbiamo che
\begin{eqnarray*}
\hat{T}_i &=& \pmatrix{c|c} \hat{T}_{i-1} & \begin{array}{c} 0\\
\vdots\\0\\ \bb_{i-1} \\ \cc_i\end{array}\\ \hline \bfo^\top  &
\bb_i\endpmatrix = \pmatrix{c|c} Q_{i-1}\hat{R}_{i-1} & \begin{array}{c} 0\\
\vdots\\0\\ \bb_{i-1} \\ \cc_i\end{array}\\ \hline \bfo^\top  &
\bb_i\endpmatrix \\
&=& \pmatrix{c|c} Q_{i-1}\\ \hline &1\endpmatrix
\pmatrix{c|c} \hat{R}_{i-1} & \bfw_i\\ \hline \bfo^\top  &
\bb_i\endpmatrix,
\end{eqnarray*}

\no dove $$\bfw_i = Q_{i-1}^\top \pmatrix{c}0\\ \vdots \\ 0 \\
\bb_{i-1}\\ \cc_i\endpmatrix \equiv \pmatrix{c} 0\\ \vdots\\
0\\ r_{i-2,i}\\r_{i-1,i}\\ w_i\endpmatrix.$$

\begin{oss} Osserviamo che, data la struttura del vettore moltiplicato, il
calcolo di $\bfw_i$ richiede solo le ultime due matrici di Givens,
$G_{i-1}^{(i-1)}$ e $G_i^{(i-1)}$, che danno $Q_{i-1}$ (vedi (\ref{Gim1})).
\end{oss}

\no Definendo una ulteriore matrice di Givens, $$G_{i+1}^{(i)} =
\pmatrix{c|rr} I_{i-1} \\ \hline &c_i & s_i\\ & -s_i
&c_i\endpmatrix,$$

\no tale che $$\pmatrix{rr}c_i & s_i\\ -s_i
&c_i\endpmatrix\pmatrix{c} w_i\\ \bb_i\endpmatrix =
\pmatrix{c} r_{ii}\\ 0\endpmatrix,$$

\no si ottiene, quindi, la fattorizzazione $QR$ di $\hat{T}_i$:
\begin{equation}\label{QRTi}
\hat{T}_i = \pmatrix{c|c} Q_{i-1}\\ \hline &1\endpmatrix
(G_{i+1}^{(i)})^\top  \pmatrix{c|c} \hat{R}_{i-1} &
\begin{array}{c} 0\\ \vdots \\ 0\\ r_{i-2,i}\\ r_{i-1,i}\\ r_{ii}
\end{array}
\\ \hline \bfo^\top  &
0\endpmatrix \equiv Q_i \hat{R}_i = Q_i \pmatrix{c} R_i\\ \bfo^\top \endpmatrix.
\end{equation}

\no Definiamo ora la matrice
\begin{equation}\label{Pconi}P_i \equiv (\bfp_1\dots\bfp_i) =
U_i R_i^{-1}.\end{equation}

\no Tenendo conto delle (\ref{gi})-(\ref{vil2}), questo implica che $$\bfx_i =\bfx_0+ P_i
(I_i~\bfo)Q_i^\top (\bb\,E_1^{(i+1)}).$$ Inoltre, dalla (\ref{Pconi}) segue che
$$P_i R_i = U_i.$$

\no Considerando le ultime colonne di ambo i membri, da quest'ultima
equazione e dalla (\ref{QRTi}) si ottiene che
$$\bfp_i = \frac{\bfu_i -r_{i-1,i}\bfp_{i-1}
-r_{i-2,i}\bfp_{i-2}}{r_{ii}}.$$

\no Pertanto, il nuovo vettore $\bfp_i$ si ottiene dagli ultimi 2
e dal vettore $\bfu_i$. Osserviamo, inoltre, che
$$R_i^{-1} = \pmatrix{c|c} R_{i-1}^{-1} & -(r_{ii} R_{i-1})^{-1}
\pmatrix{c} 0 \\ \vdots \\ 0 \\ r_{i-2,i}\\ r_{i-1,i}\endpmatrix
\\ \hline & r_{ii}^{-1}\endpmatrix,$$

\no e, pertanto, dalla (\ref{Pconi}) segue che
$$P_i \equiv (P_{i-1}~\bfp_i).$$

\no Questo ci permette di concludere, dalle
(\ref{xil})--(\ref{QRTi}), che
\begin{eqnarray*}
\bfx_i &=& \bfx_0 + U_iR_i^{-1}\bfg_{i1} = \bfx_0 + P_i
(I_i~\bfo)Q_i^\top  (\bb\,E_1^{(i+1)}) \\[1mm] &=&\bfx_0 +P_i
\pmatrix{c|c|c} I_{i-1} &\bfo&\bfo\\ \hline \bfo^\top  &1
&0\endpmatrix G_{i+1}^{(i)}\pmatrix{c|c} Q_{i-1}^\top \\ \hline
&1\endpmatrix \pmatrix{c} \bb\,E_1^{(i)}\\0\endpmatrix\\[1mm] &=&
\bfx_0 + P_{i-1}(I_{i-1}~\bfo)Q_{i-1}^\top (\bb E_1^{(i)}) +
(E_i^{(i+1)})^\top  Q_i^\top  (\bb\,E_1^{(i+1)})\bfp_i \\[1mm] &\equiv&
\bfx_{i-1} + \xi_i\bfp_i,
\end{eqnarray*}

\no dove (vedi (\ref{gi}))
$$\xi_i = (E_i^{(i+1)})^\top  Q_i^\top  (\bb\,E_1^{(i+1)}) \equiv
(E_i^{(i)})^\top  \bfg_{i1}$$ \`e l'ultima componente del vettore
$\bfg_{i1}$.

\bigskip

\no Gli argomenti su esposti, ci permettono di derivare
l'Algoritmo~\ref{algo7} per il metodo MINRES. Chiaramente, per
l'utilizzo in aritmetica finita, una versione
iterativa, invece che semi-iterativa, risulta essere, in realt\`a,
pi\`u appropriata.

\begin{eser} Scrivere un codice che implementi in modo efficiente
l'Algoritmo~\ref{algo7}.\end{eser}

\begin{figure}[p]
\begin{algo}\label{algo7} MINRES.\\ \rm
\fbox{\parbox{12cm}{
\begin{itemize}

\nulit sia assegnato un vettore iniziale $\bfx_0$

\nulit calcolo $\bfr_0 = \bfb-A\bfx_0$ ~e~ $\bb=\|\bfr_0\|$

\nulit inizializzo $\bb_0=0$, ~$\bfu_1=\bfr_0/\bb$, ~$g_0 = \bb$

\nulit per i=1,\dots,n:

\begin{itemize}

\nulit $\bfv_i=A\bfu_i$

\nulit $\cc_i = \bfu_i^\top \bfv_i$

\nulit $\bfv_i \leftarrow \bfv_i -\cc_i\bfu_i-\bb_{i-1}\bfu_{i-1}$

\nulit $\bb_i = \|\bfv_i\|$

\nulit $\bfp_i = \bfu_i$

\nulit $r_{i-1,i} = \bb_{i-1}$

\nulit $r_{ii} = \cc_i$

\nulit se $i>2$, $\pmatrix{c}r_{i-2,i}\\r_{i-1,i}\endpmatrix
\leftarrow
\pmatrix{rr} c_{i-2} &s_{i-2}\\ -s_{i-2} &c_{i-2}\endpmatrix
\pmatrix{c} 0\\ \bb_{i-1}\endpmatrix$

\nulit \qquad\qquad $\bfp_i \leftarrow \bfp_i
-r_{i-2,i}\bfp_{i-2}$, fine se

\nulit se $i>1$, $\pmatrix{c}r_{i-1,i}\\r_{ii}\endpmatrix
\leftarrow
\pmatrix{rr} c_{i-1} &s_{i-1}\\ -s_{i-1} &c_{i-1}\endpmatrix
\pmatrix{c} r_{i-1,i}\\ \cc_i \endpmatrix$

\nulit \qquad\qquad $\bfp_i \leftarrow \bfp_i
-r_{i-1,i}\bfp_{i-1}$, fine se

\nulit calcolo $c_i,s_i$ tali che $\pmatrix{c}r_{ii}\\0\endpmatrix
\leftarrow
\pmatrix{rr} c_i &s_i\\ -s_i &c_i\endpmatrix
\pmatrix{c} r_{ii}\\ \bb_i \endpmatrix$

\nulit $\bfp_i \leftarrow \bfp_i/r_{ii}$

\nulit $\pmatrix{c} \xi_i \\ g_i \endpmatrix =
\pmatrix{rr} c_i &s_i\\ -s_i &c_i\endpmatrix
\pmatrix{c} g_{i-1}\\ 0 \endpmatrix$

\nulit $\bfx_i = \bfx_{i-1} +\xi_i\bfp_i$

\nulit se $\bb_i=0$~~oppure~~$|g_i|\le tol$, $\bfx^*=\bfx_i$, esci, fine se

\nulit $\bfu_{i+1} = \bfv_i/\bb_i$

\end{itemize}

\nulit fine per

\end{itemize}

}}
\end{algo}
\end{figure}

Come esempio, si consideri il sistema lineare (\ref{sistsim}) in
cui
\begin{equation}\label{exsim}
A = T_1\otimes I_{10} + I_{10}\otimes T_1,\qquad T_1 =
\pmatrix{rrrr} 1 &-1 \\ -1 &\ddots &\ddots\\
   &\ddots &\ddots &-1\\
   &       &-1 &1\endpmatrix_{10\times 10},\end{equation}

\no la cui matrice dei coefficienti \`e simmetrica ma non
definita. Il termine noto \`e stato scelto in modo da avere come
soluzione il vettore $$\bfx^*=(1,\dots,1)^\top .$$

\begin{figure}[t]
\begin{center}
\includegraphics[width=13cm,height=10cm]{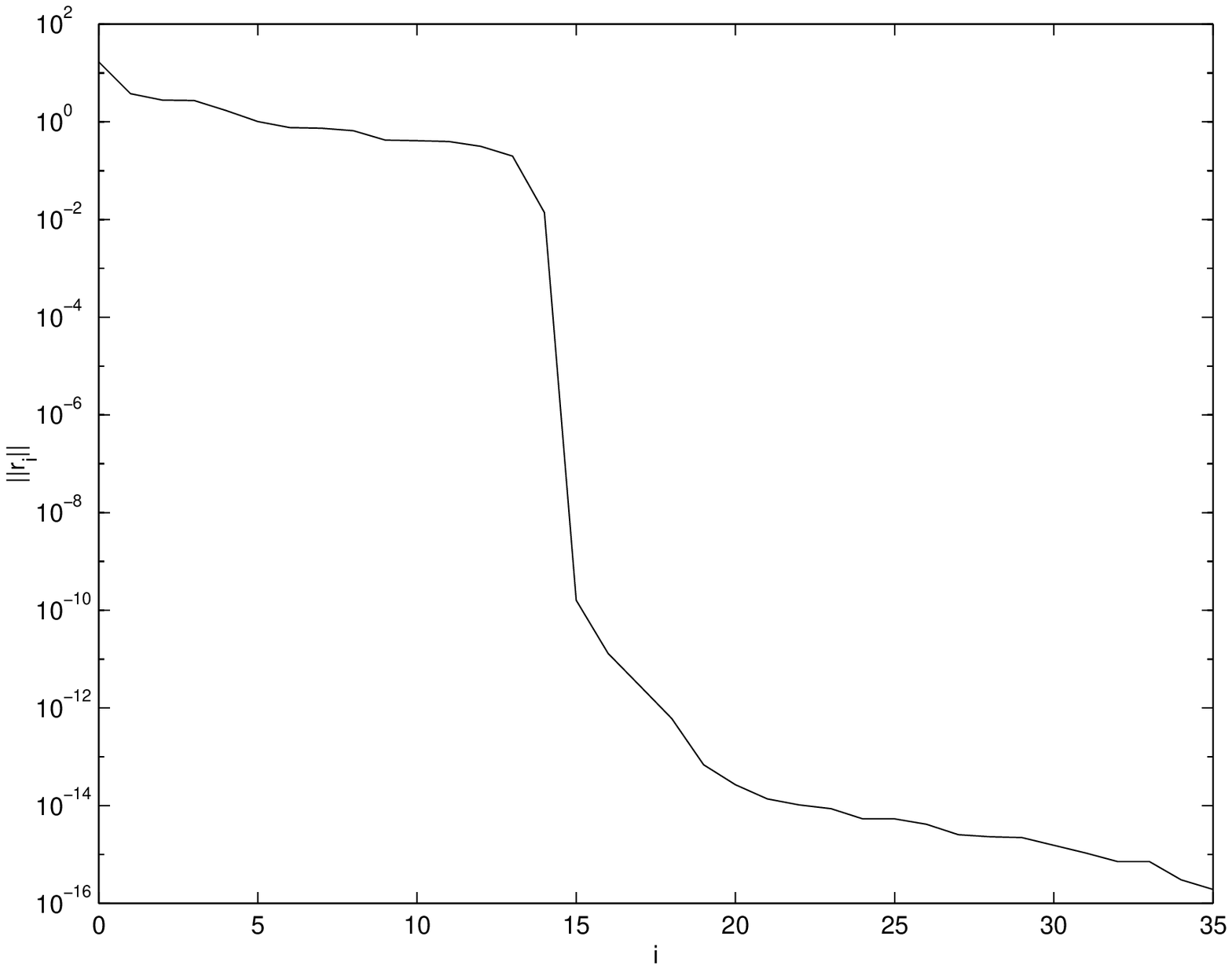}
\caption{\protect\label{fig10} MINRES applicato al problema
(\ref{sistsim})-(\ref{exsim}).}
\end{center}
\end{figure}

\no Applicando MINRES a partire dal vettore nullo, la soluzione
\`e determinata, entro il limite della precisione di macchina, in
35 iterazioni. In Figura~\ref {fig10} si riporta il grafico della
norma del residuo rispetto all'indice di iterazione.

%
%
\chapter{Il caso $A$ non simmetrica}\label{cap7}

{\em
Il caso in cui la matrice $A$ \`e non simmetrica (ma, ricordiamo,
nonsingolare) \`e quello pi\`u difficile tra quelli esaminati. In
effetti, vedremo che i metodi ite\-rativi di risoluzione o sono
particolarmente costosi o presentano la possibilit\`a di {\em
breakdown}. Quest'ultimo avviene quando la procedura non \`e pi\`u
definita, senza che la soluzione sia stata raggiunta.}

\section{Il metodo di Arnoldi}\label{Arnoldi}

Sia, dunque, \begin{equation}\label{Agen} A\bfx=\bfb\end{equation}

\no il sistema lineare da risolvere, la cui matrice dei
coefficienti, $A\in\RR^{n\times n}$, assumeremo essere nonsingolare.

Il primo metodo che andiamo ad esaminare \`e una generalizzazione
del metodo di Lanczos visto per il caso simmetrico: esso si basa
sulla trasformazione di $A$ in una matrice con struttura pi\`u
semplice, mediante una matrice ortogonale. In maggior dettaglio,
si ricerca la fattorizzazione
\begin{equation}\label{Arnold}
AU = UH \equiv(\bfu_1\dots\bfu_n)\pmatrix{rrrr}
h_{11}&\dots &\dots   &h_{1n}\\
h_{21} &      &        &\vdots\\
      &\ddots&         &\vdots\\
      &      &h_{n,n-1}&h_{nn}\endpmatrix, \qquad U^\top U=I,
      \end{equation}

\no ovvero la trasformazione di $A$ a forma di Hessemberg superiore.
Dalla (\ref{Arnold}) si ottiene che $U^\top AU=H$, che permette di
derivare il {\em metodo di Arnoldi}, di cui si riporta il
corrispondente Algoritmo~\ref{algo8}.

\begin{algo}\label{algo8} Metodo di Arnoldi.\\ \rm
\fbox{\parbox{12cm}{
\begin{itemize}

\nulit sia assegnato un vettore iniziale $\bfu_1, \|\bfu_1\|=1$

\nulit per i=1,\dots,n:

\begin{itemize}

\nulit $\bfv_i=A\bfu_i$

\nulit per j=1,\dots,i:

\begin{itemize}

\nulit $h_{ji} = \bfu_j^\top \bfv_i$

\nulit $\bfv_i \leftarrow \bfv_i-h_{ji}\bfu_j$

\end{itemize}

\nulit fine per

\nulit $h_{i+1,i} = \|\bfv_i\|$

\nulit se $h_{i+1,i}=0$, esci, fine se

\nulit $\bfu_{i+1} = \bfv_i/h_{i+1,i}$

\end{itemize}

\nulit fine per

\end{itemize}
}}\end{algo}

Si osserva che:

\begin{itemize}

\item come per tutte le fattorizzazioni esaminate nei
Capitoli~\ref{cap4}--\ref{cap6}, anche ora la fattorizzazione
(\ref{Arnold}) risulta essere univocamente determinata (avendo
fissato il segno degli elementi sottodiagonali di $H$) dalla
scelta del primo vettore, $\bfu_1$;

\item il costo dell'Algoritmo~\ref{algo8} \`e crescente con
l'indice di iterazione $i$: al passo $i$-esimo dovr\`o aver
memorizzato tutti i vettori $\bfu_1,\dots,\bfu_i$, e dovr\`o
eseguire 1 matvec, $i+1$ scal ed $i+1$ axpy;

\item se, per qualche $i<n$ dovesse risultare $h_{i+1,i}=0$ si
sarebbe determinato un sottospazio invariante destro per la
matrice $A$: dette $$
U_i = (\bfu_1\dots,\bfu_i), \qquad
H_i = \pmatrix{rrrr}
h_{11}&\dots &\dots   &h_{1i}\\
h_{21} &      &        &\vdots\\
      &\ddots&         &\vdots\\
      &      &h_{i,i-1}&h_{ii}\endpmatrix,
$$
esse verificano l'equazione
\begin{equation}\label{invArnold}
AU_i = U_iH_i.
\end{equation}
In particolare, gli autovalori di $H_i$ sono autovalori di $A$.
\end{itemize}

Relativamente all'ultimo punto, osserviamo che, qualora
$h_{i+1,i}\ne0$, allora, invece della (\ref{invArnold}), risulta
aversi, considerando le prime $i$ colonne dei membri della prima
equazione in (\ref{Arnold}), \begin{equation}\label{AUiHi}
AU_i = U_{i+1}\hat{H}_i\equiv U_{i+1} \pmatrix{cccc}
h_{11}&\dots &\dots   &h_{1i}\\
h_{21} &      &        &\vdots\\
      &\ddots&         &\vdots\\
      &      &h_{i,i-1}&h_{ii}\\ \hline & & & h_{i+1,i}\endpmatrix
      \equiv U_{i+1}\pmatrix{c} H_i \\ \hline h_{i+1,i}(E_i^{(i)})^\top \endpmatrix.
\end{equation}

\section{GMRES}\label{gmres}

Per utilizzare la fattorizzazione (\ref{Arnold}) al fine di risolvere il
sistema lineare (\ref{Agen}), le modalit\`a sono le stesse viste
per derivare MINRES dal metodo di Lanczos:

\begin{itemize}

\item data una approssimazione iniziale $\bfx_0$, si sceglie, quale vettore
iniziale,
\begin{equation}\label{u1arn}
\bfu_1 = \frac{\bfr_0}{\|\bfr_0\|}, \qquad \bfr_0 = \bfb-A\bfx_0;
\end{equation}

\item al passo $i$-esimo si ricerca
$$
\bfx_i\in S_i\equiv \bfx_0+[\bfu_1,\dots,\bfu_i],
$$

\no in modo da minimizzare la norma del corrispondente residuo, $\bfr_i=\bfb-A\bfx_i$.
\end{itemize}

\no Il metodo cos\`\i~ottenuto \`e denominato {\em GMRES
(Generalized Minimal RESi\-dual methods)}. In maggior dettaglio,
se $\bfx_i\in S_i$, allora si avr\`a
$$ \bfx_i = \bfx_0 + U_i\bfy_i,$$

\no con
$$\bfy_i = \arg\min_{\bfy\in\RR^i} \| \bfb -A(\bfx_0+U_i\bfy)\|.$$

\no Ponendo $\bb=\|\bfr_0\|$, dalle (\ref{AUiHi})-(\ref{u1arn}) si ottiene che

\begin{eqnarray}\nonumber
&&\min_{\bfy\in\RR^i} \| \bfb -A(\bfx_0+U_i\bfy)\| =\min_{\bfy\in\RR^i} \| \bfr_0-AU_i\bfy\|\\
 &&\qquad= \min_{\bfy\in\RR^i} \|\bb\bfu_1-U_{i+1}\hat{H}_i\bfy\|=\min_{\bfy\in\RR^i} \|
\bb\,E_1^{(i+1)}-\hat{H}_i\bfy\|. \label{minri}
\end{eqnarray}

\no A questo punto, ricordiamo che $\hat{H}_i$ ha rango massimo,
per cui sar\`a fattorizzabile nella forma
$$\hat{H}_i = Q_i\hat{R}_i \equiv Q_i \pmatrix{c} R_i \\
\bfo^\top \endpmatrix,$$

\no con $Q_i$ ortogonale e $R_i\in\RR^{i\times i}$ triangolare superiore e nonsingolare.
Ponendo $$\bfg_i = \bb Q_i^\top  E_1^{(i+1)} \equiv \pmatrix{c}
\bfg_{i1}\\ g_{i2}\endpmatrix,$$

\no dalla (\ref{minri}) si ottiene che
$$\|\bfr_i\| = \min_{\bfy\in\RR^i} \| Q_i(\bfg_i
-\hat{R}_i\bfy)\| = \min_{\bfy\in\RR^i} \| \bfg_i -\hat{R}_i \bfy\| =
|g_{i2}|,$$

\no in corrispondenza della scelta $\bfy = \bfy_i \equiv R_i^{-1}\bfg_{i1}$.
Osserviamo che la norma del residuo si pu\`o conoscere senza
calcolare il residuo stesso.

In modo del tutto analogo a quanto visto per MINRES (vedi il Teorema~\ref{teo6.2}
ed il Corollario~\ref{cor6.1}), si dimostrano, quindi, i seguenti risultati.

\begin{teo} $\bfu_{i+1}=\bfo ~\Rightarrow~ \bfx_i =
\bfx^*$, soluzione di (\ref{Agen}).\end{teo}

\begin{cor} GMRES converge (in aritmetica esatta) in al pi\`u $n$
ite\-razioni.\end{cor}

Al fine di ottenere un efficiente algoritmo per GMRES,
analizziamo la relazione esistente tra la fattorizzazione di
$\hat{H}_{i-1}$ e quella di $\hat{H}_i$. Supponiamo, quindi, di
aver ottenuto la fattorizzazione
$$\hat{H}_{i-1} = Q_{i-1} \hat{R}_{i-1} \equiv Q_{i-1}\pmatrix{c} R_{i-1} \\
\bfo^\top \endpmatrix,$$

\no in cui $Q_{i-1}$ sar\`a, al solito, esprimibile come il prodotto
di $i-1$ matrici elementari di Givens di dimensione $i$:
$$
Q_{i-1} = \left(G_i^{(i-1)}\cdots G_2^{(i-1)}\right)^\top .
$$

\no Si ottiene, quindi, che
\begin{eqnarray*}\hat{H}_i &=& \pmatrix{c|c} \hat{H}_{i-1} & \begin{array}{c}
h_{1i}\\ \vdots \\ h_{ii}\end{array}\\ \hline
&h_{i+1,i}\endpmatrix = \pmatrix{c|c} Q_{i-1}\hat{R}_{i-1} & \begin{array}{c}
h_{1i}\\ \vdots \\ h_{ii}\end{array}\\ \hline
&h_{i+1,i}\endpmatrix\\[1mm]
&=& \pmatrix{c|c} Q_{i-1}\\ \hline &1\endpmatrix \pmatrix{c|c} \hat{R}_{i-1}
&\begin{array}{c}
\bfr_{i1} \\w_i\end{array}\\ \hline
&h_{i+1,i}\endpmatrix,\end{eqnarray*}

\no dove $$\pmatrix{c} \bfr_{i1}\\ w_i \endpmatrix \equiv
Q_{i-1}^\top  \pmatrix{c} h_{1i}\\ \vdots \\ h_{ii}\endpmatrix.$$

\no Segue che, definendo la matrice di Givens
$$G_{i+1}^{(i)} = \pmatrix{c|rr} I_{i-1} \\
\hline &c_i & s_i\\ & -s_i &c_i\endpmatrix,$$

\no tale che $$\pmatrix{rr}c_i & s_i\\ -s_i
&c_i\endpmatrix\pmatrix{c} w_i\\ h_{i+1,i}\endpmatrix =
\pmatrix{c} r_{ii}\\ 0\endpmatrix,$$

\no si ottiene, infine, la fattorizzazione $QR$ di $\hat{H}_i$:
$$
\hat{H}_i = \pmatrix{c|c} Q_{i-1}\\ \hline &1\endpmatrix
(G_{i+1}^{(i)})^\top 
\pmatrix{c|c} \hat{R}_{i-1} &
\begin{array}{c} \bfr_{i1} \\ r_{ii}
\end{array}
\\ \hline \bfo^\top  &
0\endpmatrix \equiv Q_i \hat{R}_i = Q_i \pmatrix{c} R_i\\ \bfo^\top \endpmatrix.
$$

\begin{oss}\label{givm} Evidentemente, \`e possibile rappresentare
efficientemente la matrice $Q_i$ (ed effettuarne in modo
altrettanto efficiente i prodotti matrice-vettore e matrice
trasposta-vettore) mediante le coppie $(c_j,s_j),~j=1,\dots,i,$
che definiscono le corrispondenti matrici elementari di
Givens.\end{oss}

Questi argomenti, ci consentono di definire
l'Algoritmo~\ref{algo9} per GMRES.

\begin{eser} Scrivere un codice che implementi efficientemente
l'Algoritmo~\ref{algo9}, anche tenendo conto di quanto detto
nell'Osservazione~\ref{givm}.\end{eser}

\begin{oss} \`E possibile definire la versione precondizionata di
GMRES con le sostituzioni formali $$A\leftarrow CA, \qquad \bfb
\leftarrow C\bfb,$$ dove $C\approx A^{-1}$ \`e la matrice di
precondizionamento.\end{oss}

\begin{figure}[p]
\begin{algo}\label{algo9} GMRES.\\ \rm
\fbox{\parbox{12.5cm}{
\begin{itemize}
\setlength{\itemsep}{0cm}

\nulit sia assegnato un vettore iniziale $\bfx_0$

\nulit calcolo $\bfr_0 = \bfb-A\bfx_0$ ~e~ $\bb=\|\bfr_0\|$

\nulit inizializzo $U\equiv\bfu_1=\bfr_0/\bb$, ~$g_0 = \bb$, $Q_0^\top =1$

\nulit per i=1,\dots,n:

\begin{itemize}
\setlength{\itemsep}{0cm}

\nulit $\bfv_i=A\bfu_i$

\nulit per j=1,\dots,i:

\begin{itemize}

\nulit $h_{ji} = \bfu_j^\top \bfv_i$

\nulit $\bfv_i \leftarrow \bfv_i-h_{ji}\bfu_j$

\end{itemize}

\nulit fine per

\nulit $h_{i+1,i} = \|\bfv_i\|$

\nulit $\pmatrix{c}\bfr_{i1}\\ r_{ii}\endpmatrix = Q_{i-1}^\top \pmatrix{c}
h_{1i}\\ \vdots \\ h_{ii}\endpmatrix$

\nulit calcolo $c_i,s_i$ tali che $\pmatrix{c}r_{ii}\\0\endpmatrix
\leftarrow
\pmatrix{rr} c_i &s_i\\ -s_i &c_i\endpmatrix
\pmatrix{c} r_{ii}\\ h_{i+1,i} \endpmatrix$

\nulit se $i=1$, $R = (r_{11})$
\nulit \qquad altrimenti, $R \leftarrow \pmatrix{c|c} R &\bfr_{i1}\\ \hline \bfo^\top  & r_{ii}
\endpmatrix$
\nulit \qquad fine se

\nulit $\pmatrix{c} g_{i-1}\\ g_i\endpmatrix = \pmatrix{rr} c_i &s_i\\ -s_i &c_i\endpmatrix
\pmatrix{c} g_{i-1}\\ 0 \endpmatrix$

\nulit se $h_{i+1,i}=0$~~oppure~~$|g_i|\le tol$, esci, fine se

\nulit $\bfu_{i+1} = \bfv_i/h_{i+1,i}$

\nulit $Q_i^\top  =  \pmatrix{cc}
I_{i-1} \\ &\pmatrix{rr} c_i &s_i\\ -s_i
&c_i\endpmatrix\endpmatrix\pmatrix{cc} Q_{i-1}^\top \\ &1\endpmatrix$

\nulit $U \leftarrow (U ~\bfu_{i+1})$

\end{itemize}

\nulit fine per

\nulit risolvo $R\bfy = \bfg(0:i-1)$

\nulit $\bfx^* = \bfx_0 + U\bfy$

\end{itemize}
}}\end{algo}
\end{figure}

\begin{oss} Va sottolineato come, a causa degli errori dovuti
all'aritmetica finita, possa accadere che non si abbia convergenza
entro $n$ iterazioni. In questo caso, \`e prassi ripartire da $\bfx_n$
come soluzione iniziale ({\em restart}).\end{oss}

\begin{figure}[t]
\begin{center}
\includegraphics[width=13cm,height=10cm]{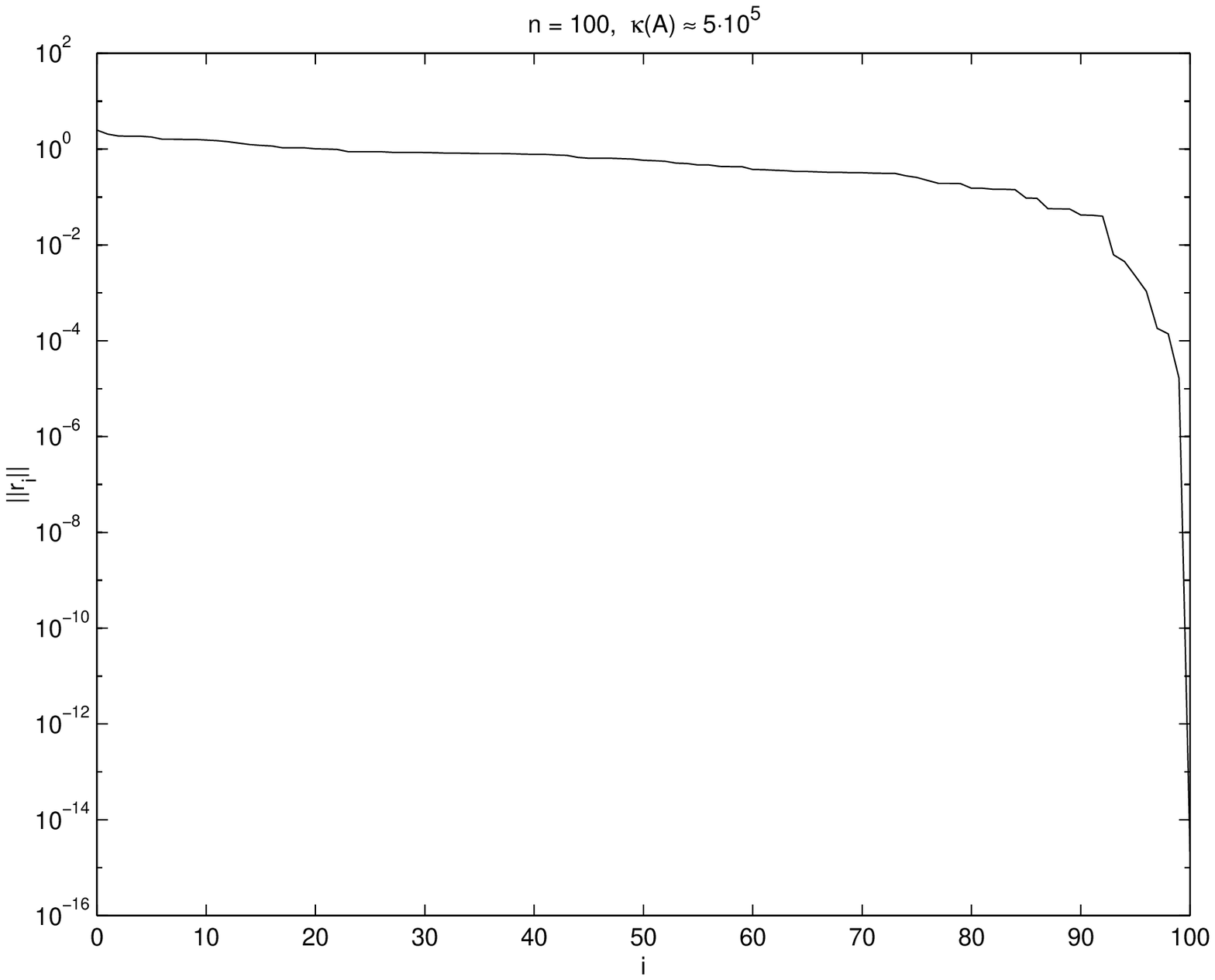}
\caption{\protect\label{fig11} GMRES applicato al problema
(\ref{Agen}).}
\end{center}
\end{figure}

Come esempio, consideriamo l'applicazione di GMRES al sistema
lineare (\ref{Agen}), la cui matrice $A$ ha dimensione $n=100$ ed
una densit\`a del $4\%$. Il corrispondente numero di condizione vale
approssimativamente $5\cdot10^5$. In Figura~\ref{fig11} \`e riportata la
norma del residuo rispetto all'indice di iterazione.

\bigskip
In ultimo, osserviamo che, se $n$ \`e molto grande, allora
l'occupazione di memoria ed il numero di operazioni
richiesti da GMRES potrebbero diventare troppo elevati. In sostanza,
il motivo di questo costo \`e legato alla struttura della matrice $H$ in
(\ref{Arnold}), che \`e di Hessemberg superiore. Il costo per iterazione era,
invece, assai pi\`u contenuto nel caso delle fattorizzazioni
(\ref{fattA}) e (\ref{fattlanc}), in cui la matrice $A$ era portata a
forma tridiagonale.

Per ovviare a questo problema, si potrebbe,
{\em in primis}, utilizzare siste\-maticamente il restart dopo un
numero (in genere assai esiguo) di iterazioni di GMRES. Tuttavia,
con questa, seppur semplice, strategia non \`e pi\`u possibile
assicurare, anche in aritmetica esatta, la propriet\`a di terminazione
finita del metodo.

Un'altra alternativa potrebbe essere quella di tentare di portare a
forma tridiagonale la matrice $A$ anche nel caso non simmetrico, sebbene
questo precluda, in pratica, l'utilizzo di una matrice di
trasformazione $U$ ortogonale. Proprio questo ultimo approccio
costituir\`a l'argomento delle prossime sezioni.

\section{Il metodo di Lanczos non simmetrico}\label{lancns}

Modificando la fattorizzazione (\ref{fattlanc}), per portare la
matrice $A$ a forma tridia\-gonale anche nel caso in cui essa non
sia simmetrica, si ottiene, considerando che in genere $U$ non \`e pi\`u
ortogonale, la seguente fattorizzazione:
\begin{eqnarray} \nonumber
AU&=&UT \equiv (\bfu_1 \dots \bfu_n)\pmatrix{cccc} \cc_1 &\bb_1\\
\aa_1 &\ddots &\ddots\\
      &\ddots &\ddots &\bb_{n-1}\\
      &       &\aa_{n-1}   &\cc_n\endpmatrix,\\ \label{lancuns}
A^\top  W &=& WT^\top  \equiv (\bfw_1\dots \bfw_n) T^\top ,\\[1mm] \nonumber
W^\top U &=& I.
\end{eqnarray}

\no Ad essa corrisponde il {\em metodo di Lanczos non simmetrico}.
In pratica, la prima equazione della (\ref{lancuns}) coincide con l'omologa
nella (\ref{fattlanc}); la seconda si ottiene riscrivendo la prima come $U^{-1}A = T U^{-1}$ e,
considerando le trasposte di ambo i membri, ponendo, infine, $W=U^{-T}$, come
sancito dall'ultima equazione.

Per ottenere la fattorizzazione (\ref{lancuns}) osserviamo che,
ponendo per convenzione $\aa_0=\bb_0=0$ e
\begin{equation}\label{lastzero}
\bfu_{n+1}=\bfw_{n+1}=\bfo,\end{equation}

\no le colonne $i$-esime delle prime due equazioni danno le seguenti relazioni:
\begin{eqnarray}\label{Auiwi}
A\bfu_i &=& \aa_i\bfu_{i+1} +\cc_i\bfu_i +\bb_{i-1}\bfu_{i-1},\\[1mm]
A^\top \bfw_i &=& \bb_i\bfw_{i+1} +\cc_i\bfw_i +\aa_{i-1}\bfw_{i-1},
\qquad i = 1,\dots,n.\nonumber
\end{eqnarray}

\no Inoltre, moltiplicando la prima equazione per $W^\top $, e tenendo conto della terza,
si ottiene la relazione $$W^\top AU = T,$$ da cui si ottengono le
seguenti espressioni per i coefficienti della matrice $T$:
\begin{eqnarray} \nonumber
\cc_i &=& \bfw_i^\top A\bfu_i, \qquad i=1,\dots,n,\\[1mm] \label{abc}
\bb_i &=& \bfw_i^\top A\bfu_{i+1},\\[1mm] \nonumber
\aa_i &=& \bfw_{i+1}^\top A\bfu_i,\qquad i=1,\dots,n-1.
\end{eqnarray}

\begin{algo}\label{algo10} Metodo di Lanczos non simmetrico.
\\ \rm
\fbox{\parbox{12cm}{
\begin{itemize}

\nulit siano assegnati due vettori $\bfr_0, \hat\bfr_0$ tali che
$\eta \equiv \bfr_0^\top \hat\bfr_0\ne0$

\nulit inizializzo $\bfu_1 = \bfr_0/\|\bfr_0\|$,~$\bfw_1 =
\hat\bfr_0/(\bfu_1^\top \hat\bfr_0)$,~$\bb_0=\aa_0=0$

\nulit per i=1,\dots,n:

\begin{itemize}

\nulit $\bfv_i = A\bfu_i$

\nulit $\cc_i = \bfw_i^\top \bfv_i$

\nulit se $i=n$, esci, fine se

\nulit $\hat\bfu_i = \bfv_i -\cc_i\bfu_i -\bb_{i-1}\bfu_{i-1}$

\nulit $\hat\bfw_i = A^\top \bfw_i -\cc_i\bfw_i-\aa_{i-1}\bfw_{i-1}$

\nulit $\aa_i = \|\hat\bfu_i\|$

\nulit se $\aa_i=0$, errore 1, esci, fine se

\nulit $\bfu_{i+1}=\hat\bfu_i/\aa_i$

\nulit $\bb_i = \bfu_{i+1}^\top \hat\bfw_i$

\nulit se $\bb_i=0$, errore 2, esci, fine se

\nulit $\bfw_{i+1}=\hat\bfw_i/\bb_i$

\end{itemize}

\nulit fine per

\end{itemize}
}}\end{algo}

\no Le (\ref{Auiwi})-(\ref{abc}) permettono di definire l'\,Algoritmo~\ref{algo10} per il metodo di Lanczos non
simmetrico.

\`E facile vedere che il precedente algoritmo corrisponde proprio
alla fatto\-rizzazione (\ref{lancuns}), a patto che i vettori
$\{\bfu_i\}$ e $\{\bfw_i\}$ siano bi-ortonormali, ovvero tali che
\begin{equation}\label{bi-ort}
\bfw_i^\top \bfu_i = 1, \qquad \bfw_j^\top \bfu_i=0, \quad i\ne j.
\end{equation}

\no Questo \`e sancito dal seguente risultato.

\begin{teo} I vettori $\{\bfu_i\}$ e $\{\bfw_i\}$ soddisfano le
(\ref{bi-ort}).\end{teo}

\proof La prima condizione \`e soddisfatta, avendosi, per definizione
di $\bb_i$,
$$\bfu_{i+1}^\top \bfw_{i+1} = \bfu_{i+1}^\top \frac{\hat\bfw_i}{\bb_i} =
\frac{\bfu_{i+1}^\top \hat\bfw_i}{\bfu_{i+1}^\top \hat\bfw_i}=1.$$

\no La seconda condizione in (\ref{bi-ort}) equivale alle
seguenti, $$\bfw_j^\top \bfu_i=0,\quad \bfu_j^\top \bfw_i=0, \qquad j<i,$$
\no che dimostriamo per induzione su $i$. Per $i=2$ si ha:
$$\bfw_1^\top \bfu_2 = \bfw_1^\top \hat\bfu_1\aa_1^{-1} =
\bfw_1^\top (A\bfu_1-\cc_1\bfu_1)\aa_1^{-1}=0,$$

\no in quanto $\cc_1=\bfw_1^\top A\bfu_1$ e $\bfw_1^\top \bfu_1=1$. Similmente,
$$\bfu_1^\top \bfw_2 = \bfu_1^\top \hat\bfw_1\bb_1^{-1} =
\bfu_1^\top (A^\top \bfw_1-\cc_1\bfw_1)\bb_1^{-1}=0.$$

\no Supponiamo ora vera la tesi per $i$ e dimostriamo per $i+1$.
Dimostriamo solo che $$\bfw_j^\top \bfu_{i+1}=0, \qquad j\le i,$$
in quanto l'altra relazione si dimostra analogamente. Possono
aversi i seguenti 3 casi:
\begin{description}

\item{$j<i-1:$~} in tal caso, dall'ipotesi di induzione segue
\begin{eqnarray*}
\bfw_j^\top \bfu_{i+1}&=&\bfw_j^\top ( A\bfu_i -\cc_i\bfu_i
-\bb_{i-1}\bfu_{i-1})\aa_i^{-1} = (\bfw_j^\top A)\bfu_i\aa_i^{-1}\\[1mm]
&=& (\bb_j\bfw_{j+1}+\cc_j\bfw_j+\aa_{j-1}\bfw_{j-1})^\top \bfu_i\aa_i^{-1} = 0;
\end{eqnarray*}

\item{$j=i-1:$~} in tal caso $$\bfw_{i-1}^\top \bfu_{i+1}=\bfw_{i-1}^\top ( A\bfu_i -\cc_i\bfu_i
-\bb_{i-1}\bfu_{i-1})\aa_i^{-1}=0,$$ per la definizione di
$\bb_{i-1}$ (vedi (\ref{abc}));

\item{$j=i:$~} in tal caso $$\bfw_i^\top \bfu_{i+1}=\bfw_i^\top ( A\bfu_i -\cc_i\bfu_i
-\bb_{i-1}\bfu_{i-1})\aa_i^{-1}=0,$$ per la definizione di
$\cc_i$ (vedi (\ref{abc})).
\end{description}\QED\\

\no La dimostrazione del seguente corollario \`e lasciata
come semplice esercizio.

\begin{cor}\label{cor7.2} Se i vettori $\{\bfu_i\}$ e $\{\bfw_i\}$ sono
bi-ortonormali, allora sia gli uni che gli altri sono linearmente
indipendenti.\end{cor}

\begin{oss} Dal precedente Corollario~\ref{cor7.2}, segue che la
convenzione (\ref{lastzero}) rispecchia un dato di fatto,
perch\`e, dopo $n$ passi, i vettori $\bfu_{n+1}$ e $\bfw_{n+1}$ sarebbero (in
aritmetica esatta) necessariamente nulli.\end{oss}

Riguardo alle condizioni di errore per l'uscita
dall'Algoritmo~\ref{algo10}, osservia\-mo che:

\begin{itemize}

\item la condizione denominata {\em errore 1}\, equivale ad aver
determinato un sottospazio invariante destro per la matrice $A$.
In questo caso, in modo del tutto analogo a quanto visto per il
metodo di Lanczos, si sar\`a ottenuta la soluzione, quando si
utilizzer\`a il metodo per risolvere il problema (\ref{Agen});

\item la condizione denominata {\em errore 2}\, \`e pi\`u seria ed
\`e legata, in generale, al problema di ridurre a forma
tridiagonale una generica matrice quadrata {\em (breakdown)}.
In tal caso, sarebbe possibile modificare l'intera procedura,
in modo da poter continuare l'iterazione {\em (look ahead)}, ma
per brevit\`a non tratteremo questo aspetto.

\end{itemize}

\begin{oss}\label{tro0} Nell'Algoritmo~\ref{algo10}, i vettori $\bfr_0$ e
$\hat\bfr_0$ sono scelti in modo da avere $\bfr_0^\top \hat\bfr_0\ne0$.
Un modo semplice per soddisfare questo requisito \`e scegliere
$\hat\bfr_0=\bfr_0\ne\bfo$.\end{oss}

\section{I gradienti bi-coniugati (Bi-CG)}\label{bicg}

Al fine di utilizzare la fattorizzazione (\ref{lancuns}) per
risolvere il problema (\ref{Agen}), data una approssimazione
iniziale $\bfx_0$, si sceglie
\begin{equation}\label{alsolito}
\bfr_0 = \bfb-A\bfx_0,
\end{equation}

\no il residuo iniziale, per definire il vettore iniziale $\bfu_1$
nell'Algoritmo~\ref{algo10} (e, seguendo quanto detto nella
Osservazione~\ref{tro0}, anche $\bfw_1$). Successivamente, al
passo $i$-esimo, si ricerca, al solito,
\begin{equation}\label{xinsi}
\bfx_i\in S_i =\bfx_0 +[\bfu_1,\dots,\bfu_i],
\end{equation}

\no in modo che siano soddisfatte le condizioni di ortogonalit\`a
(analoghe alle (\ref{rort}) viste per il metodo dei gradienti
coniugati)
\begin{equation}\label{rortbicg}
\bfw_j^\top \bfr_i = 0, \qquad j\le i.
\end{equation}

\no Le (\ref{alsolito})--(\ref{rortbicg}) definiscono il
{\em metodo dei gradienti bi-coniugati
(Bi-CG)}. Se la procedura pu\`o essere eseguita fino al passo $n$,
questo ne garantisce la terminazione, essendo i
vettori $\{\bfw_i\}$ linearmente indipendenti.
Tuttavia, vi sono vari motivi per cui questo pu\`o non avvenire.
Infatti, anche supponendo che nell'Algoritmo~\ref{algo10} non si
abbiano condizioni di errore, definendo le sottomatrici
$$U_i = (\bfu_1\dots,\bfu_i),\qquad W_i = (\bfw_1\dots\bfw_i),$$

\no e osservando che, dalla (\ref{lancuns}),
\begin{equation}\label{AUilancuns}
AU_i =
U_{i+1}\hat{T}_i \equiv U_{i+1}\pmatrix{cccc}
\cc_1 &\bb_1\\
\aa_1 &\ddots &\ddots\\
      &\ddots &\ddots &\bb_{i-1}\\
      &       &\aa_{i-1} &\cc_i\\ \hline
      &       &          &\aa_i\endpmatrix\equiv U_{i+1}\pmatrix{c} T_i \\
      \hline \aa_i (E_i^{(i)})^\top \endpmatrix,\end{equation}

\no le (\ref{xinsi})-(\ref{rortbicg}) equivalgono a ricercare
$\bfx_i=\bfx_0+U_i\bfy_i$, in modo tale che, definendo $\bb=\|\bfr_0\|$,
\begin{eqnarray}\nonumber
\bfo &=& W_i^\top \bfr_i = W_i^\top (\bfr_0-AU_i\bfy_i) = \bb E_1^{(i)}
-W_i^\top U_{i+1} \hat{T}_i \bfy_i\\[1mm] \label{bicgvi}
&=& \bb E_1^{(i)} -(I_i~\bfo)\hat{T}_i\bfy_i = \bb
E_1^{(i)}-T_i\bfy_i.
\end{eqnarray}

\no Tuttavia, potrebbe aversi $T_i$ singolare, anche se
$\hat{T}_i$ ha rango massimo.

\begin{oss} Nel caso in cui si arrivi al passo $i$-esimo e si
abbia $\bfu_{i+1}=\bfo$ (ovvero, la condizione di ``errore 1'' nell'Algoritmo~\ref{algo10}),
allora si avrebbe $AU_i=U_iT_i$, con $T_i$ nonsingolare. Questo implica che
$\bfx_i=\bfx^*$. Infatti:
\begin{eqnarray*}
\bfr_i &=& \bfb-A\bfx_i = \bfb-A(\bfx_0+U_i\bfy_i) = \bfr_0
-AU_i\bfy_i\\[1mm] &=& \bfr_0 -U_iT_i\bfy_i = U_i(\bb
E_1^{(i)}-T_i\bfy_i)=\bfo,\end{eqnarray*}

\no in virt\`u della (\ref{bicgvi}).
\end{oss}

Dai precedenti argomenti, discende il seguente risultato.

\begin{teo}\label{star} In aritmetica esatta, i Bi-CG convergono
alla soluzione, salvo breakdown, in al pi\`u $n$
iterazioni.\end{teo}

Infine, sulla scorta di quanto detto nella
Osservazione~\ref{cgfromlanc} riguardo alla derivazione del metodo dei
gradienti coniugati dal metodo di Lanczos, anche in questo caso
\`e possibile derivare un algoritmo, corrispondente
all'Algoritmo~\ref{algo4} visto per i gradienti coniugati, per il
metodo Bi-CG. Questo \`e descritto dall'Algoritmo~\ref{algo11}
(vedi Esercizi~\ref{bicges}--\ref{bicglanc}).

\begin{oss} In analogia con l'Algoritmo~\ref{algo5}, \`e possibile
definire la versione precondizionata del Bi-CG.\end{oss}

\begin{eser} Determinare il costo computazionale
dell'Algoritmo~\ref{algo11}.\end{eser}

\begin{algo}\label{algo11} Metodo dei gradienti bi-coniugati (Bi-CG).\\ \rm
\fbox{\parbox{12cm}{
\begin{itemize}

\nulit sia assegnata una approssimazione iniziale $\bfx_0$

\nulit calcolo $\bfr_0 = \bfb-A\bfx_0$

\nulit inizializzo $\hat\bfr_0=\bfp_1=\hat\bfp_1=\bfr_0$, ~$\eta_0=\hat\bfr_0^\top \bfr_0$

\nulit per i=1,\dots,n:

\begin{itemize}

\nulit $\bfv_i=A\bfp_i$

\nulit $\hat\bfv_i=A^\top \hat\bfp_i$

\nulit $\hat{d}_i = \hat\bfp_i^\top \bfv_i$

\nulit $\hat\lam_i = \eta_{i-1}/\hat{d}_i$

\nulit $\bfx_i = \bfx_{i-1}+\hat\lam_i\bfp_i$

\nulit $\bfr_i = \bfr_{i-1}-\hat\lam_i\bfv_i$

\nulit $\hat\bfr_i = \hat\bfr_{i-1}-\hat\lam_i\hat\bfv_i$

\nulit $\eta_i = \hat\bfr_i^\top \bfr_i$

\nulit $\mu_i = \eta_i/\eta_{i-1}$

\nulit $\bfp_{i+1} = \bfr_i+\mu_i\bfp_i$

\nulit $\hat\bfp_{i+1} = \hat\bfr_i+\mu_i\hat\bfp_i$

\end{itemize}

\nulit fine per

\nulit $\bfx^*=\bfx_n$
\end{itemize}

}}\end{algo}
\medskip

\begin{eser} Scrivere l'algoritmo corrispondente alla versione
precondizio\-nata del Bi-CG.\end{eser}

\begin{eser} Scrivere un codice che implementi efficientemente
l'Algoritmo~\ref{algo11}, inserendo un idoneo criterio di arresto
sul residuo ed un controllo sul breakdown. \end{eser}

\begin{oss} Se $A$ \`e sdp, gli Algoritmi~\ref{algo4} e
\ref{algo11} forniscono, in arit\-metica esatta, la stessa
successione di approssimazioni, in quanto, in tal caso, essendo
$A=A^\top $, le successioni $\{\hat\bfr_i\}$, $\{\bfr_i\}$, e
$\{\hat\bfp_i\}$, $\{\bfp_i\}$, coincidono.\end{oss}

\begin{eser}\label{bicges} Se l'Algoritmo~\ref{algo11} pu\`o essere eseguito
regolarmente fino al passo $i$-esimo, dimostrare che:
\begin{enumerate}
\item $\bfr_i = \sum_{j=0}^i b_j^{(i)}A^j\bfp_1, \quad\quad \hat\bfr_i = \sum_{j=0}^i
b_j^{(i)}(A^\top )^j\hat\bfp_1$, \\ [3mm]
$\bfp_{i+1} = \sum_{j=0}^i c_j^{(i)}A^j\bfp_1, \quad \hat\bfp_{i+1} = \sum_{j=0}^i
c_j^{(i)}(A^\top )^j\hat\bfp_1$, \quad con
\begin{equation}\label{maincoef}
b_i^{(i)} = c_i^{(i)} = (-1)^i
\prod_{r=1}^i\hat\lam_r\,;\end{equation}

\item $[\bfr_0,\dots,\bfr_i] = [\bfp_1,\dots,\bfp_{i+1}] =
[\bfp_1,A\bfp_1,\dots,A^i\bfp_1]$,\\ [3mm]
 $[\hat\bfr_0,\dots,\hat\bfr_i] = [\hat\bfp_1,\dots,\hat\bfp_{i+1}] =
[\hat\bfp_1,A^\top \hat\bfp_1,\dots,(A^\top )^i\hat\bfp_1]$;

\item $\hat\bfp_j^\top \bfr_i = \bfp_j^\top \hat\bfr_i=0,\qquad
\hat\bfp_j^\top A\bfp_{i+1} = \bfp_j^\top A^\top \hat\bfp_{i+1}=0, \qquad j\le i$.
\end{enumerate}
\end{eser}

\begin{eser} Dimostrare che, se $T$ \`e la matrice in
(\ref{lancuns}) e
\begin{equation}\label{effe}
F=\pmatrix{ccc}f_1\\&\ddots\\&&f_n\endpmatrix,\quad f_1=1,
\quad f_{i+1}=f_i\frac{\bb_i}{\aa_i}, \quad i=1,\dots,n-1,~~~
\end{equation}

\no allora \quad $F\,TF^{-1}=T^\top $.
\end{eser}

\begin{eser} Supponendo che l'Algoritmo~\ref{algo11} possa essere
eseguito regolarmente per $n$ passi, dimostrare che esso equivale ad ottenere
la fattorizzazione:
\begin{eqnarray}\nonumber
AP &=& RL\Lambda^{-1}\equiv P\hat{T},\\
A^\top \hat{P} &=& \hat{R}L\Lambda^{-1}\equiv\hat{P}\hat{T},
\label{bicgfat}\\  \nonumber
\hat{P}^\top AP &=& \hat{D},
\end{eqnarray}

\no dove (confrontare con (\ref{APPT}) e (\ref{Tcappello}))
\begin{eqnarray*}P&=&(\bfp_1,\dots,\bfp_n),\qquad
R=(\bfr_0,\dots,\bfr_{n-1}),\qquad
R=PB,\\[1mm]
\hat{P}&=&(\hat\bfp_1,\dots,\hat\bfp_n),\qquad
\hat{R}=(\hat\bfr_0,\dots,\hat\bfr_{n-1}), \qquad
 \hat{R}=\hat{P}B,
\end{eqnarray*}
$$B = \pmatrix{cccc}
1&-\mu_1\\
 &\ddots &\ddots\\
 &       &\ddots &-\mu_{n-1}\\
 &       &       &1\endpmatrix,\qquad
 L = \pmatrix{cccc}
 1\\ -1&\ddots\\
 &\ddots &\ddots\\&&-1&1\endpmatrix,$$
 $$ \Lambda = \pmatrix{ccc}
 \hat\lam_1\\ &\ddots\\ &&\hat\lam_n\endpmatrix,
 \qquad \hat{D}=\pmatrix{ccc}
 \hat{d}_1\\ &\ddots\\ &&\hat{d}_n\endpmatrix,$$
 $$\hat{T}=BL\Lambda^{-1},\qquad \hat{R}^\top R = \Sigma \equiv \pmatrix{ccc} \eta_0\\
 &\ddots\\&&\eta_{n-1}\endpmatrix.$$
\end{eser}

\begin{eser}\label{bicglanc} Dimostrare che la fattorizzazione
(\ref{bicgfat}) equivale alla fatto\-rizzazione (\ref{lancuns})
(vedi l'Algoritmo~\ref{algo10}), con
$$U = R \Delta^{-1}, \qquad W = \hat{R}\Delta\Sigma^{-1}, \qquad
\Delta = \pmatrix{ccc} \|\bfr_0\|\\ &\ddots\\
&&\|\bfr_{n-1}\|\endpmatrix.$$
(Confrontare con l'Esercizio~\ref{ex3} e l'Osservazione~\ref{cglanc}).
\end{eser}

\begin{figure}[t]
\begin{center}
\includegraphics[width=13cm,width=10cm]{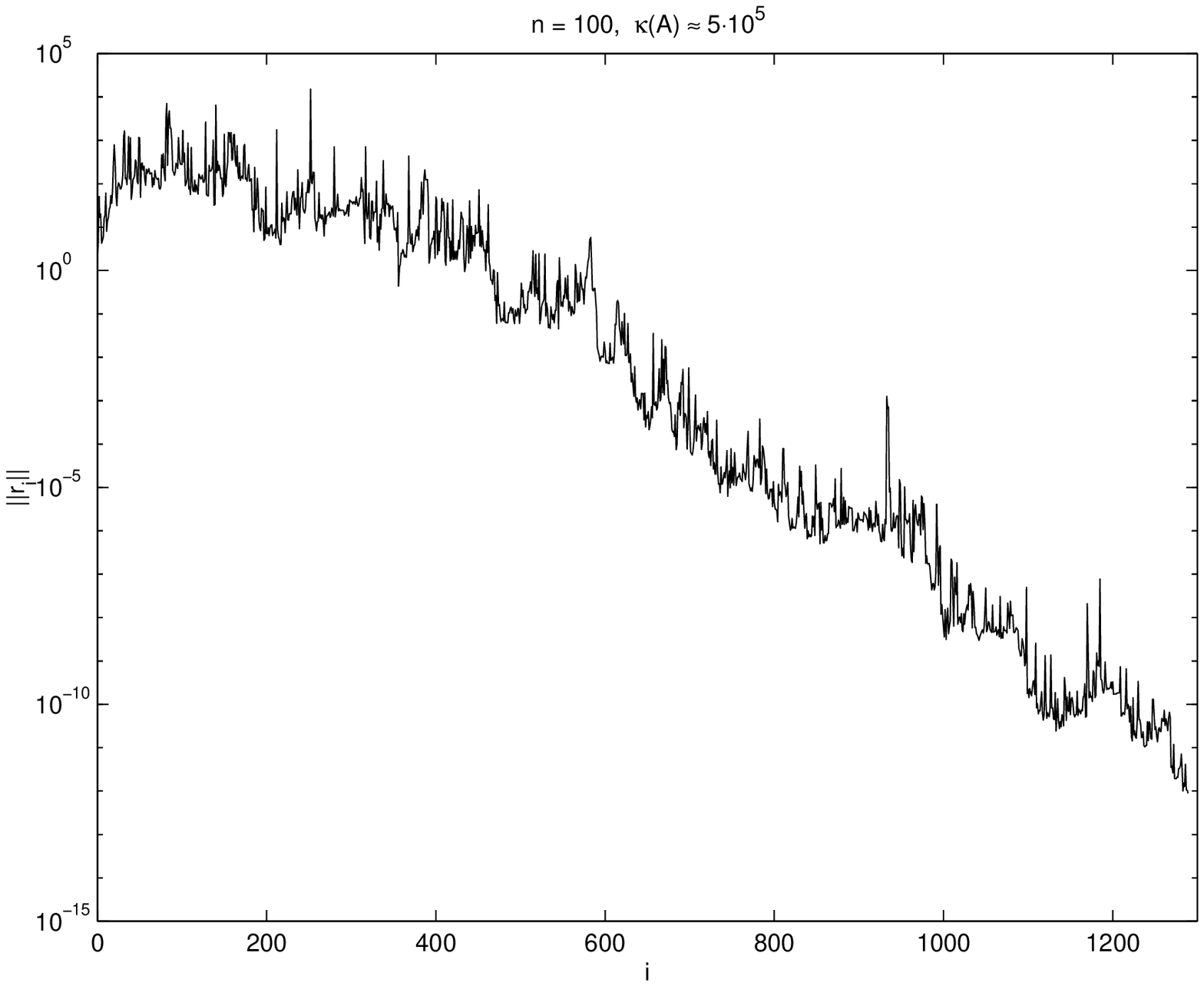}
\caption{\protect\label{fig12} Bi-CG applicato al problema
(\ref{Agen}).}
\end{center}
\end{figure}

Per concludere questa sezione, consideriamo l'applicazione del
Bi-CG al sistema lineare (\ref{Agen}), di dimensione $n=100$, densit\`a
del $4\%$ e numero di condizione $\approx5\cdot10^5$, su cui abbiamo gi\`a
testato GMRES (vedi Figura~\ref{fig11}). In Figura~\ref{fig12}
\`e riportata la norma del residuo rispetto all'indice di iterazione.
Sono ora necessarie ben 1289 iterazioni per avere il residuo con norma
inferiore a $10^{-12}$. Come si vede, in aritmetica finita la propriet\`a
di terminazione del Teorema~\ref{star} pu\`o non essere pi\`u valida.

\section{QMR}\label{qmr}

Il metodo che andiamo ora ad esaminare sta al Bi-CG come MINRES
sta ai gradienti coniugati. Infatti, sempre considerando al passo
$i$-esimo la soluzione come in (\ref{xinsi}), la si ricerca in
modo da minimizzare una maggiorazione della norma del residuo corrente.
I passi sono assai simili a quelli visti per MINRES. Si avr\`a $\bfx_i = \bfx_0
+U_i\bfy_i$, dove si vorrebbe (vedi (\ref{AUilancuns}) e, al solito, ponendo
$\bb=\|\bfr_0\|$):
\begin{eqnarray*}
\bfy_i &=& \arg\min_{\bfy\in\RR^i}\| \bfb-A(\bfx_0+U_i\bfy)\| =
           \arg\min_{\bfy\in\RR^i}\|\bfr_0-U_{i+1}\hat{T}_i\bfy\|\\
        &=&\arg\min_{\bfy\in\RR^i}\|U_{i+1}(\bb E_1^{(i+1)}-\hat{T}_i\bfy)\|.
\end{eqnarray*}

\no A questo punto, se si fattorizzasse $QR$ la matrice $U_{i+1}$, si
otterrebbe nuovamente GMRES, in quanto il prodotto del fattore $R$
(triangolare superiore) con $\hat{T}_i$ (tridiagonale), darebbe una matrice
di Hessemberg superiore che, a meno dei segni degli elementi sottodiagonali,
 risulta essere proprio la sottomatrice
$\hat{H}_i$ (vedi (\ref{AUiHi})) della matrice $H$ della
fattorizzazione (\ref{Arnold}).
Invece, al fine di evitare i problemi di costo computazionale visti per GMRES, si osserva che
$$\|U_{i+1}(\bb E_1^{(i+1)}-\hat{T}_i\bfy)\|\le \|U_{i+1}\|\cdot
\|\bb E_1^{(i+1)}-\hat{T}_i\bfy\|,$$

\no per cui si sceglier\`a
\begin{equation}\label{qmrvi}
\bfy_i = \arg\min_{\bfy\in\RR^i}\|\bb
E_1^{(i+1)}-\hat{T}_i\bfy\|.
\end{equation}

\no Il corrispondente metodo \`e denominato {\em QMR (Quasi
Minimal Residual me\-thod)}. La minimizzazione dell'argomento della
norma in (\ref{qmrvi}) si ottiene in un modo del tutto analogo a
quello visto per MINRES, con tutti i dettagli implementativi (con
le opportune, minime, modifiche) esaminati in
Sezione~\ref{effminres}. Questo permette di derivare il corrispondente
Algoritmo~\ref{algo12}.

La propriet\`a di terminazione finita, salvo breakdown, continua a
valere anche per QMR.

Osserviamo infine che, considerando la eventuale matrice
precondizionata $CA$ al posto di $A$, si ottiene la
versione precondizionata del metodo.

\begin{figure}[hp]
\begin{algo}\label{algo12} QMR.\\ \rm
\fbox{\parbox{12cm}{
\begin{itemize}
\setlength{\itemsep}{0cm}
\nulit sia assegnata una approssimazione iniziale $\bfx_0$

\nulit calcolo $\bfr_0 = \bfb-A\bfx_0$~ e ~$\bb=\|\bfr_0\|$

\nulit inizializzo
$\bfw_1=\bfu_1=\bfr_0/\bb$,~$\bb_0=\aa_0=0$, ~$g_0=\bb$

\nulit per i=1,\dots,n:

\begin{itemize}
\setlength{\itemsep}{0cm}

\nulit $\bfv_i = A\bfu_i$

\nulit $\cc_i = \bfw_i^\top \bfv_i$

\nulit $\hat\bfu_i = \bfv_i -\cc_i\bfu_i -\bb_{i-1}\bfu_{i-1}$

\nulit $\aa_i = \|\hat\bfu_i\|$

\nulit $\bfp_i = \bfu_i$

\nulit $r_{i-1,i}=\bb_{i-1}$

\nulit $r_{ii} = \cc_i$

\nulit se $i>2$, $\pmatrix{c}r_{i-2,i}\\r_{i-1,i}\endpmatrix
\leftarrow
\pmatrix{rr} c_{i-2} &s_{i-2}\\ -s_{i-2} &c_{i-2}\endpmatrix
\pmatrix{c} 0\\ \bb_{i-1}\endpmatrix$

\nulit \qquad\qquad $\bfp_i \leftarrow \bfp_i
-r_{i-2,i}\bfp_{i-2}$, fine se

\nulit se $i>1$, $\pmatrix{c}r_{i-1,i}\\r_{ii}\endpmatrix
\leftarrow
\pmatrix{rr} c_{i-1} &s_{i-1}\\ -s_{i-1} &c_{i-1}\endpmatrix
\pmatrix{c} r_{i-1,i}\\ \cc_i \endpmatrix$

\nulit \qquad\qquad $\bfp_i \leftarrow \bfp_i
-r_{i-1,i}\bfp_{i-1}$, fine se

\nulit calcolo $c_i,s_i$ tali che $\pmatrix{c}r_{ii}\\0\endpmatrix
\leftarrow
\pmatrix{rr} c_i &s_i\\ -s_i &c_i\endpmatrix
\pmatrix{c} r_{ii}\\ \aa_i \endpmatrix$

\nulit $\bfp_i \leftarrow \bfp_i/r_{ii}$

\nulit $\pmatrix{c} \xi_i \\ g_i \endpmatrix =
\pmatrix{rr} c_i &s_i\\ -s_i &c_i\endpmatrix
\pmatrix{c} g_{i-1}\\ 0 \endpmatrix$

\nulit $\bfx_i = \bfx_{i-1} +\xi_i\bfp_i$

\nulit se $\aa_i=0$~~oppure~~$|g_i|\le tol$, $\bfx^*=\bfx_i$, esci, fine se

\nulit $\bfu_{i+1}=\hat\bfu_i/\aa_i$

\nulit $\hat\bfw_i = A^\top \bfw_i -\cc_i\bfw_i-\aa_{i-1}\bfw_{i-1}$

\nulit $\bb_i = \bfu_{i+1}^\top \hat\bfw_i$

\nulit se $\bb_i=0$, breakdown, esci, fine se

\nulit $\bfw_{i+1}=\hat\bfw_i/\bb_i$

\end{itemize}

\nulit fine per

\end{itemize}

}}\end{algo}
\end{figure}

\begin{eser} Determinare il costo computazionale
dell'Algoritmo~\ref{algo12}.\end{eser}

\begin{eser} Scrivere la versione precondizionata
dell'Algoritmo~\ref{algo12}.\end{eser}

\section{Implementazione alternativa del QMR}\label{QMRimp}

Una implementazione pi\`u efficiente del QMR si ottiene
riscrivendo la seconda equazione in (\ref{lancuns}) nella forma
equivalente (vedi (\ref{effe})) $A^\top W=WFTF^{-1}$, ovvero $A^\top WF =
WFT$. Ponendo, quindi, $V=WF\equiv(\bfv_1\dots\bfv_n)$, si ottiene
la fattorizzazione equivalente $$AU=UT, \qquad
A^\top V = VT, \qquad V^\top U=F.$$

\no A questo punto, supponendo che la matrice $T$ sia
fattorizzabile $LU$, si ottiene
\begin{equation}\label{TLH}
T = LH \equiv \pmatrix{cccc} \ell_1\\ \aa_1
&\ell_2\\&\ddots&\ddots\\&&\aa_{n-1} &\ell_n\endpmatrix
\pmatrix{cccc} 1 &\phi_1\\ &\ddots&\ddots\\&&\ddots &\phi_{n-1}\\&&&1
\endpmatrix,\end{equation}

\no dove (vedi (\ref{lancuns})) \quad$\ell_1=\cc_1, \quad \phi_i =
\bb_i/\ell_i, ~ \ell_{i+1} = \cc_{i+1}-\phi_i\aa_i, ~
i=1,\dots,n-1$. Pertanto, ponendo $$Q = UH^{-1} \equiv
(\bfq_1\dots\bfq_n), \qquad Z =
VH^{-1}\equiv(\bfz_1\dots,\bfz_n),$$

\no si ottiene, infine, la fattorizzazione
\begin{equation}\label{qmrf1}
AQ = UL, \quad U = QH, \quad A^\top Z = VL, \quad V = ZH, \quad
V^\top U=F.\end{equation}

\no Per poter univocamente determinare la fattorizzazione, \`e
tuttavia necessaria una ulteriore condizione, che andiamo adesso a
definire.

\begin{eser} Dimostrare che (vedi (\ref{effe}) e (\ref{TLH}))
$$F^{-1}H^\top F = L L_D^{-1} \equiv L \pmatrix{ccc}\ell_1\\
&\ddots\\&&\ell_n\endpmatrix^{-1}.$$
\end{eser}

Dalla (\ref{qmrf1}), tenendo conto del risultato nel precedente
esercizio, si ottiene:
\begin{eqnarray*}
Z^\top AQ &=& H^{-T}V^\top AQ = H^{-T}V^\top UL = H^{-T}FL \\[1mm]
&=& F(F^{-1}H^{-T}F)L = F L_D \equiv D.\end{eqnarray*}

\no Pertanto, alle (\ref{qmrf1}) associamo l'ulteriore equazione
\begin{equation}\label{qmrf2}
Z^\top AQ = D \equiv\pmatrix{ccc} f_1\ell_1\\
&\ddots\\&&f_n\ell_n\endpmatrix.\end{equation}

\no Dalle (\ref{effe}), (\ref{TLH})--(\ref{qmrf2}), si ottiene,
quindi, l'Algoritmo~\ref{qmrfatt}, dove abbiamo posto
$\bfv_1=\bfu_1$.

\begin{figure}[t]
\begin{algo}\label{qmrfatt} Fattorizzazione (\ref{qmrf1})-(\ref{qmrf2}).
\\ \rm
\fbox{\parbox{12cm}{
\begin{itemize}

\nulit sia assegnato il vettore $\bfu_1$, $\|\bfu_1\|=1$

\nulit inizializzo $\bfv_1 = \bfq_1 = \bfz_1 = \bfu_1$, $f_1=1$,
$\ell_1 = \bfz_1^\top A\bfq_1$

\nulit per i=1,\dots,n-1:

\begin{itemize}

\nulit $\hat\bfu_i = A\bfq_i -\ell_i\bfu_i$

\nulit $\aa_i = \|\hat\bfu_i\|$

\nulit se $\aa_i=0$, errore 1, fine se

\nulit $\bfu_{i+1} = \hat\bfu_i/\aa_i$

\nulit $\bfv_{i+1} = ( A^\top \bfz_i -\ell_i\bfv_i )/\aa_i$

\nulit $f_{i+1} = \bfv_{i+1}^\top \bfu_{i+1}$

\nulit se $f_{i+1}=0$, errore 2, fine se

\nulit $\bb_i = \aa_i f_{i+1}/f_i$

\nulit $\bfq_{i+1} = \bfu_{i+1} -(\bb_i/\ell_i)\bfq_i$

\nulit $\bfz_{i+1} = \bfv_{i+1} -(\bb_i/\ell_i)\bfz_i$

\nulit $\ell_{i+1} = \bfz_{i+1}^\top A\bfq_{i+1}/f_{i+1}$

\nulit se $\ell_{i+1}=0$, errore 3, fine se

\end{itemize}

\nulit fine per

\end{itemize}
}}\end{algo}
\end{figure}

\begin{oss} Nell'Algoritmo~\ref{qmrfatt}, le condizioni denominate
~{\em errore 1}~ ed ~{\em errore 2}~ sono equivalenti alle condizioni
omonime viste nell'Algoritmo~\ref{algo10}. La condizione
denominata ~{\em errore 3}~ \`e invece legata alla possibile non
esistenza della fattorizzazione (\ref{TLH}).
\end{oss}

\begin{oss} Anche per la fattorizzazione di Lanczos (\ref{fattlanc}),
vista nel caso simmetrico, si potrebbe definire una variante
analoga alla fattorizzazione (\ref{qmrf1})-(\ref{qmrf2}). Questa
determinerebbe una corrispondente variante del metodo MINRES.\end{oss}

Utilizziamo, ora, la fattorizzazione (\ref{qmrf1})-(\ref{qmrf2})
per risolvere il sistema lineare (\ref{Agen}). Come visto per il
QMR, al passo $i$-esimo si ricercher\`a la approssimazione corrente nella
forma $\bfx_i=\bfx_0+U_i\bfy_i$, in cui il vettore $\bfy_i$ \`e
soluzione del problema di minimo (\ref{qmrvi}). Considerando che,
dalla (\ref{TLH}), $$\hat{T}_i = \hat{L}_iH_i \equiv \pmatrix{ccc}
\ell_1\\
\aa_1 &\ddots\\
      &\ddots &\ell_i\\
      &       &\aa_i\endpmatrix
      \pmatrix{cccc}
 1 &\phi_1\\
   &\ddots&\ddots\\
   &      &\ddots &\phi_{i-1}\\
   &      &       &1
\endpmatrix,$$

\no e ponendo $\bfc=H_i\bfy$, si ottiene quindi che
\begin{equation}\label{newmin}
\min_{\bfy\in\RR^i}\|\bb E_1^{(i+1)}-\hat{T}_i\bfy\| = \min_{\bfc\in\RR^i}\|\bb
E_1^{(i+1)}-\hat{L}_i\bfc\|.
\end{equation}

\no Se $\bfc_i$ \`e la soluzione del secondo problema, segue
 che $\bfy_i = H_i^{-1}\bfc_i$ \`e soluzione del primo e, pertanto,
$$\bfx_i = \bfx_0 + U_i\bfy_i = \bfx_0 +U_iH_i^{-1}\bfc_i = \bfx_0
+ Q_i\bfc_i,$$

\no dove, al solito $Q_i=(\bfq_1\dots\bfq_i)$. I rimanenti dettagli
implementativi sono simili a quelli visti per il QMR (e, precedentemente,
esaminati per MINRES); essi sono riassunti nell'Algoritmo~\ref{qmrfin}.

\begin{eser} Determinare il costo computazionale
dell'Algoritmo~\ref{qmrfin} e confrontarlo con quello
dell'Algoritmo~\ref{algo12}.\end{eser}

\begin{eser} Scrivere la versione precondizionata
dell'Algoritmo~\ref{qmrfin}.\end{eser}

\begin{eser} Scrivere un codice che implementi efficientemente
l'Algoritmo~\ref{qmrfin}, riformulandolo come un metodo iterativo
puro ed inserendo un controllo su possibili breakdown.\end{eser}

\begin{figure}[hp]
\begin{algo}\label{qmrfin} QMR, implementazione alternativa.\\ \rm
\fbox{\parbox{12cm}{
\begin{itemize}
\setlength{\itemsep}{0cm}
\nulit sia assegnata una approssimazione iniziale $\bfx_0$

\nulit calcolo $\bfr_0 = \bfb-A\bfx_0$~ e ~$\bb=\|\bfr_0\|$

\nulit inizializzo
$\bfu_1=\bfv_1=\bfq_1=\bfz_1=\bfr_0/\bb$,~$g_0=\bb$,~$f_1=1$

\nulit per i=1,\dots,n:

\begin{itemize}
\setlength{\itemsep}{0cm}
\nulit $\hat\bfq_i = A\bfq_i$

\nulit $\ell_i = \bfz_i^\top \hat\bfq_i/f_i$

\nulit se $\ell_i=0$, breakdown, esci, fine se

\nulit $\hat\bfu_i = \hat\bfq_i -\ell_i\bfu_i$

\nulit $\aa_i = \|\hat\bfu_i\|$

\nulit $r_{ii}=\ell_i$

\nulit $\bfp_i=\bfq_i$

\nulit se $i>1$, $\pmatrix{c}r_{i-1,i}\\r_{ii}\endpmatrix
\leftarrow
\pmatrix{rr} c_{i-1} &s_{i-1}\\ -s_{i-1} &c_{i-1}\endpmatrix
\pmatrix{c} 0\\ \ell_i \endpmatrix$

\nulit \qquad\qquad $\bfp_i \leftarrow \bfp_i
-r_{i-1,i}\bfp_{i-1}$, fine se

\nulit calcolo $c_i,s_i$ tali che $\pmatrix{c}r_{ii}\\0\endpmatrix
\leftarrow
\pmatrix{rr} c_i &s_i\\ -s_i &c_i\endpmatrix
\pmatrix{c} r_{ii}\\ \aa_i \endpmatrix$

\nulit $\bfp_i \leftarrow \bfp_i/r_{ii}$

\nulit $\pmatrix{c} \xi_i \\ g_i \endpmatrix =
\pmatrix{rr} c_i &s_i\\ -s_i &c_i\endpmatrix
\pmatrix{c} g_{i-1}\\ 0 \endpmatrix$

\nulit $\bfx_i = \bfx_{i-1} +\xi_i\bfp_i$

\nulit se $\aa_i=0$~~oppure~~$|g_i|\le tol$, $\bfx^*=\bfx_i$, esci, fine se

\nulit $\bfu_{i+1}=\hat\bfu_i/\aa_i$

\nulit $\bfv_{i+1} = ( A^\top \bfz_i -\ell_i\bfv_i )/\aa_i$

\nulit $f_{i+1} = \bfv_{i+1}^\top \bfu_{i+1}$

\nulit se $f_{i+1}=0$, breakdown, esci, fine se

\nulit $\phi_i = (\aa_i f_{i+1})/(\ell_i f_i)$

\nulit $\bfq_{i+1} = \bfu_{i+1} - \phi_i\bfq_i$

\nulit $\bfz_{i+1} = \bfv_{i+1} - \phi_i\bfz_i$

\end{itemize}

\nulit fine per

\end{itemize}

}}\end{algo}
\end{figure}

\begin{figure}[t]
\begin{center}
\includegraphics[width=13cm,height=10cm]{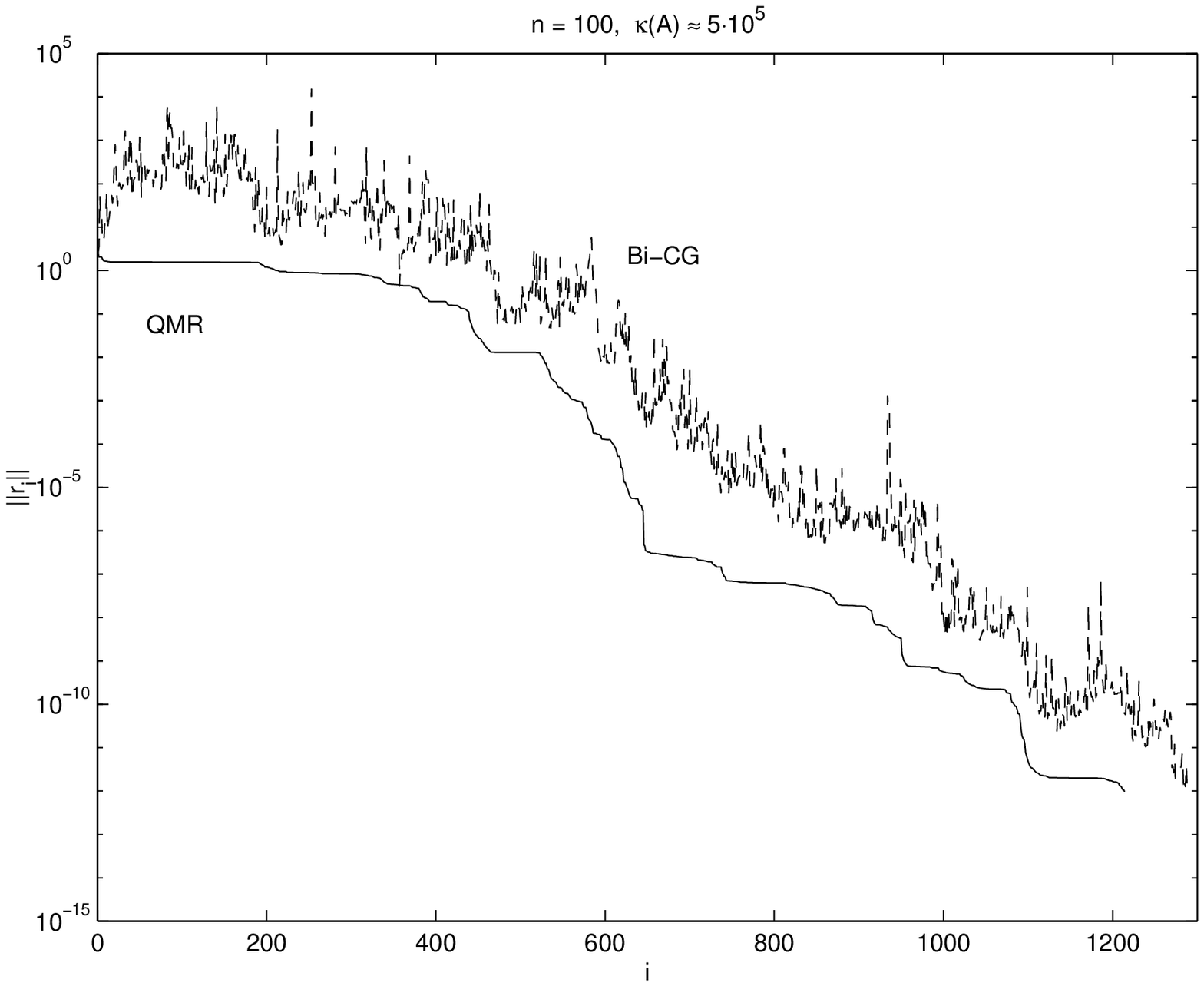}
\caption{\protect\label{fig13} QMR e Bi-CG applicati al problema
(\ref{Agen}).}\end{center}
\end{figure}

Per concludere, in Figura~\ref{fig13} si riporta la norma del residuo per QMR
e Bi-CG applicati allo stesso problema relativo alle
Figure~\ref{fig11}-\ref{fig12}. Dalla figura si evince come il comportamento del
primo metodo sia pi\`u regolare di quello del secondo, come ci si aspetta
dalla sua derivazione.

\section{Bidiagonalizzazione}\label{bidiag}

Una fattorizzazione alternativa, per ottenere un metodo di risoluzione
iterativo, \`e data da
\begin{equation}\label{bidiaf}
A V = U L, \qquad A^\top  U = V L^\top , \qquad U^\top U=V^\top V=I,
\end{equation}

\no dove $$U = (\bfu_1\dots\bfu_n), \qquad V = (\bfv_1\dots\bfv_n), \qquad L =
\pmatrix{cccc}
\aa_1\\
\bb_1 &\aa_2\\
      &\ddots&\ddots\\
      &      &\bb_{n-1} &\aa_n\endpmatrix.$$

\no Si tratta, quindi, di portare a forma bidiagonale la matrice $A$ mediante
moltiplicazione, a sinistra e a destra, per due matrici ortogonali.

Dalla (\ref{bidiaf}) si ottengono le seguenti relazioni di ricorrenza,
\begin{eqnarray}
\nonumber
A\bfv_i &=& \aa_i \bfu_i + \bb_i\bfu_{i+1},\\[1mm]
A^\top \bfu_i &=& \aa_i \bfv_i + \bb_{i-1}\bfv_{i-1},
 \qquad i=1,\dots,n,\label{ricor}
\end{eqnarray}

\no dove, al solito, $\bb_0=0$, $\bfu_{n+1}=\bfo$. Queste permettono di
ottenere il corrispondente Algoritmo~\ref{algo14}.

\begin{figure}[t]
\begin{algo}\label{algo14} Fattorizzazione (\protect\ref{bidiaf}).\\ \rm
\fbox{\parbox{12cm}{
\begin{itemize}
\setlength{\itemsep}{0cm}
\nulit sia assegnato un vettore non nullo $\hat\bfu_1$

\nulit inizializzo $\bb_0=0$, $\bfu_1 = \hat\bfu_1/\|\hat\bfu_1\|$

\nulit per $i = 1,\dots,n$:
\begin{itemize}
\setlength{\itemsep}{0cm}

\nulit $\hat\bfv_i = A^\top \bfu_i -\bb_{i-1} \bfv_{i-1}$

\nulit $\aa_i = \|\hat\bfv_i\|$

\nulit $\bfv_i = \hat\bfv_i/\aa_i$

\nulit se $i=n$, esci, fine se

\nulit $\hat\bfu_i = A\bfv_i -\aa_i\bfu_i$

\nulit $\bb_i = \|\hat\bfu_i\|$

\nulit se $\bb_i=0$, esci, fine se

\nulit $\bfu_{i+1} = \hat\bfu_i/\bb_i$

\end{itemize}

\nulit fine per

\end{itemize}
}}
\end{algo}
\end{figure}

Per utilizzare questa fattorizzazione al fine di risolvere il problema
(\ref{Agen}), assegnata una approssimazione iniziale $\bfx_0$ si
sceglie, nell'Algoritmo~\ref{algo14},

\begin{equation}\label{hu1}
\hat\bfu_1 = \bfr_0 \equiv \bfb-A\bfx_0.
\end{equation}

\no Quindi, al passo $i$-esimo, si ricercher\`a
~$\bfx_i \in \bfx_0 +[\bfv_1,\dots,\bfv_i]$~
in modo da soddisfare le seguenti condizioni di ortogonalit\`a per
il corrispondente residuo $\bfr_i$:
\begin{equation}\label{ortbid}
\bfu_j^\top \bfr_i = 0, \qquad j\le i,
\end{equation}

\no o, alternativamente, in modo da minimizzare $\|\bfr_i\|$. In entrambi i casi, si
ottiene una semplice relazione di ricorrenza per le approssimazioni $\{\bfx_i\}$.

\begin{eser} Verificare che, se si utilizza il criterio (\ref{ortbid}), allora
$$\bfx_i = \bfx_{i-1} + \xi_i\bfv_i,\qquad i = 1,2,\dots,$$ dove $$\xi_1 =
\|\bfr_0\|/\aa_1, \qquad \xi_i = -\xi_{i-1}\bb_{i-1}/\aa_i,\quad i=2,\dots,n.$$
\end{eser}

\no Inoltre, in virt\`u della scelta (\ref{hu1}), si dimostra valere il
seguente risultato.

\begin{teo}\label{bi0} Se, nell'Algoritmo~\ref{algo14}-(\ref{hu1}), si verifica
$\bb_i=0$, allora $\bfx_i=\bfx^*$, soluzione di (\ref{Agen}).\end{teo}

\medskip
\begin{eser} Dimostrare il Teorema~\ref{bi0}.\end{eser}

\medskip
\begin{eser} Dimostrare che, per l'Algoritmo~\ref{algo14},
$\aa_1>0$. Inoltre, per $i>1$, se $\bb_{i-1}\ne0$, allora
$\aa_i>0$ (suggerimento: utilizzare la fattorizzazione (\ref{bidiaf})).\end{eser}

Purtroppo, la fattorizzazione (\ref{bidiaf}) che, apparentemente non ha pi\`u
il problema del breakdown tipico del metodo di Lanczos non simmetrico, ha un
serio inconveniente. Infatti, moltiplicando la prima equazione in
(\ref{bidiaf}) per $A^\top $, e tenendo conto della seconda equazione, si ottiene:
$$A^\top A V = A^\top UL = VL^\top L \equiv V T, \qquad V^\top V=I,$$

\no dove $T\equiv L^\top L$ \`e una matrice tridiagonale e simmetrica. Si evince,
quindi, che la fattorizzazione (\ref{bidiaf}) \`e equivalente alla
fattorizzazione di Lanczos (\ref{fattlanc}) applicata alla matrice $A^\top A$.
In definitiva, il corrispondente metodo iterativo equivale a risolvere le
equazioni normali (\ref{normaleq}) associate al problema (\ref{Agen}), peggiorando
grandemente il condizionamento del problema. Pertanto, il suo utilizzo ha scarso
interesse pratico.

\section{CGS}\label{cgs}

I metodi derivati dalla fattorizzazione
(\ref{lancuns}), ovvero Bi-CG e QMR, richiedono, ad ogni
iterazione, due matvec: uno con la matrice $A$, l'altro con la
matrice $A^\top $ (vedi gli Algoritmi~\ref{algo11} e \ref{qmrfin}).
Quest'ultimo prodotto, quando $A$ \`e una matrice
sparsa di grandi dimensioni, pu\`o talora risultare difficoltoso,
se si \`e utilizzata una tecnica di memorizzazione per matrici
sparse. Per questo motivo, cercheremo ora di derivare una variante
dei Bi-CG che non necessita di matvec con $A^\top $.

Il punto di partenza della nostra analisi sar\`a
l'Algoritmo~\ref{algo11} del metodo dei Bi-CG. Infatti, se si
considerano le direzioni $\{\bfp_i\}$, $\{\hat\bfp_i\}$, i
residui $\{\bfr_i\}$ e gli {\em pseudo-residui} $\{\hat\bfr_i\}$,
si verifica facilmente che
\begin{eqnarray}\nonumber
\bfr_i &=& \phi_i(A)\bfr_0, \qquad \hat\bfr_i =
\phi_i(A^\top )\hat\bfr_0,
\\[-2mm] \label{piri}\\[-2mm] \nonumber
\bfp_i&=& \psi_i(A)\bfr_0, \qquad \hat\bfp_i =
\psi_i(A^\top )\hat\bfr_0,
\end{eqnarray}

\no dove $\phi_i,\psi_{i+1}\in\Pi_i$ sono polinomi definiti
ricorsivamente come segue (sia $\theta(x)=x$):
\begin{eqnarray}\label{fipsi}
\phi_0 &\equiv& 1, \qquad \phi_i = \phi_{i-1} - \hat\lam_i \theta
\psi_i,\\[1mm] \nonumber
\psi_1 &\equiv& 1, \qquad \psi_{i+1} = \phi_i +\mu_i\psi_i, \qquad
i = 1,2,\dots.\end{eqnarray}

\no Gli scalari $\{\hat\lam_i\}, \{\mu_i\}$ che definiscono tali
polinomi sono
\begin{equation}\label{mulam}
\mu_i = \eta_i/\eta_{i-1}, \qquad \hat\lam_i =
\eta_{i-1}/\hat{d}_i,\end{equation}

\no dove:
\begin{eqnarray}\nonumber
\eta_i &=& \hat\bfr_i^\top \bfr_i =
\left(\phi_i(A^\top )\hat\bfr_0\right)^\top \phi_i(A)\bfr_0 =
\hat\bfr_0^\top \phi_i^2(A)\bfr_0,\\[-2mm] \label{polsquare} \\[-2mm] \nonumber
\hat{d}_i &=& \hat\bfp_i^\top  A \bfp_i = \left(
\psi_i(A^\top )\hat\bfr_0\right)^\top  A\psi_i(A)\bfr_0 =
\hat\bfr_0^\top A\psi_i^2(A)\bfr_0.
\end{eqnarray}

\no Pertanto, potremmo ottenere un metodo iterativo assimilabile
ai Bi-CG, che consiste nel tener fisso il primo pseudo-residuo, $\hat\bfr_0$,
e che non richiede l'utilizzo di $A^\top $, a patto che si possano generare ricorsivamente i
polinomi $\{\psi_i^2\}$ e $\{\phi_i^2\}$. Dalla (\ref{fipsi}) si ottiene:
\begin{eqnarray}\nonumber
\phi_0^2&\equiv& 1, \qquad \psi_1^2 \equiv 1, \\[1mm] \label{fipsi2}
\phi_i^2 &=& \phi_{i-1}^2 +\hat\lam_i^2\theta^2\psi_i^2 -2\hat\lam_i\theta\phi_{i-1}\psi_i,
\\[1mm] \nonumber
\psi_{i+1}^2 &=& \phi_i^2 +\mu_i^2\psi_i^2 +2\mu_i\phi_i\psi_i,\qquad i = 1,2,\dots.
\end{eqnarray}

\no Per quanto riguarda i prodotti $\{\phi_{i-1}\psi_i\}$ e
$\{\phi_i\psi_i\}$, si ottengono le relazioni di ricorrenza:
\begin{eqnarray}\nonumber
\phi_0\psi_1 &\equiv&1, \\ \label{fimix}
\phi_i\psi_i &=& \left(\phi_{i-1}-\hat\lam_i\theta\psi_i\right)\psi_i =
\phi_{i-1}\psi_i -\hat\lam_i\theta\psi_i^2,\\ \nonumber
\phi_i\psi_{i+1} &=& \phi_i\left(\phi_i+\mu_i\psi_i\right) =
\phi_i^2 +\mu_i\phi_i\psi_i, \qquad i=1,2,\dots.
\end{eqnarray}

\no Definendo i polinomi
$$\Phi_i = \phi_i^2, \qquad \Psi_i = \psi_i^2, \qquad \Omega_i =
\phi_i\psi_i, \qquad \Gamma_i = \phi_i\psi_{i+1},$$

\no e la forma bilineare
\begin{equation}\label{bilin}
[p,q] \equiv \left(p(A^\top )\hat\bfr_0\right)^\top  q(A)\bfr_0 \equiv [1,pq]\equiv[1,qp],
\end{equation}

\no possiamo quindi formulare l'Algoritmo~\ref{cgsform} per la
generazione degli stessi.

\begin{figure}[hp]
\begin{algo}\label{cgsform} Procedura per la generazione dei
polinomi (\ref{fipsi2})-(\ref{fimix}).\\ \rm
\fbox{\parbox{12cm}{
\begin{itemize}
\setlength{\itemsep}{0cm}
\nulit Inizializzo $\Phi_0\equiv\Psi_1\equiv\Gamma_0\equiv1$,
$\eta_0=[1,\Phi_0]$

\nulit per $i = 1,\dots,n$:
\begin{itemize}
\setlength{\itemsep}{0cm}

\nulit $\hat{d}_i = [1,\theta\Psi_i]$

\nulit $\hat\lam_i = \eta_{i-1}/\hat{d}_i$

\nulit $\Omega_i = \Gamma_{i-1} -\hat\lam_i\theta\Psi_i$

\nulit $\Phi_i = \Phi_{i-1}
-\hat\lam_i\theta(\Gamma_{i-1}+\Omega_i)$

\nulit $\eta_i = [1,\Phi_i]$

\nulit $\mu_i = \eta_i/\eta_{i-1}$

\nulit $\Gamma_i = \Phi_i +\mu_i\Omega_i$

\nulit $\Psi_{i+1} = \Gamma_i +\mu_i(\Omega_i+\mu_i\Psi_i)$

\end{itemize}

\nulit fine per

\end{itemize}
}}
\end{algo}

\bigskip
\begin{algo}\label{cgsalg} Metodo CGS.\\ \rm
\fbox{\parbox{12cm}{
\begin{itemize}
\setlength{\itemsep}{0cm}

\nulit sia assegnata una approssimazione iniziale $\bfx_0$

\nulit calcolo $\bfr_0 = \bfb-A\bfx_0$

\nulit inizializzo $\hat\bfr_0 = \bfp_1 = \bfg_0 = \bfr_0$,
$\eta_0=\hat\bfr_0^\top \bfr_0$

\nulit per $i = 1,\dots,n$:
\begin{itemize}
\setlength{\itemsep}{0cm}

\nulit $\bfv_i = A\bfp_i$

\nulit $\hat{d}_i = \hat\bfr_0^\top \bfv_i$

\nulit $\hat\lam_i = \eta_{i-1}/\hat{d}_i$

\nulit $\bfw_i = \bfg_{i-1} -\hat\lam_i\bfv_i$

\nulit $\bfx_i = \bfx_{i-1} +\hat\lam_i(\bfg_{i-1}+\bfw_i)$

\nulit $\bfr_i = \bfr_{i-1}
-\hat\lam_i A (\bfg_{i-1}+\bfw_i)$

\nulit $\eta_i = \hat\bfr_0^\top \bfr_i$

\nulit $\mu_i = \eta_i/\eta_{i-1}$

\nulit $\bfg_i = \bfr_i +\mu_i\bfw_i$

\nulit $\bfp_{i+1} = \bfg_i +\mu_i(\bfw_i+\mu_i\bfp_i)$

\end{itemize}

\nulit fine per

\end{itemize}
}}
\end{algo}

\end{figure}

Tenendo conto del fatto che i residui generati sono ora dati da
$$\bfr_i = \Phi_i(A)\bfr_0 = \phi_i^2(A)\bfr_0,$$

\no si ottiene che essi sono i ``quadrati'' dei corrispondenti residui
generati dai Bi-CG (vedi (\ref{piri})). Per questo motivo, il
corrispondente metodo iterativo \`e denominato {\em CGS (CG
Squared)}: esso \`e descritto dall'Algoritmo~\ref{cgsalg}, ottenuto
dal precedente Algoritmo~\ref{cgsform}. Al solito, \`e stata utilizzata
la scelta
\begin{equation}\label{ro0tro0}\hat\bfr_0=\bfr_0=\bfb-A\bfx_0.\end{equation}

\no La versione precondizionata del metodo si ottiene mediante le
sostituzioni formali $$ \bfb \leftarrow C\bfb, \qquad A \leftarrow
CA,$$

\no dove $C\approx A^{-1}$ \`e la matrice di precondizionamento.

\medskip
\begin{oss} Dall'Algoritmo~\ref{cgsalg} si evince che non sono pi\`u
ri\-chiesti, ad ogni iterata, il matvec per $A$ e quello per $A^\top $ ma solo
due matvec con $A$. \end{oss}

\medskip
\begin{eser} Determinare il costo computazionale del metodo CGS.
\end{eser}

\medskip
\begin{eser} Riscrivere in forma pi\`u efficiente
l'Algoritmo~\ref{cgsalg}, implementando un metodo iterativo puro e
controllando eventuali breakdown.\end{eser}

\medskip
\begin{eser} Scrivere una procedura che implementi efficientemente
l'Algoritmo~\ref{cgsalg}.\end{eser}

\begin{figure}[hp]
\begin{center}
\includegraphics[width=13cm,height=8.5cm]{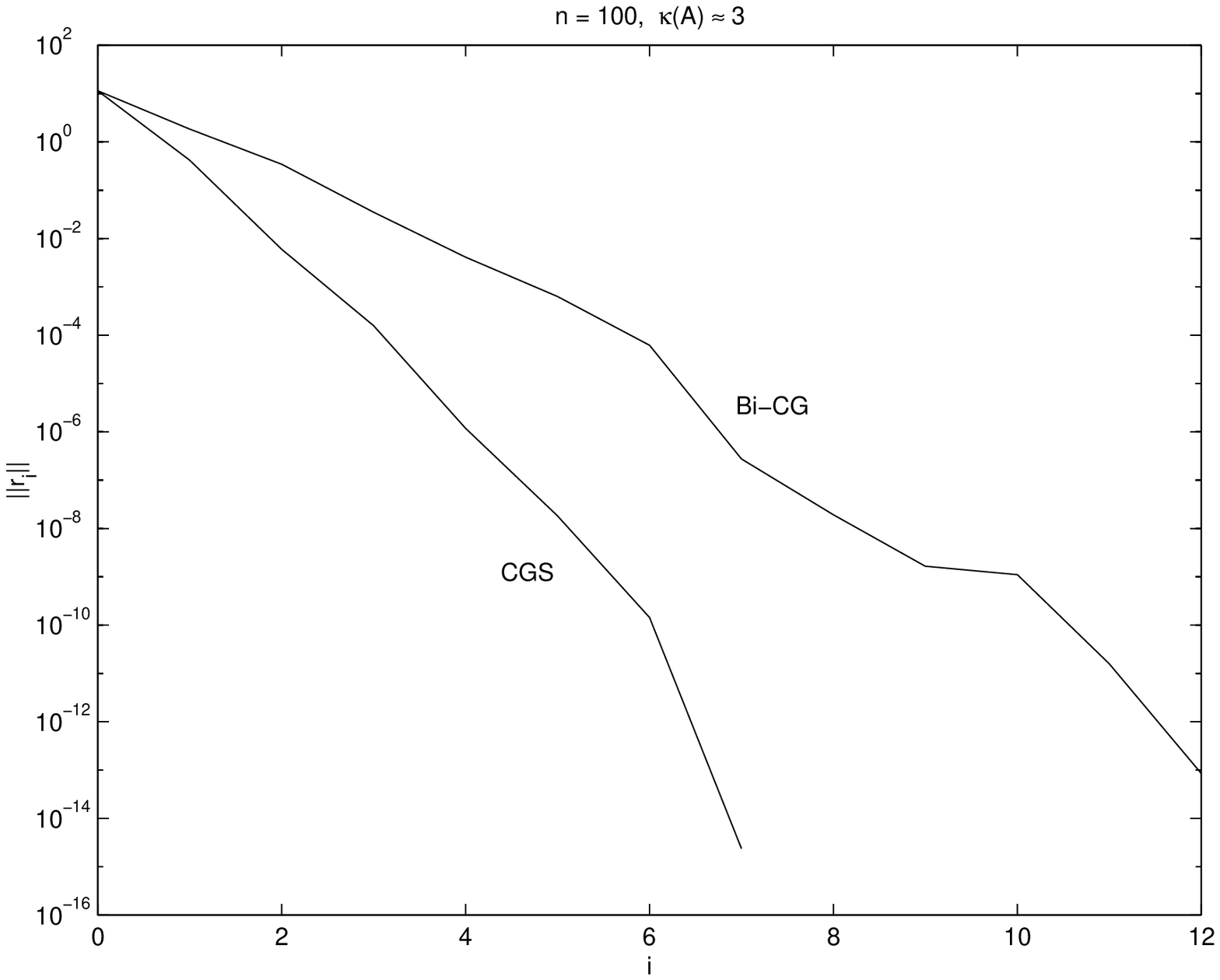}
\caption{\protect\label{fig14} CGS e Bi-CG applicati ad un problema
ben condizionato.}\end{center}
\begin{center}
\includegraphics[width=13cm,height=8.5cm]{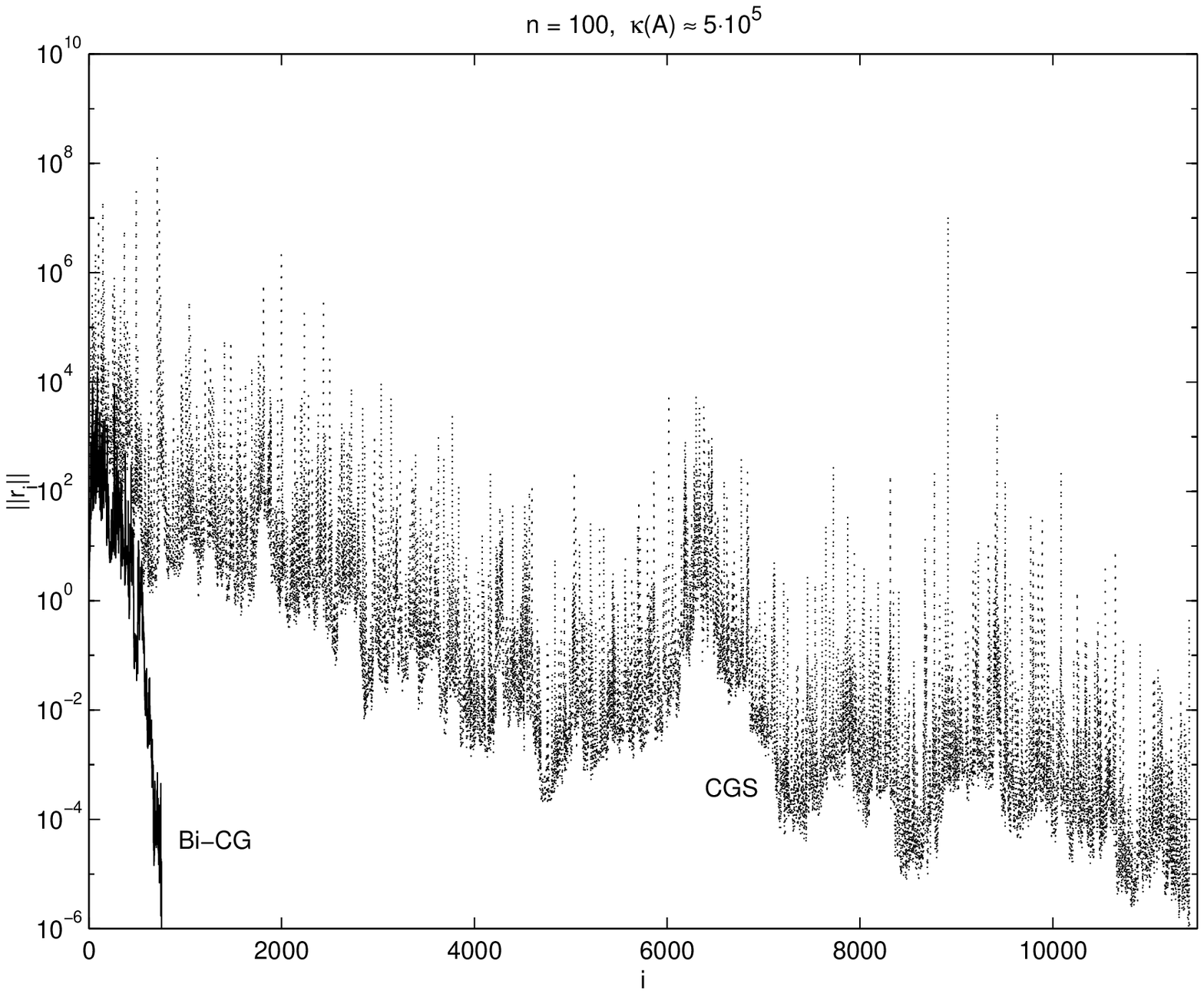}
\caption{\protect\label{fig15} CGS e Bi-CG applicati ad un problema
mal condizionato.}
\end{center}
\end{figure}

\medskip
Se il metodo Bi-CG converge monotonicamente, allora questo significa
che (vedi (\ref{piri})) $\|\phi_i(A)\bfr_0\|$ \`e una contrazione,
al crescere di $i$. In tal caso, \`e da aspettarsi che CGS, per il
quale i residui sono dati da $\phi_i^2(A)\bfr_0$, converga ancora
pi\`u rapidamente. Questo, generalmente, avviene nel caso di problemi
ben condizionati (vedi Figura~\ref{fig14}). Tuttavia, abbiamo visto
(vedi ad esempio Figura~\ref{fig12}) che talora il residuo del Bi-CG pu\`o
oscillare grandemente. In tal caso, CGS osciller\`a ancora pi\`u
violentemente, e questo generalmente comporta seri problemi,
quando si utilizza l'aritmetica finita. Ad esempio, in Figura~\ref{fig15}
si riporta la norma del residuo per CGS e Bi-CG applicati allo
stesso problema relativo alle Figure~\ref{fig11}-\ref{fig12}, da cui
si evince come il comportamento del primo metodo sia pi\`u irregolare
di quello del secondo, proprio in virt\`u delle oscillazioni
di quest'ultimo (Bi-CG). In questo caso, il Bi-CG
im\-piega 755 iterazioni per ridurre la norma del residuo al di
sotto di $10^{-6}$, mentre il CGS richiede poco meno di
11500 iterazioni per ottenere lo stesso risultato. Per ovviare a
questo problema, nella prossima sezione verr\`a esa\-minata una ulteriore
variante del metodo Bi-CG, anche questa non richiedente il matvec per $A^\top $.

\section{Bi-CGstab}\label{bicgstab}

Per evitare l'inconveniente di avere residui con una norma
elevata, come abbiamo visto pu\`o avvenire nel caso del CGS,
invece di ricercare i residui e le direzioni nella forma
~$\phi_i^2(A)\bfr_0$~ e ~$\psi_i^2(A)\bfr_0,$~ rispettivamente,
dove i polinomi $\{\phi_i\}$ e $\{\psi_i\}$ sono quelli
associati ai Bi-CG (vedi (\ref{piri})-(\ref{fipsi})), li ricercheremo
nella forma
\begin{equation}\label{ripi1}
\bfr_i = \chi_i(A)\phi_i(A)\bfr_0, \qquad \bfp_i =
\chi_{i-1}(A)\psi_i(A)\bfr_0, \qquad i=1,2,\dots,\end{equation}

\no dove
\begin{equation}\label{chi}
\chi_i(z) = (1-\omega_i z)\chi_{i-1}(z), \qquad i=1,2,\dots,\qquad \chi_0\equiv1,
\end{equation}

\no ed il parametro $\omega_i$ \`e scelto, al passo $i$-esimo, in
modo da minimizzare $\|\bfr_i\|$. Supporremo, nel seguito, che $\omega_j\ne0$, per ogni $j\ge1$:
in tal caso, $\chi_i(z)$ ha grado esatto $i$.
Il problema, al solito, \`e quello di generare i polinomi $\{\phi_i\}$ e $\{\psi_i\}$,
ovvero i parametri (\ref{mulam})-(\ref{polsquare}).

Tenendo in conto di quanto all'Esercizio~\ref{bicges}, si ottiene,
tuttavia, che, essendo $\phi_i(A)\bfr_0$ ortogonale al sottospazio
$$[\hat\bfr_0,A^\top \hat\bfr_0,\dots,(A^\top )^{i-1}\hat\bfr_0]\equiv
[\hat\bfr_0,\chi_1(A^\top )\hat\bfr_0,\dots,\chi_{i-1}(A^\top )\hat\bfr_0],$$
allora (vedi (\ref{bilin})):
$$\eta_i \equiv [\phi_i,\phi_i] = b_i^{(i)}\,\hat\bfr_0^\top A^i\phi_i(A)\bfr_0,$$

\no dove il coefficiente $b_i^{(i)}$ \`e dato da (\ref{maincoef}).
Tuttavia, noi possiamo in realt\`a calcolare
$$\tilde\eta_i \equiv [\chi_i,\phi_i] =
\tilde{b}_i^{(i)}\,\hat\bfr_0^\top A^i\phi_i(A)\bfr_0,$$

\no dove
\begin{equation}\label{maincoef1}
\tilde{b}_i^{(i)} = (-1)^i\prod_{j=1}^i \omega_j
\end{equation}

\no \`e il coefficiente principale di $\chi_i(z)$ (vedi (\ref{chi})).
Dalle precedenti equazioni, si ricava quindi che
$$\mu_i = \frac{\eta_i}{\eta_{i-1}} =
\frac{\tilde\eta_i}{\tilde\eta_{i-1}}\,\frac{\hat\lam_i}{\omega_i}.$$

In modo del tutto analogo, sempre dall'Esercizio~\ref{bicges}, si
ha che, essendo $\psi_i(A)\bfr_0$ $A$-ortogonale a
$$[\hat\bfr_0,A^\top \hat\bfr_0,\dots,(A^\top )^{i-2}\hat\bfr_0],$$

\no allora (vedi (\ref{maincoef}))
$$\hat{d_i} \equiv [\psi_i,\theta\psi_i] = b_{i-1}^{(i-1)}\,\hat\bfr_0^\top A^i\psi_i(A)\bfr_0,$$

\no mentre noi possiamo calcolare (vedi (\ref{maincoef1}))
$$\tilde{d}_i \equiv [\chi_{i-1},\theta\psi_i] =
\tilde{b}_{i-1}^{(i-1)}\, \hat\bfr_0^\top A^i \psi_i(A)\bfr_0.$$

\no Ne consegue che
$$\hat\lam_i = \frac{\eta_{i-1}}{\hat{d}_i} =
\frac{\tilde\eta_{i-1}}{\tilde{d}_i}.$$

\no Rimane da vedere la scelta del parametro $\omega_i$ (vedi
(\ref{chi})), che serve a minimizzare la norma di $\bfr_i$ in
(\ref{ripi1}). Si ottiene:
\begin{eqnarray*}
\bfr_i^\top \bfr_i &=&
(\chi_{i-1}(A)\phi_i(A)\bfr_0)^\top (I-\omega_iA)^\top (I-\omega_iA)(\chi_{i-1}(A)\phi_i(A)\bfr_0)\\[1mm]
&\equiv& \bfr_{i-\frac{1}2}^\top (I +\omega_i^2A^\top A
-2\omega_iA)\bfr_{i-\frac{1}2},\end{eqnarray*}

\no che \`e minimizzata dalla scelta
$$\omega_i = \frac{\bfr_{i-\frac{1}2}^\top \bft_i}{\bft_i^\top \bft_i}, \qquad \bft_i
= A\bfr_{i-\frac{1}2}.$$

\no Osserviamo che la precedente quantit\`a \`e definita, a meno
che il vettore
$$\bfr_{i-\frac{1}2} = \chi_{i-1}(A)\phi_i(A)\bfr_0 = \chi_{i-1}(A)(\phi_{i-1}(A)-\hat\lam_iA\psi_i(A))\bfr_0
=\bfr_{i-1}-\hat\lam_i A\bfp_i$$

\no non si annulli. Ma questo altri non \`e che il residuo del corrispondente passo
di Bi-CG, premoltiplicato per $\chi_{i-1}(A)$. Pertanto, se esso \`e nullo, vuol dire che abbiamo
raggiunto la soluzione. Proprio per la caratteristica di rendere
pi\`u {\em smooth} il residuo, il metodo appena definito \`e stato
denominato {\em Bi-CGStab (Bi-CG Stabilized)}. Esso \`e descritto
dall'Algoritmo~\ref{BiCGS}, in cui, al solito, si \`e utilizzata
la scelta (\ref{ro0tro0}) per $\hat\bfr_0$.

\begin{figure}[t]
\begin{algo}\label{BiCGS} Bi-CGStab.\\ \rm
\fbox{\parbox{12cm}{
\begin{itemize}
\setlength{\itemsep}{0cm}

\nulit sia assegnata una approssimazione iniziale $\bfx_0$

\nulit calcolo $\bfr_0 = \bfb-A\bfx_0$

\nulit inizializzo $\hat\bfr_0 = \bfp_1 = \bfr_0$,
$\eta_0=\hat\bfr_0^\top \bfr_0$

\nulit per $i = 1,\dots,n$:
\begin{itemize}
\setlength{\itemsep}{0cm}

\nulit $\bfv_i = A\bfp_i$

\nulit $\hat{d}_i = \hat\bfr_0^\top \bfv_i$

\nulit $\hat\lam_i = \eta_{i-1}/\hat{d}_i$

\nulit $\bfx_{i-\frac{1}2} = \bfx_{i-1} +\hat\lam_i \bfp_i$

\nulit $\bfr_{i-\frac{1}2} = \bfr_{i-1} -\hat\lam_i\bfv_i$

\nulit $\bft_i = A\bfr_{i-\frac{1}2}$

\nulit $\omega_i = \bft_i^\top \bfr_{i-\frac{1}2}/\bft_i^\top \bft_i$

\nulit $\bfx_i = \bfx_{i-\frac{1}2}+\omega_i\bfr_{i-\frac{1}2}$

\nulit $\bfr_i = \bfr_{i-\frac{1}2} -\omega_i \bft_i$

\nulit $\eta_i = \hat\bfr_0^\top \bfr_i$

\nulit $\mu_i = (\eta_i/\eta_{i-1})(\hat\lam_i/\omega_i)$

\nulit $\bfp_{i+1} = \bfr_i +\mu_i(\bfp_i -\omega_i \bfv_i)$

\end{itemize}

\nulit fine per

\end{itemize}
}}
\end{algo}

\end{figure}

\begin{oss} Il modo pi\`u semplice per ottenere una versione precondizionata
del Bi-CGStab \`e, al solito, quello di effettuare le sostituzioni formali
$$A\leftarrow CA, \qquad \bfb\leftarrow C\bfb,$$ dove $C$ \`e la matrice di
precondizionamento.\end{oss}

\begin{oss} Un criterio di arresto basato sulla norma del residuo pu\`o
controllare, alla generica iterata $i$-esima, sia $\|\bfr_i\|$ che
$\|\bfr_{i-\frac{1}2}\|$.\end{oss}

\begin{eser} Determinare il costo computazionale dell'Algoritmo~\ref{BiCGS}.
\end{eser}

\begin{eser} Scrivere un codice che implementi efficientemente
l'Algoritmo~\ref{BiCGS}, riformulandolo come un metodo iterativo puro.
\end{eser}

\begin{figure}[t]
\begin{center}
\includegraphics[width=13cm,height=10cm]{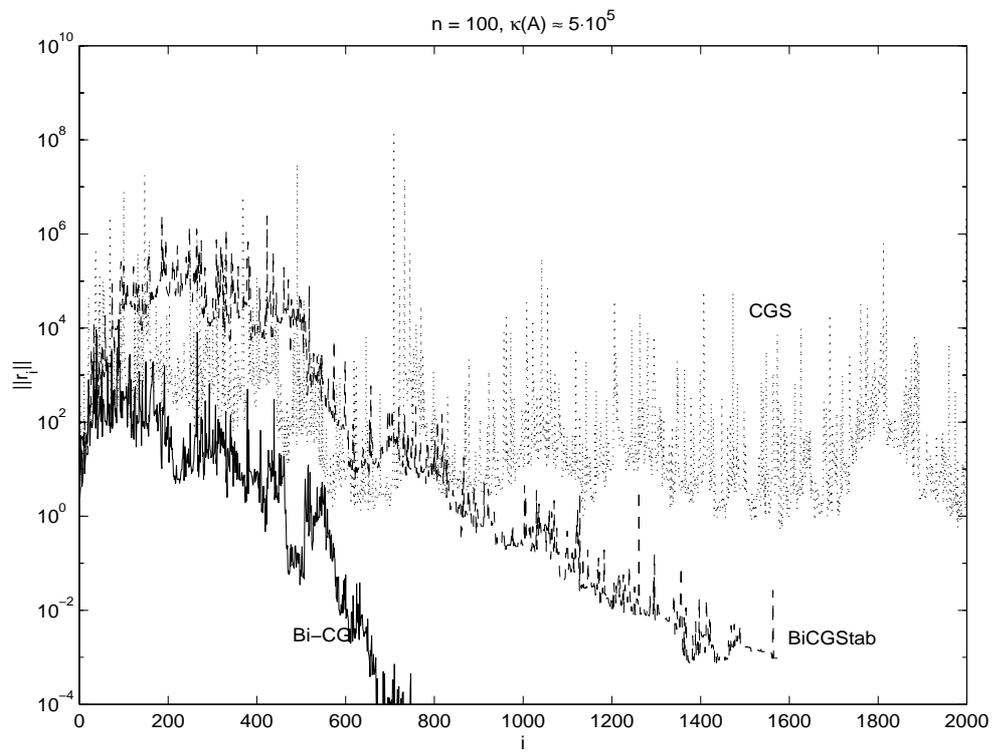}
\caption{\protect\label{fig16} CGS, Bi-CGStab e Bi-CG applicati al problema di Figura~\ref{fig15}.}
\end{center}
\end{figure}

Per concludere, consideriamo l'applicazione del Bi-CGStab allo stesso problema
di Figura~\ref{fig15}. In tal caso, il residuo del Bi-CGStab ``stagna'' a circa 
1600 iterazioni (vedi Figura~\ref{fig16}), sebbene il suo comportamento risulti essere pi\`u favorevole di quello di CGS.
Generalmente, la performance di questo metodo migliora grandemente con l'utilizzo di un idoneo
precondizionatore.

%
%
\appendix
\chapter{Memorizzazione di matrici sparse}\label{appA}

{\em
In questa appendice si descrivono brevemente alcune delle pi\`u
comuni tecni\-che di memorizzazione per matrici sparse, con il
corrispondente algoritmo per effettuare il matvec. Si sottolinea,
comunque, che spesso la memorizzazione di matrici con struttura
pu\`o richiedere tecniche ibride o {\em ad hoc}. Inoltre,
sostanziosi miglioramenti della performance possono derivare
dall'utilizzo efficiente del linguaggio di programmazione (ad
esempio, memorizzare le matrici nella forma pi\`u conveniente,
rispetto alla allocazione interna utilizzata).}

\section{Memorizzazione compressa per righe}

Questo tipo di memorizzazione \`e indicato per matrici sparse che
hanno un (basso) numero di elementi non nulli su ciascuna riga.
Sia, dunque:

\begin{description}
\item{$k = $} massimo numero di elementi non nulli su ciascuna riga
della matrice, sia essa $A$, di dimensione $n\times n$.
\end{description}

\no Gli elementi della matrice saranno quindi memorizzati
utilizzando due array:

\begin{description}

\item{\tt A(1..n,1..k)}: {\tt A(i,j)} contiene il $j$-esimo elemento
significativo (ovvero non nullo) sulla riga $i$-esima di $A$. Gli
elementi non utilizzati sono posti uguali a 0;

\item{\tt COL(1..n,1..k)}: {\tt COL(i,j)} contiene l'indice di colonna del
$j$-esimo elemento si\-gnificativo sulla riga $i$-esima di $A$.

\end{description}

\begin{oss} Pertanto, se $A=(a_{ij})$, allora {\tt A(i,j)} conterr\`a
l'elemento $a_{i,COL(i,j)}$. \end{oss}

\begin{oss} Al fine di minimizzare gli spostamenti da memoria, \`e
conveniente utilizzare la seguente strategia: se $${\tt
A(i,\ell)}\ne0\qquad e \qquad {\tt A(i,j)}=0,\quad {\tt
j=\ell+1,\dots,k},$$ allora si pone $${\tt COL(i,j)} = {\tt
COL(i,\ell)}, \qquad {\tt j=\ell+1,\dots,k}.$$
\end{oss}

\no Di seguito, si riporta l'algoritmo per effettuare il
corrispondente matvec.

\begin{algo}\label{A1} Matvec $\bfy=A\bfx$: memorizzazione compressa per righe.\\ \rm
\fbox{\parbox{12cm}{

\begin{itemize}

\nulit per i = 1,\dots,n:

\begin{itemize}
\setlength{\itemsep}{0cm}

\nulit y(i) = 0

\nulit per j = 1,\dots,k:

\begin{itemize}
\setlength{\itemsep}{0cm}

\nulit y(i) = y(i) + A(i,j)\,$*$\,x(COL(i,j))

\end{itemize}

\nulit fine per

\end{itemize}

\nulit fine per

\end{itemize}

}}\end{algo}

\section{Memorizzazione compressa per colonne}

Questo tipo di memorizzazione \`e, in un certo senso,
complementare a quello precedentemente esaminato. Esso \`e,
infatti, indicato per matrici sparse che hanno un (basso) numero
di elementi non nulli su ciascuna colonna. Sia, dunque:

\begin{description}
\item{$k = $} massimo numero di elementi non nulli su ciascuna
colonna della matrice, sia essa $A$, di dimensione $n\times n$.
\end{description}

\no Gli elementi della matrice saranno quindi memorizzati
utilizzando due array:

\begin{description}

\item{\tt A(1..k,1..n)}: {\tt A(i,j)} contiene l'$i$-esimo elemento
significativo sulla colonna $j$-esima di $A$. Gli elementi non
utilizzati sono posti uguali a 0;

\item{\tt ROW(1..k,1..n)}: {\tt ROW(i,j)} contiene l'indice di riga
dell'$i$-esimo elemento si\-gnificativo sulla colonna $j$-esima di
$A$.

\end{description}

\begin{oss} Pertanto, se $A=(a_{ij})$, allora {\tt A(i,j)} conterr\`a
l'elemento $a_{ROW(i,j),j}$. \end{oss}

\begin{oss} Al fine di minimizzare gli spostamenti da memoria, \`e
conveniente utilizzare la seguente strategia: se $${\tt
A(\ell,j)}\ne0\qquad e \qquad {\tt A(i,j)}=0,\quad {\tt
i=\ell+1,\dots,k},$$ allora si pone $${\tt ROW(i,j)} = {\tt
ROW(\ell,j)}, \qquad {\tt i=\ell+1,\dots,k}.$$
\end{oss}

\no Di seguito, si riporta l'algoritmo per effettuare il
corrispondente matvec.

\begin{algo} Matvec $\bfy=A\bfx$: memorizzazione compressa per colonne.\\ \rm
\fbox{\parbox{12cm}{

\begin{itemize}

\nulit per i = 1,\dots,n:

\begin{itemize}
\setlength{\itemsep}{0cm}

\nulit y(i) = 0

\end{itemize}

\nulit fine per

\nulit per j = 1,\dots,n:

\begin{itemize}
\setlength{\itemsep}{0cm}

\nulit per i = 1,\dots,k:

\begin{itemize}
\setlength{\itemsep}{0cm}

\nulit y(ROW(i,j)) = y(ROW(i,j)) + A(i,j)\,$*$\,x(j)

\end{itemize}

\nulit fine per

\end{itemize}

\nulit fine per

\end{itemize}

}}\end{algo}

\section{Memorizzazione compressa per diagonali}

Un caso abbastanza rilevante nelle applicazioni \`e quello in cui
la matrice $A$ ha struttura di sparsit\`a {\em diagonale}, ovvero
i suoi elementi significativi sono situati lungo (poche) diagonali
della matrice stessa. Gli esempi dei problemi esaminati nel
Capitolo~\ref{cap1} sono, per l'appunto, di questo tipo.

Osserviamo che, se $A=(a_{ij})$, allora consideriamo la sua
$\nu$-esima diagonale, con la convenzione che:
\begin{itemize}

\item a $\nu=0$ corrisponde la diagonale principale;

\item a $\nu>0$ corrisponde la $\nu$-esima sopradiagonale;

\item a $\nu<0$ corrisponde la ($-\nu$)-esima sottodiagonale.

\end{itemize}

\no Osserviamo che $\nu$ coincide con la differenza $j-i$ degli
indici degli elementi $a_{ij}$ lungo la corrispondente diagonale.
Pertanto, se $A\in\RR^{n\times n}$ ha gli elementi significativi
di\-sposti lungo $k$ (generalmente, $k\ll n$) diagonali, potremo
memorizzare questi elementi nei due seguenti array:

\begin{description}

\item{\tt A(1..n,1..k)}: {\tt A(i,r)} contiene l'elemento sulla riga $i$-esima
lungo la dia\-gonale memorizzata in colonna $r$;

\item{\tt NU(1..k)}: {\tt NU(r)} contiene l'indice $\nu\equiv j-i$ relativo
alla diagonale memorizzata in {\tt A(1..n,r)}.

\end{description}

\no Di seguito si riporta l'algoritmo per effettuare il
corrispondente matvec.

\begin{algo} Matvec $\bfy=A\bfx$: memorizzazione compressa per diagonali.\\ \rm
\fbox{\parbox{12cm}{

\begin{itemize}

\nulit per i = 1,\dots,n:

\begin{itemize}
\setlength{\itemsep}{0cm}

\nulit y(i) = 0

\end{itemize}

\nulit fine per

\nulit per r = 1,\dots,k:

\begin{itemize}
\setlength{\itemsep}{0cm}

\nulit start = max(1,1-NU(r))

\nulit end    \,~= min(n,n-NU(r))

\nulit per i = start,\dots,end:

\begin{itemize}
\setlength{\itemsep}{0cm}

\nulit y(i) = y(i) + A(i,r)\,$*$\,x(i+NU(r))

\end{itemize}

\nulit fine per

\end{itemize}

\nulit fine per

\end{itemize}

}}\end{algo}

\medskip
\begin{oss} Laddove sia convenientemente applicabile, la
memorizzazione compressa per diagonali risulta essere quella pi\`u
conveniente, dal punto di vista del costo computazionale, tra
quelle esaminate.\end{oss}

%
%

\chapter*{~}
\thispagestyle{empty}
~
\newpage

\thispagestyle{empty}
\parbox{12cm}{\small

\rule{12cm}{1mm}
\bigskip

Il presente volume ha l'obiettivo di fornire le nozioni principali
riguardo ai metodi iterativi per la risoluzione di sistemi
di equazioni lineari sparsi di grandi dimensioni.

La trattazione, sebbene non esaustiva, ha comunque l'obiettivo di
fornire, in un contesto il pi\`u possibile omogeneo ed
autoconsistente, le idee principali che hanno ispirato le tecniche
di risoluzione pi\`u note.

Il materiale presentato \`e sufficiente a coprire i contenuti di
un modulo di 40-45 ore di lezione. Inoltre, gli esercizi proposti
possono essere argomento di ulteriori ore di esercitazione. Esso
\`e convenientemente utilizzabile per un corso avanzato, o di
dottorato, di Matematica e/o Informatica. In alternativa, i primi
capitoli possono essere utilizzati anche all'interno di corsi
pi\`u introduttivi di Calcolo Numerico o Analisi Numerica.

\vskip 7cm

\em Luigi Brugnano \`e nato a Brindisi il 16.12.1962. Si \`e
laureato in Scienze dell'Informazione nel 1985 presso
l'Universit\`a degli Studi di Bari. Nel 1990 \`e ricercatore di
Analisi Numerica presso la stessa Universit\`a. Nel 1992 \`e
professore associato di Analisi Numerica presso l'Universit\`a
degli Studi di Firenze e, dal 2001, professore ordinario di
Analisi Numerica presso la stessa Universit\`a.

\bigskip

Cecilia Magherini (22.6.1972--19.11.2020$\,^\dag$) si \`e laureata
in Scienze dell'Informazione nel 1999 presso l'Universit\`a degli
Studi di Firenze. Negli anni 1999-2000 \`e stata borsista presso
il CINECA, dove si \`e occupata di calcolo pa\-ral\-lelo. Ha quindi 
conseguito il Dottorato di Ricerca in Mate\-matica (XV ciclo),
presso l'Universit\`a degli Studi di Firenze, con una tesi
riguardante metodi e software numerici per equazioni
differenziali. Dal 2010 \`e stata prima ricercatore e poi, dal 2020, 
professore  associato di Analisi Numerica presso l'Universit\`a di Pisa.
\bigskip

\rule{12cm}{1mm}
}

\end{document}